\documentclass[11pt]{article}
\usepackage{amsfonts}
\usepackage{amssymb}
\usepackage{amsthm}
\usepackage{amsmath}
\usepackage{graphicx}
\usepackage{empheq}
\usepackage{indentfirst}
\usepackage{cite}
\usepackage{xcolor}
\usepackage{mathrsfs}
\usepackage{graphics}
\usepackage{mathabx}   
 \usepackage{esint}
 \allowdisplaybreaks[4]
\usepackage{latexsym}
\def\be{\begin{equation}}
\def\ee{\end{equation}}
\def\ba{\begin{array}}
\def\ea{\end{array}}
\def\bd{\begin{definition}}
\def\ed{\end{definition}}
\def\bt{\begin{theorem}}
\def\et{\end{theorem}}
\def\bc{\begin{corollary}}
\def\ec{\end{corollary}}
\def\bl{\begin{lemma}}
\def\el{\end{lemma}}
\def\bdm{\begin{displaymath}}
\def\edm{\end{displaymath}}

 \newtheorem{theorem}{Theorem}[section]
 \newtheorem{definition}{Definition}[section]
 \newtheorem{coro}{Corollary}[section]
 \newtheorem{lemma}{Lemma}[section]
 \newtheorem{prop}{Proposition}[section]
 \theoremstyle{definition}
 \newtheorem{remark}{Remark}[section]
 \numberwithin{equation}{section}

\newcommand{\beq}{\begin{equation}}
\newcommand{\eeq}{\end{equation}}

\makeatletter 
\def\Eqlfill@{\arrowfill@\Relbar\Relbar\Relbar} 
\newcommand{\extendEql}[1][]{\ext@arrow 0359\Eqlfill@{#1}} 
\makeatother

\usepackage[utf8]{inputenc}
\begin{document}
\bibliographystyle{plain}
\title{Weighted estimates for bilinear fractional integral operators and their commutators on Morrey spaces}
\author{Qianjun He\footnote{ {\it Email}: heqianjun16@mails.ucas.ac.cn\newline \hspace*{0.52cm} {\it Address}: School of Mathematics, University of Chinese Academy of Sciences, Beijing 100049, China}\qquad
Mingquan Wei\footnote{{\it Email}: weimingquan11@mails.ucas.ac.cn \newline \hspace*{0.52cm} {\it Address}: Xinyang Normal University, Xinyang, Henan 464000, China}\qquad
Dunyan Yan\footnote{{\it Email}: ydunyan@ucas.ac.cn  \newline \hspace*{0.52cm} {\it Address}: School of Mathematics, University of Chinese Academy of Sciences, Beijing 100049, China}
}
\date{}

\maketitle
\begin{abstract}
	This paper mainly dedicates to prove a plethora of weighted estimates on Morrey spaces for bilinear fractional integral operators and their general commutators with BMO functions of the form
	$$B_{\alpha}(f,g)(x)=\int_{\mathbb{R}^{n}}\frac{f(x-y)g(x+y)}{|y|^{n-\alpha}}dy,\qquad 0<\alpha<n.$$
	We also prove some maximal function control theorems for these operators, that is, the weighted Morrey norm is bounded by the weighted Morrey norm of a natural maximal operator when the weight belongs to $A_{\infty}$. As a corollary, some new weighted estimates for the bilinear maximal function associated to the bilinear Hilbert transform are obtained. Furthermore, we formulate a bilinear version of Stein-Weiss inequality on Morrey spaces for fractional integrals.
\end{abstract}
\medskip

 
 
\maketitle

\vspace*{0.5mm}
\section{Introduction}
\setcounter{equation}{0}
In the paper, we will consider the family of bilinear fractional integral operators that first was introduced by Kenig and Stein \cite{KS1999}, and its formula is as follows:
\begin{equation}
B_{\alpha}(f,g)(x):=\int_{\mathbb{R}^{n}}\frac{f(x-y)g(x+y)}{|y|^{n-\alpha}}dy,\qquad 0<\alpha<n.
\end{equation}
In fact, in 1992, Grafakos \cite{Grafakos1992} extended Adams's result to cover a multilinear analog of classical fractional integral operator. He studied a class of multilinear operators $I_{\alpha,k}$ (as general case for $B_{\alpha}$) which was defined by
$$I_{\alpha,k}(f_{1},\ldots,f_{k})(x)=\int_{\mathbb{R}^{n}}f_{1}(x-\theta_{1}y)\cdots f_{k}(x-\theta_{k}y)|y|^{\alpha-n}dy,$$
where $\theta_{j}\in\mathbb{R}\backslash 0$, $1\leq j\leq k$.  The following  two theorems are owning to Grafakos \cite{Grafakos1992}.

{\noindent}{\bf Theorem A}. Let $0<\alpha<n$, $1<p_{1}\ldots, p_{k}<\infty$ and $1/s=1/p_{1}+\cdots+1/p_{k}$, and $r$ be such that the identity $1/r+\alpha/n=1/s$ holds. Then $I_{\alpha,k}$ maps $L^{p_{1}}\times\cdots\times L^{p_{k}}$ into $L^{r}$ for $n/(n+\alpha)\leq s<n/\alpha$ (equivalently $1\leq r<\infty$).	

{\noindent}{\bf Theorem B}.	Let $p\in(1,\infty)$, ${1}/{p}=\sum_{j=1}^{k}1/p_{j}$, $p_{j}\in(1,\infty]$ and  $\alpha=n/p$. Assume the real numbers $\theta_{j}$ are distinct. Let $B$ be a ball in $\mathbb{R}^{n}$
and let $f_{j}\in L^{p_{j}}$ be supported in $B$. Then for any $\gamma<1$, there exists a constant $C_{0}(\gamma)$ depending only on $n$, $\alpha$, $\theta_{j}$ and $\gamma$ such that
\begin{equation}\label{endpoint estimate}
\frac{1}{|B|}\int_{B}\exp\left(\frac{n}{\omega_{n-1}}\gamma\left|\frac{LI_{\alpha,k}(f_{1}\ldots,f_{k})(x)}{\|f_{1}\|_{L^{p_{1}}}\cdots\|f_{k}\|_{L^{p_{k}}}}\right|^{p^{\prime}}\right)dx\leq C_{0}(\gamma),
\end{equation}
where $L=\prod_{j=1}^{k}|\theta_{j}|^{\frac{n}{p_{j}}}$. Furthermore, $\eqref{endpoint estimate}$ fails if $\gamma>1$.  In \cite{Grafakos1992}, Grafakos leaves an open problem for $\gamma=1$ in Theorem B. Later, Bak \cite{Bak1992} proved this open problem due to Grafakos \cite{Grafakos1992}. In \cite{GK2001,KS1999}, Grafakos-Kalton and Kenig-Stein
showed by very elementary considerations that $B_{\alpha}$ is bounded from $L^{p_{1}}\times L^{p_{2}}$ to $L^{r}$ for the full range $1<p_{1}, p_{2}<\infty$ and $1/r=1/p_{1}+1/p_{2}-\alpha/n$ with $r<\infty$ instead of the crucial condition $1\leq r<\infty$ for boundedness of $I_{\alpha,k}$. In \cite{KL2012}, Kuk and Lee proved that the endpoint weak boundedness of $I_{\alpha,k}$ in $L^{r,\infty}$ $(r=\frac{n}{2n-\alpha})$. To be more specific, they extended the index of Theorem A to $r>\frac{n}{2n-\alpha}$.
In \cite{Moen2014}, Moen studied  a plethora of new weighted estimates for $B_{\alpha}$ in Lebesgue spaces. In \cite{HM2016}, Hoang and Moen proved several weighted estimates for bilinear fractional integral operators and their commutators with BMO functions. Recently, the authors \cite{HY2018} obtained weighted estimates for $B_{\alpha}$ on Morrey spaces.

To make the present paper more readable, we announce some standard notations. For a measurable function $f$ the average of $f$ over a set $E$ is given by
$$m_{E}(f)=\fint_{E}fdx=\frac{1}{|E|}\int_{E}fdx.$$	
The Euclidian norm of a point $x=(x_{1},\ldots,x_{n})\in\mathbb{R}^{n}$ is given by $|x|=(x_{1}^{2}+\cdots+x_{n}^{2})^{1/2}$. We also use the $\ell^{\infty}$ norm $|x|_{\infty}=\max(|x_{1}|,\ldots,|x_{n}|)$. It is clear that $|x|_{\infty}\leq|x|\leq\sqrt{n}|x|_{\infty}$ for all  $x\in\mathbb{R}^{n}$. A cube with center $x_{0}$ and side length $d$, denoted by $Q=Q(x_{0},d)$, will be all points $x\in\mathbb{R}^{n}$ such that $|x-x_{0}|_{\infty}\leq\frac{d}{2}$. For an arbitrary cube $Q$, $c_{Q}$ will be its center and $\ell(Q)$ its side length, that is, $Q=Q(c_{Q},\ell(Q))$. Given $\lambda>0$ and a cube $Q$, we let $\lambda Q=Q(c_{Q},\lambda \ell(Q))$. The set of dyadic cubes, denoted $\mathcal{D}$, is all cubes of the form $2^{k}(m+[0,1)^{n})$ where $k\in\mathbb{Z}$ and $m\in\mathbb{Z}^{n}$. Finally, for $k\in\mathbb{Z}$ we let $\mathcal{D}_{k}$ denote the cubes of level $2^{k}$, that is, $\mathcal{D}_{k}=\{Q\in\mathcal{D}:\,\ell(Q)=2^{k}\}$. Throughout this paper, the letter $C$ will always denote a constant which may be different in each occasion, but is independent of the essential variables.

We first recall the definition of the Morrey (quasi-)norms \cite{Peetre1969}. For $0<q\leq p<\infty$, the Morrey norm is given by
\begin{equation}\label{define Morrey norm}
\|f\|_{\mathcal{M}_{q}^{p}}=\sup_{Q\in\mathcal{D}(\mathbb{R}^{n})}|Q|^{{1}/{p}}\left(\fint_{Q}|f(x)|^{q}dx\right)^{\frac{1}{q}}.
\end{equation}
The two-weight norm inequality for $B_{\alpha}$ on Morrey spaces were unknow untill the authors made some progress in \cite{HY2018}. This paper aim to study  weighted norm inequalities for $B_{\alpha}$ and its commutators on Morrey spaces. Given a linear operator $T$ and a function $b$, the commutator $[b,T]$ is defined by
$$[b,T]f=bT(f)-T(bf).$$
Coifman, Rochberg and Weiss \cite{CRW1976} introduced commutators of singular integral operators as a tool to extend the classical fractorization theory of $H^{p}$ spaces. They proved that if $b\in \rm BMO$ and $T$ is a singular operator, then $[b,T]$ in bounded on $L^{p}$ for $1<p<\infty$.

For a locally integrable function $b$, the $\rm BMO$ norm is defined by
$$\|b\|_{ \rm BMO}=\sup_{Q}\fint_{Q}|b(x)-m_{Q}(b)|dx<\infty,$$
where the supreme is taken from all cubes with sides parallel to the axises.
Let $\vec{b}=(b_{1},\cdots,b_{N})$ and $\vec{b}=(\beta_{1},\cdots,\beta_{N})\in\{1,2\}^{N}$, and then the iterated product commutators of $B_{\alpha}$ is given by
\begin{equation}\label{define general commutators}
[\vec{b},B_{\alpha}]_{\vec{\beta}}=[b_{N},[b_{N-1}\cdots,[b_{2},[b_{1},B_{\alpha}]_{\beta_{1}}]_{\beta_{2}}\cdots]_{\beta_{N-1}}]_{\beta_{N}}.
\end{equation}
It should be pointed out that in \cite{PRR2017} P\'{e}rez and Rivera-Rios studied these commutators in the linear case.
\section{The two-weight inequalities for bilinear fractinal integral operators and its commutators}

In the whole paper, we will work with two different cases. The first case is when $t\leq1$, while the second case is when $t>1$. Our departure is the following result inspired by the first and third author \cite{HY2018}.
\begin{theorem}\label{main}
	Let $v$ be a weight on $\mathbb{R}^{n}$ and $\vec{w}=(w_{1},w_{2})$ be a collection of two weights on $\mathbb{R}^{n}$. Assume that
	$$0<\alpha<n,\,\,\vec{q}=(q_{1},q_{2}),\,\, 1<q_{1},q_{2}<\infty,\,\,0<q\leq p<\infty,\,\,0<t\leq s<\infty,\,\,0<r\leq\infty$$
	and $1<a<\min(q_{1}, q_{2})$. Here, $q$ is given by ${1}/{q}={1}/{q_{1}}+{1}/{q_{2}}$. Suppose that
	$$\frac{\alpha}{n}>\frac{1}{r},\quad\frac{1}{s}=\frac{1}{p}+\frac{1}{r}-\frac{\alpha}{n},\quad\frac{t}{s}=\frac{q}{p},\quad 0<t\leq1$$
	and the weights ${v}$ and $\vec{w}$ satisfy the following two conditions:
	
	(i) If $0<s<1$,
	\begin{equation}\label{two weight condition_1}
	[v,\vec{w}]_{t,\vec{q}/{a}}^{r,as}:=\mathop{\sup_{Q,Q^{\prime}\in\mathcal{D}(\mathbb{R}^{n})}}_{Q\subset Q^{\prime}}\left(\frac{|Q|}{|Q^{\prime}|}\right)^{\frac{1-s}{as}}|Q^{\prime}|^{\frac{1}{r}}\left(\fint_{Q}v^{\frac{at}{1-t}}\right)^{\frac{1-t}{at}}\prod_{i=1}^{2}\left(\fint_{Q^{\prime}}w_{i}^{-(q_{i}/a)^{\prime}}\right)^{\frac{1}{(q_{i}/a)^{\prime}}}<\infty.
	\end{equation}
	
	(ii) If $s\geq1$,
	\begin{equation}\label{two weight condition_1_1}
	[v,\vec{w}]_{t,\vec{q}/{a}}^{r,as}:=\mathop{\sup_{Q,Q^{\prime}\in\mathcal{D}(\mathbb{R}^{n})}}_{Q\subset Q^{\prime}}\left(\frac{|Q|}{|Q^{\prime}|}\right)^{\frac{1-as}{as}}|Q^{\prime}|^{\frac{1}{r}}\left(\fint_{Q}v^{\frac{at}{1-t}}\right)^{\frac{1-t}{at}}\prod_{i=1}^{2}\left(\fint_{Q^{\prime}}w_{i}^{-(q_{i}/a)^{\prime}}\right)^{\frac{1}{(q_{i}/a)^{\prime}}}<\infty.
	\end{equation}
	In the above inequaity, we denote $\left(\fint_{Q}v^{\frac{at}{1-t}}\right)^{\frac{1-t}{at}}=\|v\|_{L^{\infty}(Q)}$ when $t=1$. Then we have
	$$\|B_{\alpha}(f,g)v\|_{\mathcal{M}_{t}^{s}}\leq C	[v,\vec{w}]_{t,\vec{q}/{a}}^{r,as}\sup_{Q\in\mathcal{D}(\mathbb{R}^{n})}|Q|^{1/p}\left(\fint_{Q}(|f|w_{1})^{q_{1}}\right)^{1/q_{1}}\left(\fint_{Q}(|g|w_{2})^{q_{2}}\right)^{1/q_{2}}.$$
\end{theorem}
\quad

By using a different technique, we are able to prove a similar result in the case $t>1$.\vspace{-0.2cm}
\begin{theorem}\label{main_1}
	Let $v$ be a weight on $\mathbb{R}^{n}$ and $\vec{w}=(w_{1},w_{2})$ be a collection of two weights on $\mathbb{R}^{n}$. Assume that
	$$0<\alpha<n,\, 0<q\leq p<\infty,\,1<t\leq s<r\leq\infty,\,\,\vec{q}=(q_{1},q_{2}),\, 1<r_{i}<q_{i}<\infty,\,i=1,2$$
	and $1<a<\min(q_{1},q_{2})$. Here, $r_{1}$, $r_{2}$ satisfy condition $1/r_{1}+1/r_{2}=1$ and $q$ is given by ${1}/{q}={1}/{q_{1}}+{1}/{q_{2}}$. Suppose that
	$$\frac{\alpha}{n}>\frac{1}{r},\quad\frac{1}{s}=\frac{1}{p}+\frac{1}{r}-\frac{\alpha}{n},\quad\frac{t}{s}=\frac{q}{p}$$
	and the weights ${v}$ and $\vec{w}$ satisfy the following condition
	\begin{equation}\label{two weight condition_2}
	[v,\vec{w}]_{t,\vec{q}/{a}}^{r,as}:=\mathop{\sup_{Q,Q^{\prime}\in\mathcal{D}(\mathbb{R}^{n})}}_{Q\subset Q^{\prime}}\left(\frac{|Q|}{|Q^{\prime}|}\right)^{\frac{1}{as}}|Q^{\prime}|^{\frac{1}{r}}\left(\fint_{Q}v^{at}\right)^{\frac{1}{at}}\prod_{i=1}^{2}\left(\fint_{Q^{\prime}}w_{i}^{-(q_{i}/a)^{\prime}}\right)^{\frac{1}{(q_{i}/a)^{\prime}}}<\infty.
	\end{equation}
	Then we have
	$$\|B_{\alpha}(f,g)v\|_{\mathcal{M}_{t}^{s}}\leq C	[v,\vec{w}]_{t,\vec{q}/{a}}^{r,as}\sup_{Q\in\mathcal{D}(\mathbb{R}^{n})}|Q|^{1/p}\left(\fint_{Q}(|f|w_{1})^{q_{1}}\right)^{1/q_{1}}\left(\fint_{Q}(|g|w_{2})^{q_{2}}\right)^{1/q_{2}}.$$	
\end{theorem}	
\quad

For the commutators for $B_{\alpha}$ define as $\eqref{define general commutators}$, we have the following theorems.
\vspace{-0.3cm}
\begin{theorem}\label{main_2}
	Let $v$ be a weight on $\mathbb{R}^{n}$ and $\vec{w}=(w_{1},w_{2})$ be a collection of two weights on $\mathbb{R}^{n}$. Assume that
	$$0<\alpha<n,\,\,\vec{q}=(q_{1},q_{2}),\,\, 1<q_{1},q_{2}<\infty,\,\,0<q\leq p<\infty,\,\,1/2\leq t\leq s<r\leq\infty$$
	and $1<a<\min(q_{1},q_{2})$. Here, $q$ is given by ${1}/{q}={1}/{q_{1}}+{1}/{q_{2}}$. Suppose that
	$$\frac{\alpha}{n}>\frac{1}{r},\quad\frac{1}{s}=\frac{1}{p}+\frac{1}{r}-\frac{\alpha}{n},\quad\frac{t}{s}=\frac{q}{p},\quad 0<t\leq1,\quad\vec{b}\in{ \rm BMO}^{N}$$
	and the weights ${v}$ and $\vec{w}$ satisfy the conditions $\eqref{two weight condition_1}$ and $\eqref{two weight condition_1_1}$. Then we have
	$$\|[\vec{b},B_{\alpha}]_{\vec{\beta}}(f,g)v\|_{\mathcal{M}_{t}^{s}}\leq C\|\vec{b}\|_{ {\rm BMO}^{N}}[v,\vec{w}]_{t,\vec{q}/{a}}^{r,as}\sup_{Q\in\mathcal{D}(\mathbb{R}^{n})}|Q|^{1/p}\left(\fint_{Q}(|f|w_{1})^{q_{1}}\right)^{1/q_{1}}\left(\fint_{Q}(|g|w_{2})^{q_{2}}\right)^{1/q_{2}}.$$
\end{theorem}

\begin{theorem}\label{main_3}
	Let $v$ be a weight on $\mathbb{R}^{n}$ and $\vec{w}=(w_{1},w_{2})$ be a collection of two weights on $\mathbb{R}^{n}$. Assume that
	$$0<\alpha<n,\, 0<q\leq p<\infty,\,1<t\leq s<r\leq\infty,\,\,\vec{q}=(q_{1},q_{2}),\, 1<r_{i}<q_{i}<\infty,\,i=1,2$$
	and $1<a<\min(q_{1},q_{2})$. Here, $r_{1}$, $r_{2}$ satisfy condition $1/r_{1}+1/r_{2}=1$ and $q$ is given by ${1}/{q}={1}/{q_{1}}+{1}/{q_{2}}$. Suppose that
	$$\frac{\alpha}{n}>\frac{1}{r},\quad\frac{1}{s}=\frac{1}{p}+\frac{1}{r}-\frac{\alpha}{n},\quad\frac{t}{s}=\frac{q}{p},\quad\vec{b}\in{ \rm BMO}^{N}$$
	and the weights ${v}$ and $\vec{w}$ satisfy the condition $\eqref{two weight condition_2}$. Then we have
	$$\|[\vec{b},B_{\alpha}]_{\vec{\beta}}(f,g)v\|_{\mathcal{M}_{t}^{s}}\leq C\|\vec{b}\|_{{\rm BMO}^{N}}[v,\vec{w}]_{t,\vec{q}/{a}}^{r,as}\sup_{Q\in\mathcal{D}(\mathbb{R}^{n})}|Q|^{1/p}\left(\fint_{Q}(|f|w_{1})^{q_{1}}\right)^{1/q_{1}}\left(\fint_{Q}(|g|w_{2})^{q_{2}}\right)^{1/q_{2}}.$$
\end{theorem}
\quad

We note here that the H\"{o}lder pairs $(r_{1},r_{2})$ in Theorems $\ref{main_1}$ and $\ref{main_3}$ exist and there are many such pairs. To name a few, we can start with $r_{1}=\frac{q_{1}}{q}$ and $r_{2}=\frac{q_{2}}{q}$, and then we use the facts that $r_{1}<q_{1}$ and $r_{2}<q_{2}$ to obtain more choices by considering either $r_{1}=\frac{q_{1}}{q}+\varepsilon$ or $r_{2}=\frac{q_{2}}{q}+\varepsilon$, for small $\varepsilon>0$.

While studying $B_{\alpha}$ and $[\vec{b},B_{\alpha}]_{\vec{\beta}}$, we need the following maximal operators which is defined by
\begin{equation}\label{define maximal operator}
M_{\alpha,\vec{R}}(f,g)(x)=\sup_{\mathcal{D}(\mathbb{R}^{n})\ni Q\ni x}|Q|^{\frac{\alpha}{n}}\left(\fint_{Q}|f(y)|^{r_{1}}dy\right)^{\frac{1}{r_{1}}}\left(\fint_{Q}|g(z)|^{r_{2}}dz\right)^{\frac{1}{r_{2}}},
\end{equation}
for a given vector $\vec{R}=(r_{1},r_{2})$.
When $\alpha=0$ we write $M_{\vec{R}}=M_{\alpha,\vec{R}}$ for short. The controls that we have mentioned above are stated in the theorems below.\vspace{-0.3cm}
\begin{theorem}\label{main_4}
	Assume that $0<\alpha<n$, $0<q<\infty$, $0<q\leq p<\infty$ and $(r_{1},r_{2})$ is a H\"{o}lder pair. If the weight $w\in A_{\infty}$, then
	$$\|B_{\alpha}(f,g)\|_{\mathcal{M}_{q}^{p}(w)}\leq C\|M_{\alpha,\vec{R}}(f,g)\|_{\mathcal{M}_{q}^{p}(w)},$$
	whenever the left-hand side is finite.
	
	Furthermore, if $w\in RH_{\nu}$ and $q\leq p\leq q\nu$, we have
	$$\|B_{\alpha}(f,g)\|_{\mathcal{M}_{q}^{p}(w)}\leq C\|M_{\alpha,\vec{R}}(f,g)\|_{\mathcal{M}_{q}^{p}(w)},$$
	whenever the left-hand side is finite.	
\end{theorem}

\begin{theorem}\label{main_5}
	Assume that $0<\alpha<n$, $0<q<\infty$, $0<q\leq p<\infty$, $\vartheta_{1},\vartheta_{2}>1$ and $(r_{1},r_{2})$ is a H\"{o}lder pair. If the weight $w\in A_{\infty}$, then
	$$\|[\vec{b},B_{\alpha}(f,g)]_{\vec{\beta}}\|_{\mathcal{M}_{q}^{p}(w)}\leq C\|\vec{b}\|_{{\rm BMO}^{N}}\|M_{\alpha,\vec{R_{\vartheta}}}(f,g)\|_{\mathcal{M}_{q}^{p}(w)},$$
	whenever the left-hand side is finite. Here, $\vec{R_{\vartheta}}$ is given by $\vec{R_{\vartheta}}=(\vartheta_{1}r_{1},\vartheta_{2}r_{2})$.
	
	Furthermore, if $w\in RH_{\nu}$ and $q\leq p\leq q\nu$, we have
	$$\|[\vec{b},B_{\alpha}(f,g)]_{\vec{\beta}}\|_{\mathcal{M}_{q}^{p}(w)}\leq C\|\vec{b}\|_{{\rm BMO}^{N}}\|M_{\alpha,\vec{R_{\vartheta}}}(f,g)\|_{\mathcal{M}_{q}^{p}(w)},$$
	whenever the left-hand side is finite. Here, $\vec{R_{\vartheta}}$ is given by $\vec{R_{\vartheta}}=(\vartheta_{1}r_{1},\vartheta_{2}r_{2})$.
\end{theorem}

Similar to $\eqref{two weight condition_2}$ in Theorem $\ref{main_1}$, we define another condition for the weights $(v,\vec{w})$ by
\begin{equation}\label{two weight condition_3}
[v,\vec{w}]_{t,\vec{q}}^{r,s}=\mathop{\sup_{Q,Q^{\prime}\in\mathcal{D}(\mathbb{R}^{n})}}_{Q\subseteq Q^{\prime}}\left(\frac{|Q|}{|Q^{\prime}|}\right)^{\frac{1}{s}}|Q^{\prime}|^{\frac{1}{r}}\left(\fint_{Q}v^{t}\right)^{\frac{1}{t}}\prod_{i=1}^{2}\left(\fint_{Q^{\prime}}w_{i}^{-r_{i}\left(\frac{q_{i}}{r_{i}}\right)^{\prime}}\right)^{\frac{1}{r_{i}\left(\frac{q_{i}}{r_{i}}\right)^{\prime}}}<\infty,
\end{equation}
where 
$$\left(\fint_{Q^{\prime}}w_{i}^{-r_{i}\left(\frac{q_{i}}{r_{i}}\right)^{\prime}}\right)^{\frac{1}{r_{i}\left(\frac{q_{i}}{r_{i}}\right)^{\prime}}}=\|w_{i}^{-1}\|_{L^{\infty}(Q^{\prime})}$$
 when $q_{i}=r_{i}$, $i=1,2$.

It turns out that condition $\eqref{two weight condition_3}$ can be characterized via the weak type and the strong type weighted boundedness of the maximal operator $M_{\alpha,\vec{R}}$, not only for $0 < \alpha< n$ but also for $\alpha= 0$. These results are stated in the following theorems.\vspace{-0.3cm}
\begin{theorem}\label{main_6}
	Let $v$ be a weight on $\mathbb{R}^{n}$ and $\vec{w}=(w_{1},w_{2})$ be a collection of two weights on $\mathbb{R}^{n}$. Assume that
	$$0\leq\alpha<n,\, 0<q\leq p<\infty,\,0<t\leq s<r\leq\infty,\,\,\vec{q}=(q_{1},q_{2})\,\,\text{and}\,\, 0<r_{i}\leq q_{i}<\infty,\,i=1,2.$$ Here, $q$ is given by ${1}/{q}={1}/{q_{1}}+{1}/{q_{2}}$. Suppose that
	$$\frac{\alpha}{n}\geq\frac{1}{r},\quad\frac{1}{s}=\frac{1}{p}+\frac{1}{r}-\frac{\alpha}{n}\quad\text{and}\quad\frac{1}{t}=\frac{1}{q}+\frac{1}{r}-\frac{\alpha}{n}.$$
	Then, for every $Q_{0}\in\mathcal{D}(\mathbb{R}^{n})$, the weighted inequality
	\begin{align}\label{the weak type boundedness}
	&\sup_{\lambda>0}|Q_{0}|^{\frac{1}{s}-\frac{1}{t}}\lambda (v^{t}(\{x\in Q_{0}:\,|M_{\alpha,\vec{R}}(f,g)(x)|>\lambda\}))^{\frac{1}{t}}\nonumber\\
	&\qquad\qquad\leq C[v,\vec{w}]_{t,\vec{q}}^{r,s} \mathop{\sup_{Q\in\mathcal{D}(\mathbb{R}^{n})}}_{Q\supseteq Q_{0}}|Q|^{1/p}\left(\fint_{Q}(|f|w_{1})^{q_{1}}\right)^{1/q_{1}}\left(\fint_{Q}(|g|w_{2})^{q_{2}}\right)^{1/q_{2}}
	\end{align}
	holds if and only if the weights ${v}$ and $\vec{w}$ satisfy the condition $\eqref{two weight condition_3}$.
\end{theorem}

The corresponding strong type result appears in following theorem, but we only obtain the sufficient conditions.
\begin{theorem}\label{main_7}
	Let $v$ be a weight on $\mathbb{R}^{n}$ and $\vec{w}=(w_{1},w_{2})$ be a collection of two weights on $\mathbb{R}^{n}$. Assume that
	$$0\leq\alpha<n,\, 0<q\leq p<\infty,\,0<t\leq s<r\leq\infty,\,\,\vec{q}=(q_{1},q_{2}),\,\, 0<r_{i}< q_{i}<\infty,\,i=1,2$$
	and $1<a<\min(q_{1}/r_{1},q_{2}/r_{2})$. Here, $q$ is given by ${1}/{q}={1}/{q_{1}}+{1}/{q_{2}}$. Suppose that
	$$\frac{\alpha}{n}\geq\frac{1}{r},\quad\frac{1}{s}=\frac{1}{p}+\frac{1}{r}-\frac{\alpha}{n},\quad\frac{t}{s}=\frac{q}{p}.$$
	and the weights ${v}$ and $\vec{w}$ satisfy the condition
	\begin{equation}\label{two weight condition_4}
	[v,\vec{w}]_{t,\vec{q}/a}^{r,s}=\mathop{\sup_{Q,Q^{\prime}\in\mathcal{D}(\mathbb{R}^{n})}}_{Q\subseteq Q^{\prime}}\left(\frac{|Q|}{|Q^{\prime}|}\right)^{\frac{1}{s}}|Q^{\prime}|^{\frac{1}{r}}\left(\fint_{Q}v^{t}\right)^{\frac{1}{t}}\prod_{i=1}^{2}\left(\fint_{Q^{\prime}}w_{i}^{-r_{i}\left(\frac{q_{i}}{ar_{i}}\right)^{\prime}}\right)^{\frac{1}{r_{i}\left(\frac{q_{i}}{ar_{i}}\right)^{\prime}}}<\infty.
	\end{equation}
	Then we have
	$$\|M_{\alpha,\vec{R}}(f,g)v\|_{\mathcal{M}_{t}^{s}}\leq C[v,\vec{w}]_{t,\vec{q}/a}^{r,s} \sup_{Q\in\mathcal{D}(\mathbb{R}^{n})}|Q|^{1/p}\left(\fint_{Q}(|f|w_{1})^{q_{1}}\right)^{1/q_{1}}\left(\fint_{Q}(|g|w_{2})^{q_{2}}\right)^{1/q_{2}}.$$
\end{theorem}

Taking $v=u_{1}^{\frac{1}{q_{1}}}u_{2}^{\frac{1}{q_{2}}}$ in Theorem $\ref{main_7}$,  we have the following theorem for two weights case.
\begin{theorem}\label{main_8}
	Let $u_{1}$ and $u_{2}$ are two weights on $\mathbb{R}^{n}$. Assume that 	
	$$0\leq\alpha<n,\, 0<q\leq p<\infty,\,0<t\leq s<\infty,\,\,\vec{q}=(q_{1},q_{2})\,\,\text{and}\,\, 0<r_{i}< q_{i}<\infty,\,i=1,2.$$
	Here, $q$ is given by ${1}/{q}={1}/{q_{1}}+{1}/{q_{2}}$. Suppose that
	$$\frac{1}{s}=\frac{1}{p}+\frac{1}{r}-\frac{\alpha}{n},\quad\frac{t}{s}=\frac{q}{p}$$
	and the weights $u_{1}$ and $u_{2}$ satisfy the condition
	\begin{align}\label{two weight condition_5}
	[u_{1},u_{2}]_{\vec{q}}^{s}:=\sup_{Q\in\mathcal{D}(\mathbb{R}^{n})}\left(\fint_{Q}u_{1}^{\frac{s}{q_{1}}}u_{2}^{\frac{s}{q_{2}}}\right)^{\frac{1}{s}}\left(\fint_{Q}u_{1}^{-\frac{r_{1}}{q_{1}-r_{1}}}\right)^{\frac{q_{1}-r_{1}}{r_{1}q_{1}}}\left(\fint_{Q}u_{2}^{-\frac{r_{2}}{q_{2}-r_{2}}}\right)^{\frac{q_{2}-r_{2}}{r_{2}q_{2}}}<\infty.
	\end{align}
	Then we have
	$$\|M_{\alpha,\vec{R}}(f,g)u_{1}^{\frac{1}{q_{1}}}u_{2}^{\frac{1}{q_{2}}}\|_{\mathcal{M}_{t}^{s}}\leq C[u_{1},u_{2}]_{\vec{q}}^{s} \sup_{Q\in\mathcal{D}(\mathbb{R}^{n})}|Q|^{1/p}\left(\fint_{Q}|f|^{q_{1}}u_{1}\right)^{1/q_{1}}\left(\fint_{Q}|g|^{q_{2}}u_{2}\right)^{1/q_{2}}.$$
\end{theorem}

Using Theorem $\ref{main_4}$ and Theorem $\ref{main_8}$, we can immediately obtain the following corollary.
\begin{coro}
	Let the indexs be as same as in Theorem $\ref{main_8}$, then the condition $\eqref{two weight condition_5}$ implies the following inequality
	$$\|B_{\alpha}(f,g)u_{1}^{\frac{1}{q_{1}}}u_{2}^{\frac{1}{q_{2}}}\|_{\mathcal{M}_{t}^{s}}\leq C[u_{1},u_{2}]_{\vec{q}}^{s} \sup_{Q\in\mathcal{D}(\mathbb{R}^{n})}|Q|^{1/p}\left(\fint_{Q}|f|^{q_{1}}u_{1}\right)^{1/q_{1}}\left(\fint_{Q}|g|^{q_{2}}u_{2}\right)^{1/q_{2}}.$$
\end{coro}

Finally, we end this section by giving an application of our estimates. The associated maximal operator to the bilinear Hilbert transform is defined as
$$BH(f,g)(x)=\sup_{r>0}\frac{1}{(2r)^{n}}\int_{[-r,r]^{n}}|f(x-y)g(x+y)|dy.$$
For the 1-dimensional case, Lacey \cite{Lacey2000} proved the boundedness of $BH$. He show that for $\frac{1}{p}=\frac{1}{p_{1}}+\frac{1}{p_{2}}$ and $p>3/2$, there holds
$$BH:\,L^{p_{1}}(\mathbb{R})\times L^{p_{2}}(\mathbb{R})\rightarrow L^{p}(\mathbb{R}).$$
Surprisingly, and contrary to usual paradigm in harmonic analysis, the boundedness of the bilinear Hilbert transform was shown first and used to proved the boundedness of Hilbert transform. For the $n$-dimensional case, Hoang and Moen \cite{HM2016} obtained new weighted estimates for this operator. By means of the inequality (see \cite{HM2016})
$$BH(f,g)(x)\leq{M}_{\vec{R}}(f,g)(x) \quad\text{for any}\quad{1}/{r_{1}}+{1}/{r_{2}}=1,$$
we have the following corollary.
\begin{coro}\label{corolary}
	Let $v$ be a weight on $\mathbb{R}^{n}$ and $\vec{w}=(w_{1},w_{2})$ be a collection of two weights on $\mathbb{R}^{n}$. Assume that
	$$ 0<q\leq p<\infty,\,0<t\leq s<\infty,\,\,\vec{q}=(q_{1},q_{2}),\,\, 1/r_{1}+1/r_{2}=1,\,\,1<r_{i}< q_{i}<\infty,\,i=1,2$$
	and $1<a<\min(q_{1}/r_{1},q_{2}/r_{2})$. Here, $q$ is given by ${1}/{q}={1}/{q_{1}}+{1}/{q_{2}}$. Suppose that
	$$\frac{1}{s}=\frac{1}{p}-\frac{\alpha}{n},\quad\frac{t}{s}=\frac{q}{p}.$$
	and the weights ${v}$ and $\vec{w}$ satisfies condition
	\begin{equation*}
	[v,\vec{w}]_{t,\vec{q}/a}^{s}=\mathop{\sup_{Q,Q^{\prime}\in\mathcal{D}(\mathbb{R}^{n})}}_{Q\subseteq Q^{\prime}}\left(\frac{|Q|}{|Q^{\prime}|}\right)^{\frac{1}{s}}\left(\fint_{Q}v^{t}\right)^{\frac{1}{t}}\prod_{i=1}^{2}\left(\fint_{Q^{\prime}}w_{i}^{-r_{i}\left(\frac{q_{i}}{ar_{i}}\right)^{\prime}}\right)^{\frac{1}{r_{i}\left(\frac{q_{i}}{ar_{i}}\right)^{\prime}}}<\infty.
	\end{equation*}
	Then we have
	$$\|BH(f,g)v\|_{\mathcal{M}_{t}^{s}}\leq C[v,\vec{w}]_{t,\vec{q}/a}^{s} \sup_{Q\in\mathcal{D}(\mathbb{R}^{n})}|Q|^{1/p}\left(\fint_{Q}(|f|w_{1})^{q_{1}}\right)^{1/q_{1}}\left(\fint_{Q}(|g|w_{2})^{q_{2}}\right)^{1/q_{2}}.$$
\end{coro}

Next, we consider a vector weight theorem. In this case, we choice $r_{1}=\frac{q_{1}}{q}$ and $r_{2}=\frac{q_{2}}{q}$.
\begin{coro}
	Let $u_{1}$ and $u_{2}$ are two weights on $\mathbb{R}^{n}$. Assume that 	
	$$0<q\leq p<\infty,\,0<t\leq s<\infty\,\,\text{and}\,\, 0<r_{i}< q_{i}<\infty,\,i=1,2.$$
	Here, $q$ is given by ${1}/{q}={1}/{q_{1}}+{1}/{q_{2}}$. Suppose that
	$$\frac{1}{s}=\frac{1}{p}+\frac{1}{r}-\frac{\alpha}{n},\quad\frac{t}{s}=\frac{q}{p}$$
	and the weights $u_{1}$ and $u_{2}$ satisfy the condition
	\begin{align}\label{two weight condition_9}
	[u_{1},u_{2}]_{{q}}^{s}:=\sup_{Q\in\mathcal{D}(\mathbb{R}^{n})}\left(\fint_{Q}u_{1}^{\frac{s}{q_{1}}}u_{2}^{\frac{s}{q_{2}}}\right)^{\frac{1}{s}}\left(\fint_{Q}u_{1}^{-\frac{1}{1-q}}\right)^{\frac{1-q}{q_{1}}}\left(\fint_{Q}u_{2}^{\frac{-1}{1-q}}\right)^{\frac{1-q}{q_{2}}}<\infty.
	\end{align}
	Then we have
	$$\|BH(f,g)u_{1}^{\frac{1}{q_{1}}}u_{2}^{\frac{1}{q_{2}}}\|_{\mathcal{M}_{t}^{s}}\leq C[u_{1},u_{2}]_{\vec{q}}^{s} \sup_{Q\in\mathcal{D}(\mathbb{R}^{n})}|Q|^{1/p}\left(\fint_{Q}|f|^{q_{1}}u_{1}\right)^{1/q_{1}}\left(\fint_{Q}|g|^{q_{2}}u_{2}\right)^{1/q_{2}}.$$
\end{coro}
\begin{remark}
	In \cite{HM2016}, Hoang and Moen shown that the condition $\eqref{two weight condition_9}$ is more general than classical weight condition $(u_{1},u_{2})\in A_{p}\times A_{p}$.
\end{remark}
In the forthcoming sections of this paper, we mainly give proofs of Theorems $\ref{main_2}$, $\ref{main_3}$, $\ref{main_5}$, $\ref{main_6}$, $\ref{main_7}$ and $\ref{main_8}$. Theorems $\ref{main_1}$ and $\ref{main_4}$ can be proved using similar techniques as those in the proof of Theorems $\ref{main_3}$ and $\ref{main_5}$, respectively.
At the end of this paper, we show an application of our results: a bilinear Stein-Weiss inequality.

\section{Preliminaries and some lemmas}

In this section, we present some preliminaries and lemmas for proving main results. We may rearrange the commutators in any order as the following Proposition states, which can be found in the paper \cite{HM2016}.

\begin{prop}\label{proposition_1}
	For any permutation $\sigma$ on $\{1,\cdots,N\}$,
	\begin{equation}\label{rearrange the commutators}
	[\sigma(\vec{b}),B_{\alpha}]_{\sigma(\vec{\beta})}=[\vec{b},B_{\alpha}]_{\vec{\beta}},
	\end{equation}
	where $\sigma(\vec{b})=(b_{\sigma(1)},\cdots,b_{\sigma(N)})$ and $\sigma(\vec{\beta})=(\beta_{\sigma(1)},\cdots,\beta_{\sigma(N)})$. In particular, equality $\eqref{rearrange the commutators}$ is valid for any permutation $\sigma_{0}$ be such that $\sigma_{0}(\vec{\beta})=(1,\cdots,1,2,\cdots,2)$.
\end{prop}

For simplicity in the notation and proof, from now on we will always assume that $$\vec{\beta}=(1,\cdots,1,2,\cdots,2)$$ and reserve the notation $m=m(\vec{\beta})$ to denote the number of 1's in $\vec{\beta}$.

We recall the following John-Nirenberg inequality (see \cite{Duoandikoetxea2001,Grafakos2008,LDY2007}).
\begin{lemma}\label{lemma_1}
	Let $1\leq p<\infty$. Assume that $Q$ be a cube and $b\in {\rm BMO}$. Then
	$$\left(\fint_{Q}|b(x)-m_{Q}(b)|^{p}dx\right)^{\frac{1}{p}}\leq C\|b\|_{\rm BMO}.$$
\end{lemma}
\quad

A dyadic grid $\mathcal{D}(\mathbb{R}^{n})$ is a countable collection of cubes  that satisfies the following properties:

\qquad(1) $Q\in\mathcal{D}(\mathbb{R}^{n})\Rightarrow\ell(Q)=2^{k}$ for some $k\in\mathbb{Z}$.

\qquad(2) For each $k\in\mathbb{Z}$, the set $\{Q\in\mathcal{D}(\mathbb{R}^{n}):\,\ell(Q)=2^{k}\}$ forms a partition of $\mathbb{R}^{n}$.

\qquad(3) $Q, R\in\mathcal{D}(\mathbb{R}^{n})\Rightarrow Q\cap R\in\{\emptyset, P, R\}$.

One very clear point for this concept is the dyadic grid that is formed by translating and then dilating the nuit cube $[0,1)^{n}$ all over $\mathbb{R}^{n}$. More precisely, it is formulated as
$$\mathcal{D}(\mathbb{R}^{n})=\{2^{-k}([0,1)^{n}+m):\,k\in\mathbb{Z}, m\in\mathbb{Z}^{n}\}.$$	
We invoke the following decomposition which is derived in \cite{Perez1994,Perez1995_1,Perez1995_2}.

Fix a cube $Q_{0}\in\mathcal{D}(\mathbb{R}^{n})$. Let $\mathcal{D}(Q_{0})$ be the collection of all dyadic subcubes of $Q_{0}$, that is, all those cubes obtained by dividing $Q_{0}$ into $2^{n}$ congruent cubes of half its side-length, dividing each of those into $2^{n}$ congruent cubes, and so on. By convention, $Q_{0}$ itself belongs to $\mathcal{D}(Q_{0})$.

\begin{lemma}\label{lemma_2}
	For $\theta_{1},\theta_{2}>1$, let $A=(4\cdot 18^{n})^{\frac{1}{\theta_{1}}+\frac{1}{\theta_{2}}}$ and
	$$\gamma:=\left(\fint_{3Q_{0}}|f(y)|^{\theta_{1}}dy\right)^{\frac{1}{\theta_{1}}}\left(\fint_{3Q_{0}}|g(z)|^{\theta_{2}}dz\right)^{\frac{1}{\theta_{2}}}.$$
	For $k=1,2,\cdots$, we take
	$$D_{k}:=\bigcup\left\{Q\in\mathcal{D}(Q_{0}):\,\left(\fint_{3Q}|f(y)|^{\theta_{1}}dy\right)^{\frac{1}{\theta_{1}}}\left(\fint_{3Q}|g(z)|^{\theta_{2}}dz\right)^{\frac{1}{\theta_{2}}}>\gamma A^{k}\right\}.$$
	Considering the maximality cube, we have
	$$D_{k}=\bigcup_{j}Q_{j}^{k}.$$
	Then we have
	$$\gamma A^{k}<\left(\fint_{3Q_{j}^{k}}|f(y)|^{\theta_{1}}dy\right)^{\frac{1}{\theta_{1}}}\left(\fint_{3Q_{j}^{k}}|g(z)|^{\theta_{2}}dz\right)^{\frac{1}{\theta_{2}}}\leq 2^{n(\frac{1}{\theta_{1}}+\frac{1}{\theta_{2}})}\gamma A^{k}.$$
	Let $E_{0}:=Q_{0}\backslash D_{1}$ and $E_{j}^{k}:=Q_{j}^{k}\backslash D_{k+1}$. Moreover we obtain
	$$|Q_{0}|\leq 2|E_{0}|\quad\text{and}\quad\left|Q_{j}^{k}\right|\leq 2\left|E_{j}^{k}\right|.$$
\end{lemma}
{\noindent}$Proof$. By the maximality of $Q_{j}^{k}$, we obtain the following:
\begin{equation}\label{the maximality of Q_{j}^{k}}
\gamma A^{k}<\left(\fint_{3Q_{j}^{k}}|f(y)|^{\theta_{1}}dy\right)^{\frac{1}{\theta_{1}}}\left(\fint_{3Q_{j}^{k}}|g(z)|^{\theta_{2}}dz\right)^{\frac{1}{\theta_{2}}}\leq 2^{n(\frac{1}{\theta_{1}}+\frac{1}{\theta_{2}})}\gamma A^{k}.
\end{equation}
Let $E_{0}=Q_{0}\backslash D_{1}$ and $E_{j}^{k}=Q_{j}^{k}\backslash D_{k+1}$. Then $\{E_{0}\}$ and $\{E_{j}^{k}\}$ are disjoint and satisfy $$E_{0}\bigcup\Big(\bigcup_{k,j}E_{j}^{k}\Big)=Q_{0}.$$
For any fix fixed $Q_{j}^{k}$, we set	
\begin{equation*}
A_{1}:=\left[\left(\int_{3Q_{j}^{k}}|f(y)|^{\theta_{1}}dy\right)^{\frac{1}{\theta_{1}}}\left(\int_{3Q_{j}^{k}}|g(z)|^{\theta_{2}}dz\right)^{\frac{1}{\theta_{2}}}\right]^{-\frac{\theta_{2}}{\theta_{1}+\theta_{2}}}(\gamma A^{k+1})^{\frac{\theta_{2}}{\theta_{1}+\theta_{2}}}\left(\int_{3Q_{j}^{k}}|f(y)|^{\theta_{1}}dy\right)^{\frac{1}{\theta_{1}}}
\end{equation*}	
and
\begin{equation*}
A_{2}:=\left[\left(\int_{3Q_{j}^{k}}|f(y)|^{\theta_{1}}dy\right)^{\frac{1}{\theta_{1}}}\left(\int_{3Q_{j}^{k}}|g(z)|^{\theta_{2}}dz\right)^{\frac{1}{\theta_{2}}}\right]^{-\frac{\theta_{1}}{\theta_{1}+\theta_{2}}}(\gamma A^{k+1})^{\frac{\theta_{1}}{\theta_{1}+\theta_{2}}}\left(\int_{3Q_{j}^{k}}|f(y)|^{\theta_{2}}dy\right)^{\frac{1}{\theta_{2}}}.
\end{equation*}	
Observe that $A_{1}A_{2}=\gamma A^{k+1}$. Define
$$M_{\theta_{1}}(f)(x):=\sup_{\mathcal{D}(\mathbb{R}^{n})\ni Q\ni x}\left(\fint_{Q}|f(y)|^{\theta_{1}}dy\right)^{\frac{1}{\theta_{1}}}\,\,\text{and}\,\,M_{\theta_{2}}(g)(x):=\sup_{\mathcal{D}(\mathbb{R}^{n})\ni Q\ni x}\left(\fint_{Q}|g(z)|^{\theta_{2}}dz\right)^{\frac{1}{\theta_{2}}}$$
and
$$M_{\theta_{1},\theta_{2}}(f,g)(x):=\sup_{\mathcal{D}(\mathbb{R}^{n})\ni Q\ni x}\left(\fint_{Q}|f(y)|^{\theta_{1}}dy\right)^{\frac{1}{\theta_{1}}}\left(\fint_{Q}|g(z)|^{\theta_{2}}dz\right)^{\frac{1}{\theta_{2}}}.$$
By $\eqref{the maximality of Q_{j}^{k}}$ we see that
\begin{align*}
Q_{j}^{k}\cap D_{k+1}&\subset\left\{x\in Q_{j}^{k}:\,M_{\theta_{1},\theta_{2}}(\chi_{Q_{j}^{k}}f,\chi_{Q_{j}^{k}}g)(x)>\gamma A^{k+1}\right\}\\
&\subset\left\{x\in Q_{j}^{k}:\,M_{\theta_{1}}(\chi_{Q_{j}^{k}}f)(x)M_{\theta_{2}}(\chi_{Q_{j}^{k}}g)(x)>\gamma A^{k+1}\right\}\\
&\subset\left\{x\in Q_{j}^{k}:\,M_{\theta_{1}}(\chi_{Q_{j}^{k}}f)(x)>A_{1}\right\}\bigcup\left\{x\in Q_{j}^{k}:\,M_{\theta_{2}}(\chi_{Q_{j}^{k}}g)(x)>A_{2}\right\}\\
&\subset\left\{x\in\mathbb{R}^{n}:\,M(\chi_{Q_{j}^{k}}f^{\theta_{1}})(x)>A_{1}^{\theta_{1}}\right\}\bigcup\left\{x\in\mathbb{R}^{n}:\,M(\chi_{Q_{j}^{k}}g^{\theta_{2}})(x)>A_{2}^{\theta_{2}}\right\}.
\end{align*}
Owning to the weak-$(1,1)$ boundedness of $M$, we have
\begin{align*}
\left|Q_{j}^{k}\cap D_{k+1}\right|&\leq\left|\left\{x\in\mathbb{R}^{n}:\,M(\chi_{Q_{j}^{k}}f^{\theta_{1}})(x)>A_{1}^{\theta_{1}}\right\}\right|+\left|\left\{x\in\mathbb{R}^{n}:\,M(\chi_{Q_{j}^{k}}g^{\theta_{2}})(x)>A_{2}^{\theta_{2}}\right\}\right|\\
&\leq\frac{3^{n}}{A_{1}^{\theta_{1}}}\int_{3Q_{j}^{k}}|f(y)|^{\theta_{1}}dy+\frac{3^{n}}{A_{2}^{\theta_{2}}}\int_{3Q_{j}^{k}}|f(z)|^{\theta_{2}}dz\\
&=2\cdot 3^{n}\left[\frac{1}{\gamma A^{k+1}}\left(\int_{3Q_{j}^{k}}|f(y)|^{\theta_{1}}dy\right)^{\frac{1}{\theta_{1}}}\left(\int_{3Q_{j}^{k}}|g(z)|^{\theta_{2}}dz\right)^{\frac{1}{\theta_{2}}}\right]^{\frac{\theta_{1}\theta_{2}}{\theta_{1}+\theta_{2}}},
\end{align*}
where we have used the defintions of $A_{1}$ and $A_{2}$. From $\eqref{the maximality of Q_{j}^{k}}$ we further have
\begin{align*}
\left|Q_{j}^{k}\cap D_{k+1}\right|&\leq2\cdot 3^{n}\left[\frac{1}{\gamma A^{k+1}}\left(\fint_{3Q_{j}^{k}}|f(y)|^{\theta_{1}}dy\right)^{\frac{1}{\theta_{1}}}\left(\fint_{3Q_{j}^{k}}|g(z)|^{\theta_{2}}dz\right)^{\frac{1}{\theta_{2}}}\right]^{\frac{\theta_{1}\theta_{2}}{\theta_{1}+\theta_{2}}}\left|3Q_{j}^{k}\right|\\
&\leq\frac{1}{2}\left|3Q_{j}^{k}\right|.
\end{align*}
Similarly, we can get
$$|Q_{0}|\leq 2|E_{0}|, \quad\left|Q_{0}^{\prime}\right|\leq2\left|E_{0}^{\prime}\right|\quad\text{and}\quad\left|{Q_{j}^{k}}^{\prime}\right|\leq2\left|{E_{j}^{k}}^{\prime}\right|.$$
This finishes the proof of Lemma $\ref{lemma_2}$. $\hfill$ $\Box$
\quad

The following two lemmas are essentially due to Iida \cite{Iida 2016}.
\begin{lemma}\label{lemma_3}
	Let $q_{1}$, $q_{2}$ and $a$ satisfy the conditions of Theorem $\ref{main_2}$, and then we can choose auxiliary indices $\theta_{1}$, $\theta_{2}$,
	$\theta_{3}$, $\theta_{4}$ and $\theta_{5}$ so that the following conditions hold:
	
	{\noindent}\quad{\rm 1.} $\theta_{1},\theta_{4}\in(1,q_{1})$, $\theta_{2},\theta_{5}\in(1,q_{2})$ and $\theta_{3}>1$.
	
	{\noindent}\quad{\rm 2.} For the indices $\theta_{1}\in(1,q_{1})$ and $\theta_{2}\in(1,q_{2})$, we can choose $a_{*}>1$ such that $$a_{*}\theta_{1}<q_{1}\quad\text{and}\quad a_{*}\theta_{2}<q_{2}.$$
	Assume in addition that, for these indices,
	$$a\geq\max\left\{\theta_{3},\frac{q_{1}}{(\theta_{1}(\frac{q_{1}}{a_{*}\theta_{1}})^{\prime})^{\prime}},\frac{q_{2}}{(\theta_{2}(\frac{q_{2}}{a_{*}\theta_{2}})^{\prime})^{\prime}},\frac{q_{1}}{(\theta_{4}(\frac{q_{1}}{\theta_{4}})^{\prime})^{\prime}},\frac{q_{2}}{(\theta_{5}(\frac{q_{2}}{\theta_{5}})^{\prime})^{\prime}}\right\}>1.$$
	Then we obtain
	$$\max\left\{\theta_{1}\left(\frac{q_{1}}{a_{*}\theta_{1}}\right)^{\prime},\theta_{4}\left(\frac{q_{1}}{\theta_{4}}\right)^{\prime}\right\}\leq\left(\frac{q_{1}}{a}\right)^{\prime}\quad\text{and}\quad
	\max\left\{\theta_{2}\left(\frac{q_{2}}{a_{*}\theta_{2}}\right)^{\prime},\theta_{5}\left(\frac{q_{2}}{\theta_{5}}\right)^{\prime}\right\}\leq\left(\frac{q_{2}}{a}\right)^{\prime}.$$
\end{lemma}
\quad

\begin{lemma}\label{lemma_4}
	Let $a>1$, $1<r_{1}<q_{1}$, $1<r_{2}<q_{2}$, $1/r_{1}+1/r_{2}=1$ and $1<t<r$, and then we can choose auxiliary indices $\vartheta_{1}$, $\vartheta_{2}$,
	$\vartheta_{3}$, $\vartheta_{4}$ and $\vartheta_{5}$ so that the following conditions hold:
	
	{\noindent}\quad{\rm 1.} $\vartheta_{1}r_{1},\vartheta_{4}r_{1}\in(r_{1},q_{1})$, $\vartheta_{2}r_{2},\vartheta_{5}r_{2}\in(r_{2},q_{2})$ and $\vartheta_{3}>1$.
	
	{\noindent}\quad{\rm 2.} Let $L>1$ and $e\in(t,r)$ such that $e\vartheta_{3}<Lt$ and $e^{\prime}\vartheta_{3}<t^{\prime}$.
	
	{\noindent}\quad{\rm 3.} For the indices $\vartheta_{1}r_{1}\in(r_{1},q_{1})$ and $\vartheta_{2}r_{2}\in(r_{2},q_{2})$, we can choose $a_{*}>1$ such that $$a_{*}\vartheta_{1}r_{1}<q_{1}\quad\text{and}\quad a_{*}\vartheta_{2}r_{2}<q_{2}.$$
	Assume in addition that, for these indices,
	$$a\geq\max\left\{\theta_{3},L,\frac{q_{1}}{(\vartheta_{1}r_{1}(\frac{q_{1}}{a_{*}\vartheta_{1}r_{1}})^{\prime})^{\prime}},\frac{q_{2}}{(\vartheta_{2}r_{2}(\frac{q_{2}}{a_{*}\vartheta_{2}r_{2}})^{\prime})^{\prime}},\frac{q_{1}}{(\vartheta_{4}r_{1}(\frac{q_{1}}{\vartheta_{4}r_{1}})^{\prime})^{\prime}},\frac{q_{2}}{(\vartheta_{5}r_{2}(\frac{q_{2}}{\vartheta_{5}r_{2}})^{\prime})^{\prime}}\right\}>1.$$
	Then we obtain
	$$\max\left\{\vartheta_{1}r_{1}\left(\frac{q_{1}}{a_{*}\vartheta_{1}r_{1}}\right)^{\prime},\vartheta_{4}r_{1}\left(\frac{q_{1}}{\vartheta_{4}r_{1}}\right)^{\prime}\right\}\leq\left(\frac{q_{1}}{a}\right)^{\prime}$$
	and
	$$
	\max\left\{\vartheta_{2}r_{2}\left(\frac{q_{2}}{a_{*}\vartheta_{2}r_{2}}\right)^{\prime},\vartheta_{5}r_{2}\left(\frac{q_{2}}{\vartheta_{5}r_{2}}\right)^{\prime}\right\}\leq\left(\frac{q_{2}}{a}\right)^{\prime}.$$
\end{lemma}

\quad

The following lemma is a simple property for BMO functions.
\begin{lemma}\label{lemma_5}
	Let $b\in \rm BMO$. Then for $Q_{0}, Q\in\mathcal{D}(\mathbb{R}^{n})$ and $Q\supsetneq Q_{0}$, these exists a natural number $k$ such that	
	$$|m_{Q_{0}}(b)-m_{Q}(b)|\leq k2^{n}\|b\|_{\rm BMO}.$$
\end{lemma}
{\noindent}Proof. Since $Q\supsetneq Q_{0}$, then there exists $k=1,2,\cdots$ such that
$$Q_{k}:=Q,\,\,Q_{j}\in\mathcal{D}(\mathbb{R}^{n}),\,\,Q_{j}\supsetneq Q_{j-1}\,\,\text{and}\,\,\,|Q_{j}|=2^{n}|Q_{j-1}|,\,\,j=1,2,\cdots,k.$$
By triangle inequality, we obtain
\begin{align*}
|m_{Q_{0}}(b)-m_{Q}(b)|&=\Big|\sum_{j=1}^{k}(m_{Q_{j-1}}(b)-m_{Q_{j}}(b))\Big|\leq\sum_{j=1}^{k}|m_{Q_{j-1}}(b)-m_{Q_{j}}(b)|.
\end{align*}
Moreover, we have
\begin{align*}
|m_{Q_{0}}(b)-m_{Q}(b)|&\leq\sum_{j=1}^{k}\left|\fint_{Q_{j-1}}b(y)dy-m_{Q_{j}}(b)\right|\\
&\leq\sum_{j=1}^{k}\fint_{Q_{j-1}}|b(y)-m_{Q_{j}}(b)|dy\\
&\leq\sum_{j=1}^{k}\frac{2^{n}}{|Q_{j}|}\int_{Q_{j}}|b(y)-m_{Q_{j}}|dy\\
&\leq k2^{n}\|b\|_{\rm BMO}.
\end{align*}
This completes the proof of Lemma $\ref{lemma_5}$. $\hfill$ $\Box$

To prove the results in Section 2, we also need the following lemma.\vspace{-0.3cm}
\begin{lemma}\label{lemma_6}
	Let $0\leq\alpha<2n$ and $\vec{R}=(r_{1},r_{2})$. Assume that $0<r_{i}<q_{i}<\infty$, $i=1,2$. If $0<t\leq s<\infty$ and $0<q\leq p<\infty$ satisfy
	\begin{equation}\label{condition of maximal operator on Morrey space }
	\frac{1}{s}=\frac{1}{p}-\frac{\alpha}{n}\quad\text{and}\quad\frac{t}{s}=\frac{q}{p},
	\end{equation}
	where $q$ is given by $1/q=1/q_{1}+1/q_{2}$, then
	$$\|M_{\alpha,\vec{R}}(f,g)\|_{\mathcal{M}_{t}^{s}}\leq C\sup_{Q\in\mathcal{D}(\mathbb{R}^{n})}|Q|^{\frac{1}{p}}\left(\fint_{Q}|f(y)|^{q_{1}}dy\right)^{\frac{1}{q_{1}}}\left(\fint_{Q}|g(z)|^{q_{2}}dz\right)^{\frac{1}{q_{2}}}.$$
\end{lemma}
\quad

Finally we define $A_{\infty}$ as the union of all $A_{p}$ classes for $p>1$. From \cite{Duoandikoetxea2001}, we know the following facts.\vspace{-0.3cm}
\begin{lemma}\label{lemma_7}
	If $w\in A_{\infty}$, then the following hold:
	\begin{enumerate}
		\item[\rm i)] For every $\eta\in(0,1)$, there exists $\kappa\in(0,1)$ such that: given a cube $Q$ and $E\subseteq Q$ with $|E|\leq\eta|Q|$, we will also have $w(E)\leq\kappa w(Q)$.
		\item[\rm ii)] There exist an $\vartheta_{3}>1$ such that
		$$\left(\fint_{Q}w(x)^{\vartheta_{3}}dx\right)^{\frac{1}{\vartheta_{3}}}\leq C\fint_{Q}w(x)dx.$$
	\end{enumerate}
\end{lemma}
\quad

Let a nonnegative locally integrable function $w$ on $\mathbb{R}^{n}$ belong to the reverse H\"{o}lder class $RH_{\nu}$ for $1<\nu<\infty$ if it satisfies the reverse H\"{o}lder inequality with exponent $\nu$, i.e.
$$\left(\fint_{Q}w(x)^{\nu}dx\right)^{\frac{1}{\nu}}\leq C\fint_{Q}w(x)dx,$$
where the constant $C$ is universal for all cubes $Q\subset\mathbb{R}^{n}$.

Duoandikoetxea and Rosenthal \cite{DR2017} proved an extrapolation theorem for $A_{\infty}$ weights. Namely,\vspace{-0.3cm}
\begin{lemma}\label{lemma_8}
	Let $0<q_{0}<\infty$ and let $\mathcal{F}$ be a collection of nonnegative measurable pairs of functions. Assume that for every $(f,g)\in\mathcal{F}$ every $w\in A_{\infty}$, we have
	$$\int_{\mathbb{R}^{n}}|f|^{q_{0}}w\leq C\int_{\mathbb{R}^{n}}|g|^{q_{0}}w,$$
	whenever the left-hand side is finite. Then for every $0<q<\infty$, every $0<q\leq p<\infty$ and every $w\in A_{\infty}$, we have
	$$\|f\|_{\mathcal{M}_{q}^{p}(w)}\leq C\|g\|_{\mathcal{M}_{q}^{p}(w)},$$
	whenever the left-hand side is finite.
	
	Furthermore, if $w\in RH_{\nu}$ and $q\leq p\leq q\nu$, we have
	$$\|g\|_{\mathcal{M}_{q}^{p}(w)}\leq C\|f\|_{\mathcal{M}_{q}^{p}(w)},$$
	whenever the left-hand side is finite.
\end{lemma}

We also need the following a characterization of a multiple weights given by Iida \cite{Iida2012}.
\begin{lemma}\label{Characterization of a multiple weights}
	Let $1<q_{1},q_{2}<\infty$ and $\hat{t}\geq q$ with $1/q=1/q_{1}+1/q_{2}$. Then, for two weights $w_{1}$, $w_{2}$, the inequality
	$$\sup_{Q\in\mathcal{D}(\mathbb{R}^{n})}\left(\fint_{Q}(w_{1}w_{2})^{\hat{t}}\right)^{1/\hat{t}}\prod_{i=1}^{2}\left(\fint_{Q}w_{i}^{-q_{i}^{\prime}}\right)^{1/q_{i}^{\prime}}<\infty$$
	holds if and only if
	$$
	\left\{
	\begin{aligned}
	&(w_{1}w_{2})^{\hat{t}}\in A_{1+\hat{t}(2-1/q)}, \\
	&w_{i}^{-q_{i}^{\prime}}\in A_{q_{i}^{\prime}(1/\hat{t}+2-1/q)},\,\,i=1,2.
	\end{aligned}
	\right.
	$$
\end{lemma}

\section{The proof of Theorem $\ref{main_2}$ }

Without loss of generality, we may assume that $f$ and $g$ are nonnegative, bounded and compactly supported. By induction, it is obvious that
\begin{equation}
[\vec{b},B_{\alpha}]_{\vec{\beta}}(f,g)(x)=\int_{\mathbb{R}^{n}}\prod_{i=1}^{m}(b_{i}(x)-b_{i}(x-y))\prod_{i=m+1}^{N}(b_{i}(x)-b_{i}(x+y))\frac{f(x-y)g(x+y)}{|y|^{n-\alpha}}dy.
\end{equation}

Fix  a dyadic cube $Q\in\mathcal{D}(\mathbb{R}^{n})$. Let $\mathcal{D}_{\nu}$ be the collection of dyadic cubes and the volume of the elements of $\mathcal{D}_{\nu}$ is $2^{n\nu}$. For $x\in Q_{0}$, we have
\begin{align}\label{decomposition of commutator for bilinear fractional integral}
|[\vec{b},B_{\alpha}]_{\vec{\beta}}(f,g)(x)|&\leq C\sum_{\nu\in\mathbb{Z}}\mathop{\sum_{Q\in\mathcal{D}(\mathbb{R}^{n})}}_{|Q|=2^{n\nu}}2^{-\nu(n-\alpha)}\chi_{Q}(x)\mathop{\int}_{|y|\leq \ell(Q)}\prod_{i=1}^{m}|b_{i}(x)-b_{i}(x-y)|\nonumber\\
&\qquad\qquad\qquad\times\prod_{i=m+1}^{N}|b_{i}(x)-b_{i}(x+y)|{f(x-y)g(x+y)}dy\nonumber\\
&=C\sum_{\nu\in\mathbb{Z}}\Big(\mathop{\sum_{Q\in\mathcal{D}_{\nu}}}_{Q_{0}\subseteq Q}+\mathop{\sum_{Q\in\mathcal{D}_{\nu}}}_{Q\supsetneq Q_{0}}\Big)2^{-\nu(n-\alpha)}\chi_{Q}(x)\mathop{\int}_{|y|\leq \ell(Q)}\prod_{i=1}^{m}|b_{i}(x)-b_{i}(x-y)|\nonumber\\
&\qquad\times\prod_{i=m+1}^{N}|b_{i}(x)-b_{i}(x+y)|{f(x-y)g(x+y)}dy=:C(\mathscr{T}_{1}(x)+\mathscr{T}_{2}(x)).
\end{align}
For each $Q\in\mathcal{D}(\mathbb{R}^{n})$, let
\begin{equation}\label{define as average of BMO function}
\lambda_{i}=\lambda_{i}(Q)=\fint_{3Q}b_{i}(x)dx,
\end{equation}
where $i=1,\cdots,N$, and then we have
\begin{align}\label{part one}
\prod_{i=1}^{m}(b_{i}(x)-b_{i}(x-y))&=\prod_{i=1}^{m}\left[(b_{i}(x)-\lambda_{i})+(\lambda_{i}-b_{i}(x-y))\right]\nonumber\\
&=\sum_{A\subseteq\{1,\cdots,m\}}\prod_{i\in A}(b_{i}(x)-\lambda_{i})\prod_{i\in\bar{A}}(\lambda_{i}-b_{i}(x-y)).
\end{align}
Similarly, there holds
\begin{align}\label{part two}
\prod_{i=m+1}^{N}(b_{i}(x)-b_{i}(x+y))=\sum_{B\subseteq\{m+1,\cdots,N\}}\prod_{i\in B}(b_{i}(x)-\lambda_{i})\prod_{i\in\bar{B}}(\lambda_{i}-b_{i}(x+y)).
\end{align}
\quad

Now, we start to estimate $\mathscr{T}_{1}(x)$ and $\mathscr{T}_{2}(x)$.

{\noindent}{\bf Estimate of $\mathscr{T}_{1}$.} By the defintion of $\mathscr{T}_{1}$, trangle inequality, the equalities $\eqref{define as average of BMO function}$, $\eqref{part one}$ and $\eqref{part two}$, we obtain
\begin{align}\label{estimate for T_{1}}
\mathscr{T}_{1}(x)&\leq \sum_{A\subseteq\{1,\cdots,m\}}\sum_{B\subseteq\{m+1,\cdots,N\}}\sum_{Q\in\mathcal{D}(Q_{0})}|Q|^{\frac{\alpha}{n}-1}\mathop{\int}_{|y|_{\infty}\leq \ell(Q)}\prod_{i\in A\cup B}|b_{i}(x)-\lambda_{i}|\prod_{i\in\bar{A}}|b_{i}(x-y)-\lambda_{i}|\nonumber\\
&\quad\times\prod_{i\in\bar{B}}|b_{i}(x+y)-\lambda_{i}|f(x-y)g(x+y)dy\chi_{Q}(x).
\end{align}
Since $t\leq1$, we have
\begin{align*}
\int_{Q_{0}}(|\mathscr{T}_{1}(x)|v(x))^{t}dx&\lesssim\sum_{A\subseteq\{1,\cdots,m\}}\sum_{B\subseteq\{m+1,\cdots,N\}}\sum_{Q\in\mathcal{D}(Q_{0})}|Q|^{(\frac{\alpha}{n}-1)t}\int_{Q}\left[\int_{|y|_{\infty}\leq \ell(Q)}\prod_{i\in\bar{A}}|b_{i}(x-y)-\lambda_{i}|\right.\\
&\quad\times\left.\prod_{i\in\bar{B}}|b_{i}(x+y)-\lambda_{i}|f(x-y)g(x+y)dy\right]^{t}\left[\prod_{i\in A\cup B}|b_{i}(x)-\lambda_{i}|\right]^{t}v(x)^{t}dx.
\end{align*}
If we use H\"{o}lder's inequality with the pair $(\frac{1}{t},\frac{1}{1-t})$, we will arrive at the inequality
\begin{align*}
&\int_{Q_{0}}(|\mathscr{T}_{1}(x)|v(x))^{t}dx\lesssim\sum_{A\subseteq\{1,\cdots,m\}}\sum_{B\subseteq\{m+1,\cdots,N\}}\sum_{Q\in\mathcal{D}(Q_{0})}|Q|^{(\frac{\alpha}{n}-1)t}\left[\int_{Q}\int_{|y|_{\infty}\leq \ell(Q)}\prod_{i\in\bar{A}}|b_{i}(x-y)-\lambda_{i}|\right.\\
&\qquad\qquad\quad\times\left.\prod_{i\in\bar{B}}|b_{i}(x+y)-\lambda_{i}|f(x-y)g(x+y)dydx\right]^{t}\left[\int_{Q}\prod_{i\in A\cup B}|b_{i}(x)-\lambda_{i}|^{\frac{t}{1-t}}v(x)^{\frac{t}{1-t}}dx\right]^{1-t}.
\end{align*}
By a change of variables, we obtain
\begin{align}\label{change of variables}
&\int_{Q_{0}}(|\mathscr{T}_{1}(x)|v(x))^{t}dx\leq C\sum_{A\subseteq\{1,\cdots,m\}}\sum_{B\subseteq\{m+1,\cdots,N\}}\sum_{Q\in\mathcal{D}(Q_{0})}|Q|^{\frac{\alpha}{n}t+1}\left[\fint_{3Q}\prod_{i\in\bar{A}}|b_{i}(y)-\lambda_{i}|f(y)dy\right.\nonumber\\
&\qquad\qquad\qquad\qquad\times\left.\fint_{3Q}\prod_{i\in\bar{B}}|b_{i}(z)-\lambda_{i}|g(z)dz\right]^{t}\left[\fint_{Q}\prod_{i\in A\cup B}|b_{i}(x)-\lambda_{i}|^{\frac{t}{1-t}}v(x)^{\frac{t}{1-t}}dx\right]^{1-t}.
\end{align}
Now by using H\"{o}lder's inequality for $\theta_{1}>1$ and Lemma $\ref{lemma_1}$, we obtain the following estimates:
\begin{align*}
II_{0}:=\fint_{3Q}\prod_{i\in\bar{A}}|b_{i}(y)-\lambda_{i}|f(y)dy&\leq\left(\fint_{3Q}\prod_{i\in\bar{A}}|b_{i}(y)-\lambda_{i}|^{\theta_{1}^{\prime}}dy\right)^{\frac{1}{\theta_{1}^{\prime}}}\left(\fint_{3Q}|f(y)|^{\theta_{1}}dy\right)^{\frac{1}{\theta_{1}}}\\
&\leq\prod_{i\in\bar{A}}\left(\fint_{3Q}|b_{i}(y)-\lambda_{i}|^{|\bar{A}|\theta_{1}^{\prime}}dy\right)^{\frac{1}{|\bar{A}|\theta_{1}^{\prime}}}\left(\fint_{3Q}|f(y)|^{\theta_{1}}dy\right)^{\frac{1}{\theta_{1}}}\\
&\leq C\prod_{i\in\bar{A}}\|b_{i}\|_{\rm BMO}\left(\fint_{3Q}|f(y)|^{\theta_{1}}dy\right)^{\frac{1}{\theta_{1}}}.
\end{align*}
For $\theta_{2}>1$, similar to the estimate for $II_{0}$, we have
$$\fint_{3Q}\prod_{i\in\bar{B}}|b_{i}(z)-\lambda_{i}|g(z)dz\leq C\prod_{i\in\bar{B}}\|b_{i}\|_{\rm BMO}\left(\fint_{3Q}|g(z)|^{\theta_{2}}dz\right)^{\frac{1}{\theta_{2}}}.$$
Then we obtain
\begin{align*}
\fint_{Q}\prod_{i\in A\cup B}|b_{i}(x)-\lambda_{i}|^{\frac{t}{1-t}}v(x)^{\frac{t}{1-t}}dx&\leq\left(\fint_{Q}\prod_{i\in A\cup B}|b_{i}(x)-\lambda_{i}|^{\frac{\theta_{3}^{\prime}t}{1-t}}dx\right)^{\frac{1}{\theta_{3}^{\prime}}}\left(\fint_{Q}v(x)^{\frac{\theta_{3}t}{1-t}}dx\right)^{\frac{1}{\theta_{3}}}\\
&\leq\prod_{i\in A\cup B}\|b_{i}\|_{\rm BMO}^{\frac{t}{1-t}}\left(\fint_{Q}v(x)^{\frac{\theta_{3}t}{1-t}}dx\right)^{\frac{1}{\theta_{3}}},
\end{align*}
where this estimate is based on the facts $\theta_{3}>1$ and $t\geq\frac{1}{2}$.

Substituting these estimates into $\eqref{change of variables}$, we come up with the following estimates:
\begin{align}\label{final estimate_1}
\int_{Q_{0}}(|\mathscr{T}_{1}(x)|v(x))^{t}dx&\leq C\|\vec{b}\|_{{\rm BMO}^{N}}^{t}\sum_{A\subseteq\{1,\cdots,m\}}\sum_{B\subseteq\{m+1,\cdots,N\}}\sum_{Q\in\mathcal{D}(Q_{0})}|Q|^{\frac{\alpha}{n}t+1}\nonumber\\
&\qquad\times\left(\fint_{3Q}|f(y)|^{\theta_{1}}dy\right)^{\frac{t}{\theta_{1}}}\left(\fint_{3Q}|g(z)|^{\theta_{2}}dz\right)^{\frac{t}{\theta_{2}}}\left(\fint_{Q}v(x)^{\frac{\theta_{3}t}{1-t}}dx\right)^{\frac{1-t}{\theta_{3}}}.
\end{align}
From here, it suffices to estimate inner most sum of the last expression in $\eqref{final estimate_1}$ for a general dyadic grid $\mathcal{D}(\mathbb{R}^{n})$. We denote this sum as
\begin{equation*}
\mathcal{S}:=\sum_{Q\in\mathcal{D}(Q_{0})}|Q|^{\frac{\alpha}{n}t+1}\left(\fint_{3Q}|f(y)|^{\theta_{1}}dy\right)^{\frac{t}{\theta_{1}}}\left(\fint_{3Q}|g(z)|^{\theta_{2}}dz\right)^{\frac{t}{\theta_{2}}}\left(\fint_{Q}v(x)^{\frac{\theta_{3}t}{1-t}}dx\right)^{\frac{1-t}{\theta_{3}}}.
\end{equation*}
Nxet we will estimate for $\mathcal{S}$. Let
$$\mathcal{D}_{0}(Q_{0}):=\left\{Q\in\mathcal{D}(Q_{0}):\,\left(\fint_{3Q}|f(y)|^{\theta_{1}}dy\right)^{\frac{1}{\theta_{1}}}\left(\fint_{3Q}|g(z)|^{\theta_{2}}dz\right)^{\frac{1}{\theta_{2}}}\leq \gamma A\right\}$$
and
$$\mathcal{D}_{j}^{k}(Q_{0}):=\left\{Q\in\mathcal{D}(Q_{0}):\,Q\subseteq Q_{j}^{k},\gamma A^{k}<\left(\fint_{3Q}|f(y)|^{\theta_{1}}dy\right)^{\frac{1}{\theta_{1}}}\left(\fint_{3Q}|g(z)|^{\theta_{2}}dz\right)^{\frac{1}{\theta_{2}}}\leq \gamma A^{k+1}\right\},$$
where $Q_{j}^{k}$ is as in Lemma $\ref{lemma_2}$. Then we have
$$\mathcal{D}(Q_{0})=\mathcal{D}_{0}(Q_{0})\bigcup\Big(\bigcup_{k,j}\mathcal{D}_{j}^{k}(Q_{0})\Big).$$
From the definition of $\mathcal{S}$, we obtain
\begin{align*}
\mathcal{S}&=\Big(\sum_{Q\in\mathcal{D}(Q_{0})}+\sum_{k,j}\sum_{Q\in\mathcal{D}_{j}^{k}(Q_{0})}\Big)|Q|^{\frac{\alpha}{n}t+1}\left(\fint_{3Q}|f(y)|^{\theta_{1}}dy\right)^{\frac{t}{\theta_{1}}}\left(\fint_{3Q}|g(z)|^{\theta_{2}}dz\right)^{\frac{t}{\theta_{2}}}\left(\fint_{Q}v(x)^{\frac{\theta_{3}t}{1-t}}dx\right)^{\frac{1-t}{\theta_{3}}}\\
&=:\mathcal{S}_{0}+\sum_{k,j}\mathcal{S}_{j}^{k}.
\end{align*}
Now we begin to evaluate $\mathcal{S}_{j}^{k}$. If $Q\in\mathcal{D}_{j}^{k}(Q_{0})$, we have
$$\left(\fint_{3Q}|f(y)|^{\theta_{1}}dy\right)^{\frac{1}{\theta_{1}}}\left(\fint_{3Q}|g(z)|^{\theta_{2}}dz\right)^{\frac{1}{\theta_{2}}}\leq \gamma A^{k+1}.$$
Therefore, we obtain
\begin{align*}
\mathcal{S}_{j}^{k}\leq(\gamma A^{k+1})^{t}\sum_{Q\in\mathcal{D}_{j}^{k}(Q_{0})}|Q|^{\frac{\alpha}{n}t+1}\left(\fint_{Q}v(x)^{\frac{\theta_{3}t}{1-t}}dx\right)^{\frac{1-t}{\theta_{3}}}.
\end{align*}
Since
$$\sum_{Q\in\mathcal{D}_{j}^{k}(Q_{0})}=\mathop{\sum_{Q\in\mathcal{D}_{j}^{k}(Q_{0})}}_{Q\subseteq Q_{j}^{k}}$$
and
\begin{align*}
&\mathop{\sum_{Q\in\mathcal{D}_{j}^{k}(Q_{0})}}_{Q\subseteq Q_{j}^{k}}|Q|^{\frac{\alpha}{n}t+1}\left(\fint_{Q}v(x)^{\frac{\theta_{3}t}{1-t}}dx\right)^{\frac{1-t}{\theta_{3}}}\\
&=\sum_{i=0}^{\infty}\mathop{\sum_{Q\in\mathcal{D}_{j}^{k}(Q_{0}),\,Q\subseteq Q_{j}^{k}}}_{\ell(Q)=2^{-i}\ell(Q_{j}^{k})}|Q|^{\frac{\alpha}{n}t+1-\frac{1-t}{\theta_{3}}}\left(\int_{Q}v(x)^{\frac{\theta_{3}t}{1-t}}dx\right)^{\frac{1-t}{\theta_{3}}}\\
&=|Q_{j}^{k}|^{\frac{\alpha}{n}t+1-\frac{1-t}{\theta_{3}}}\sum_{i=0}^{\infty}2^{-\frac{\alpha}{n}t-1+\frac{1-t}{\theta_{3}}}\mathop{\sum_{Q\in\mathcal{D}_{j}^{k}(Q_{0}),\,Q\subseteq Q_{j}^{k}}}_{\ell(Q)=2^{-i}\ell(Q_{j}^{k})}\left(\int_{Q}v(x)^{\frac{\theta_{3}t}{1-t}}dx\right)^{\frac{1-t}{\theta_{3}}}\\
&\leq|Q_{j}^{k}|^{\frac{\alpha}{n}t+1-\frac{1-t}{\theta_{3}}}\sum_{i=0}^{\infty}2^{-\frac{\alpha}{n}t-1+\frac{1-t}{\theta_{3}}}\left[\mathop{\sum_{Q\in\mathcal{D}_{j}^{k}(Q_{0}),\,Q\subseteq Q_{j}^{k}}}_{\ell(Q)=2^{-i}\ell(Q_{j}^{k})}\int_{Q}v(x)^{\frac{\theta_{3}t}{1-t}}dx\right]^{\frac{1-t}{\theta_{3}}}\left[\mathop{\sum_{Q\in\mathcal{D}_{j}^{k}(Q_{0}),\,Q\subseteq Q_{j}^{k}}}_{\ell(Q)=2^{-i}\ell(Q_{j}^{k})}1\right]^{1-\frac{1-t}{\theta_{3}}}\\
&=|Q_{j}^{k}|^{\frac{\alpha}{n}t+1-\frac{1-t}{\theta_{3}}}\left(\int_{Q_{j}^{k}}v(x)^{\frac{\theta_{3}t}{1-t}}dx\right)^{\frac{1-t}{\theta_{3}}}\sum_{i=0}^{\infty}2^{-\alpha ti}=\frac{2^{\alpha t}}{2^{\alpha t}-1}|Q_{j}^{k}|^{\frac{\alpha}{n}t+1}\left(
\fint_{Q_{j}^{k}}v(x)^{\frac{\theta_{3}t}{1-t}}dx\right)^{\frac{1-t}{\theta_{3}}},
\end{align*}
we get
\begin{align*}
\mathcal{S}_{j}^{k}&\leq C(\gamma A^{k+1})^{t}|Q_{j}^{k}|^{\frac{\alpha}{n}t+1}\left(
\fint_{Q_{j}^{k}}v(x)^{\frac{\theta_{3}t}{1-t}}dx\right)^{\frac{1-t}{\theta_{3}}}\\
&\leq CA|Q_{j}^{k}|^{\frac{\alpha}{n}t+1}\left(\fint_{3Q_{j}^{k}}|f(y)|^{\theta_{1}}dy\right)^{\frac{t}{\theta_{1}}}\left(\fint_{3Q_{j}^{k}}|g(z)|^{\theta_{2}}dz\right)^{\frac{t}{\theta_{2}}}\left(
\fint_{Q_{j}^{k}}v(x)^{\frac{\theta_{3}t}{1-t}}dx\right)^{\frac{1-t}{\theta_{3}}},
\end{align*}
where we have used the inequality
$$\left(\fint_{3Q_{j}^{k}}|f(y)|^{\theta_{1}}dy\right)^{\frac{1}{\theta_{1}}}\left(\fint_{3Q_{j}^{k}}|g(z)|^{\theta_{2}}dz\right)^{\frac{1}{\theta_{2}}}>\gamma A^{k}$$
in the last step. By Lemma $\ref{lemma_2}$, we have
\begin{align*}
\mathcal{S}_{j}^{k}&\leq C|Q_{j}^{k}|^{\frac{\alpha}{n}t}\left(\fint_{3Q_{j}^{k}}|f(y)|^{\theta_{1}}dy\right)^{\frac{t}{\theta_{1}}}\left(\fint_{3Q_{j}^{k}}|g(z)|^{\theta_{2}}dz\right)^{\frac{t}{\theta_{2}}}\left(
\fint_{Q_{j}^{k}}v(x)^{\frac{\theta_{3}t}{1-t}}dx\right)^{\frac{1-t}{\theta_{3}}}|Q_{j}^{k}|\\
&\leq C\int_{E_{j}^{k}}|M_{\alpha,\theta_{3}t}^{\theta_{1},\theta_{2}}(f,g,v)(x)|^{t}dx,
\end{align*}
where
$$M_{\alpha,\theta_{3}t}^{\theta_{1},\theta_{2}}(f,g,v)(x)=\sup_{Q\ni x}|Q|^{\frac{\alpha}{n}}\left(\fint_{3Q}|f(y)|^{\theta_{1}}dy\right)^{\frac{1}{\theta_{1}}}\left(\fint_{3Q}|g(z)|^{\theta_{2}}dz\right)^{\frac{1}{\theta_{2}}}\left(
\fint_{Q}v(x)^{\frac{\theta_{3}t}{1-t}}dx\right)^{\frac{1-t}{\theta_{3}t}}.$$
A similar argument gives us the following estimate:
$$\mathcal{S}_{0}\leq C\int_{E_{0}}|M_{\alpha,\theta_{3}t}^{\theta_{1},\theta_{2}}(f,g,v)(x)|^{t}dx.$$
Summing up $\mathcal{S}_{0}$ and $\mathcal{S}_{j}^{k}$, we obtain
$$\mathcal{S}=\mathcal{S}_{0}+\sum_{k,j}\mathcal{S}_{j}^{k}\leq C\int_{Q_{0}}|M_{\alpha,\theta_{3}t}^{\theta_{1},\theta_{2}}(f,g,v)(x)|^{t}dx.$$
By using H\"{o}lder's inequality for $\frac{q_{1}}{a_{*}\theta_{1}}>1$, $\frac{q_{2}}{a_{*}\theta_{2}}>1$ and $\theta_{3}\leq a$ as in Lemma $\ref{lemma_3}$, it yields
\begin{align*}
&M_{\alpha,\theta_{3}t}^{\theta_{1},\theta_{2}}(f,g,v)(x)\leq C\sup_{Q\ni x}|Q|^{\frac{\alpha}{n}}\left(\fint_{3Q}|f(y)w_{1}(y)|^{\frac{q_{1}}{a_{*}}}dy\right)^{\frac{a_{*}}{q_{1}}}\left(\fint_{3Q}|g(z)w_{2}(z)|^{\frac{q_{2}}{a_{*}}}dz\right)^{\frac{a_{*}}{q_{2}}}\\
&\quad\times\left(\fint_{3Q}w_{1}(y)^{-\theta_{1}\left(\frac{q_{1}}{a_{*}\theta_{1}}\right)^{\prime}}dy\right)^{\frac{1}{\theta_{1}\left(\frac{q_{1}}{a_{*}\theta_{1}}\right)^{\prime}}}\left(\fint_{3Q}w_{2}(z)^{-\theta_{2}\left(\frac{q_{2}}{a_{*}\theta_{2}}\right)^{\prime}}dz\right)^{\frac{1}{\theta_{2}\left(\frac{q_{2}}{a_{*}\theta_{2}}\right)^{\prime}}}\left(
\fint_{Q_{j}^{k}}v(x)^{\frac{at}{1-t}}dx\right)^{\frac{1-t}{at}}.
\end{align*}
From Lemma $\ref{lemma_3}$, we have $\theta_{1}\left(\frac{q_{1}}{a_{*}\theta_{1}}\right)^{\prime}\leq\left(\frac{q_{1}}{a}\right)^{\prime}$ and $\theta_{2}\left(\frac{q_{2}}{a_{*}\theta_{2}}\right)^{\prime}\leq\left(\frac{q_{2}}{a}\right)^{\prime}$.

Utilizing H\"{o}lder's inequality, we obtain
\begin{align}\label{estimate for maximal function_1}
&M_{\alpha,\theta_{3}t}^{\theta_{1},\theta_{2}}(f,g,v)(x)\leq C\sup_{Q\ni x}|Q|^{\frac{\alpha}{n}}\left(\fint_{3Q}|f(y)w_{1}(y)|^{\frac{q_{1}}{a_{*}}}dy\right)^{\frac{a_{*}}{q_{1}}}\left(\fint_{3Q}|g(z)w_{2}(z)|^{\frac{q_{2}}{a_{*}}}dz\right)^{\frac{a_{*}}{q_{2}}}\left(\frac{|3Q|}{|Q|}\right)^{\frac{1-s}{as}}\nonumber\\
&\times|3Q|^{\frac{1}{r}}\left(\frac{|Q|}{|3Q|}\right)^{\frac{1-s}{as}}|3Q|^{-\frac{1}{r}}\left(\fint_{3Q}w_{1}(y)^{\left(\frac{q_{1}}{a}\right)^{\prime}}dy\right)^{\frac{1}{\left(\frac{q_{1}}{a}\right)^{\prime}}}\left(\fint_{3Q}w_{2}(z)^{\left(\frac{q_{2}}{a}\right)^{\prime}}dz\right)^{\frac{1}{\left(\frac{q_{2}}{a}\right)^{\prime}}}\left(
\fint_{Q_{j}^{k}}v(x)^{\frac{at}{1-t}}dx\right)^{\frac{1-t}{at}}\nonumber\\
&\leq C[v,\vec{w}]_{t,\vec{q}/a}^{q,as}M_{\alpha-n/r,\vec{q}/a_{*}}(fw_{1},gw_{2})(x),
\end{align}
where we have used the condition $\eqref{two weight condition_1}$ in the last step. For $1\leq s<\infty$, by condition $\eqref{two weight condition_1_1}$, we obtain an estimate which is similar to $\eqref{estimate for maximal function_1}$
\begin{equation}\label{estimate for maximal function_2}
M_{\alpha,\theta_{3}t}^{\theta_{1},\theta_{2}}(f,g,v)(x)\leq C[v,\vec{w}]_{t,\vec{q}/a}^{q,as}M_{\alpha-n/r,\vec{q}/a_{*}}(fw_{1},gw_{2})(x).
\end{equation}
From Lemma $\ref{lemma_5}$, the inequalities $\eqref{estimate for maximal function_1}$ and $\eqref{estimate for maximal function_2}$, it yields that
\begin{align*}
|Q_{0}|^{\frac{1}{s}}\left(\fint_{Q_{0}}|M_{\alpha,\theta_{3}t}^{\theta_{1},\theta_{2}}(f,g,v)(x)|^{t}dx\right)^{\frac{1}{t}}&\leq C[v,\vec{w}]_{t,\vec{q}/a}^{q,as}\|M_{\alpha-n/r,\vec{q}/a_{*}}(fw_{1},gw_{2})\|_{\mathcal{M}_{t}^{s}}\\
&\leq C[v,\vec{w}]_{t,\vec{q}/a}^{q,as}\sup_{Q\in\mathcal{D}(\mathbb{R}^{n})}|Q|^{\frac{1}{p}}\left(\fint_{Q}(fw_{1})^{q_{1}}\right)^{\frac{1}{q_{1}}}\left(\fint_{Q}(fw_{2})^{q_{2}}\right)^{\frac{1}{q_{2}}},
\end{align*}
where
$$\frac{\alpha}{n}>\frac{1}{r},\quad\frac{1}{s}=\frac{1}{p}-\frac{\alpha-n/r}{n}\quad\text{and}\quad\frac{t}{s}=\frac{q}{p}.$$
Therefore we have
\begin{equation}\label{final estimate}
\|\mathscr{T}_{1}v\|_{\mathcal{M}_{t}^{s}}\leq C[v,\vec{w}]_{t,\vec{q}/a}^{q,as}\sup_{Q\in\mathcal{D}(\mathbb{R}^{n})}|Q|^{\frac{1}{p}}\left(\fint_{Q}(fw_{1})^{q_{1}}\right)^{\frac{1}{q_{1}}}\left(\fint_{Q}(fw_{2})^{q_{2}}\right)^{\frac{1}{q_{2}}}.
\end{equation}

{\noindent}{\bf Estimate of $\mathscr{T}_{2}(x)$.} A normalization allows us to assume that
\begin{equation}\label{normalization}
\sup_{Q\in\mathcal{D}(\mathbb{R}^{n})}|Q|^{\frac{1}{p}}\left(\fint_{Q}|f(y)w_{1}(y)|^{q_{1}}dy\right)^{\frac{1}{q_{1}}}\left(\fint_{Q}|g(z)w_{2}(z)|^{q_{2}}dz\right)^{\frac{1}{q_{2}}}=1.
\end{equation}
Combining the equalities $\eqref{define as average of BMO function}$-$\eqref{part two}$, H\"{o}lder's inequality with the definition of $\mathscr{T}_{2}(x)$, we have the following estimate:
\begin{align*}
&\int_{Q_{0}}(|\mathscr{T}_{2}(x)|v(x))^{t}dx\leq\sum_{A\subseteq\{1,\cdots,m\}}\sum_{B\subseteq\{m+1,\cdots,N\}}\mathop{\sum_{Q\supsetneq Q_{0}}}_{Q\in\mathcal{D}(\mathbb{R}^{n})}|Q|^{(\frac{\alpha}{n}-1)t}\left[\int_{Q}\int_{|y|_{\infty}\leq \ell(Q)}\prod_{i\in\bar{A}}|b_{i}(x-y)-\lambda_{i}|\right.\\
&\qquad\qquad\quad\times\left.\prod_{i\in\bar{B}}|b_{i}(x+y)-\lambda_{i}|f(x-y)g(x+y)dydx\right]^{t}\left[\int_{Q_{0}}\prod_{i\in A\cup B}|b_{i}(x)-\lambda_{i}|^{\frac{t}{1-t}}v(x)^{\frac{t}{1-t}}dx\right]^{1-t}\\
&\leq \sum_{A\subseteq\{1,\cdots,m\}}\sum_{B\subseteq\{m+1,\cdots,N\}}\mathop{\sum_{Q\supsetneq Q_{0}}}_{Q\in\mathcal{D}(\mathbb{R}^{n})}|Q|^{(\frac{\alpha}{n}+1)t}\left[\fint_{3Q}\prod_{i\in\bar{A}}|b_{i}(y)-\lambda_{i}|f(y)dy\fint_{3Q}\prod_{i\in\bar{B}}|b_{i}(z)-\lambda_{i}|g(z)dz\right]^{t}\nonumber\\
&\qquad\qquad\quad\times\left[\int_{Q_{0}}\prod_{i\in A\cup B}|b_{i}(x)-\lambda_{i}|^{\frac{t}{1-t}}v(x)^{\frac{t}{1-t}}dx\right]^{1-t}.
\end{align*}
Using H\"{o}lder's inquality for $\theta_{4}>1$ and $\frac{q_{1}}{\theta_{4}}>1$ as in Lemma $\ref{lemma_3}$, we obtain
\begin{align*}
\fint_{3Q}\prod_{i\in\bar{A}}|b_{i}(y)-\lambda_{i}|f(y)dy&\leq\left(\fint_{3Q}\prod_{i\in\bar{A}}|b_{i}(y)-\lambda_{i}|^{\theta_{4}^{\prime}}dy\right)^{\frac{1}{\theta_{4}^{\prime}}}\left(\fint_{3Q}|f(y)|^{\theta_{4}}dy\right)^{\frac{1}{\theta_{4}}}\\
&\leq\prod_{i\in\bar{A}}\|b_{i}\|_{\rm BMO}\left(\fint_{3Q}|f(y)w_{1}(y)|^{q_{1}}dy\right)^{\frac{1}{q_{1}}}\left(\fint_{3Q}w_{1}(y)^{-\theta_{4}\left(\frac{q_{1}}{\theta_{4}}\right)^{\prime}}dy\right)^{\frac{1}{\theta_{4}\left(\frac{q_{1}}{\theta_{4}}\right)^{\prime}}}.
\end{align*}
Similarly, we have
$$\fint_{3Q}\prod_{i\in\bar{B}}|b_{i}(z)-\lambda_{i}|g(z)dz\leq\prod_{i\in\bar{B}}\|b_{i}\|_{\rm BMO}\left(\fint_{3Q}|g(z)w_{2}(z)|^{q_{2}}dz\right)^{\frac{1}{q_{2}}}\left(\fint_{3Q}w_{2}(z)^{-\theta_{5}\left(\frac{q_{5}}{\theta_{5}}\right)^{\prime}}dz\right)^{\frac{1}{\theta_{5}\left(\frac{q_{2}}{\theta_{5}}\right)^{\prime}}}.$$
By $\eqref{normalization}$, we obtain
\begin{align*}
&\fint_{3Q}\prod_{i\in\bar{A}}|b_{i}(y)-\lambda_{i}|f(y)dy\fint_{3Q}\prod_{i\in\bar{B}}|b_{i}(z)-\lambda_{i}|g(z)dz\\
&\leq C|3Q|^{-\frac{1}{p}}\prod_{i\in\bar{A}\cup\bar{B}}\|b_{i}\|_{\rm BMO}\left(\fint_{3Q}w_{1}(y)^{-\theta_{4}\left(\frac{q_{1}}{\theta_{4}}\right)^{\prime}}dy\right)^{\frac{1}{\theta_{4}\left(\frac{q_{1}}{\theta_{4}}\right)^{\prime}}}\left(\fint_{3Q}w_{2}(z)^{-\theta_{5}\left(\frac{q_{5}}{\theta_{5}}\right)^{\prime}}dz\right)^{\frac{1}{\theta_{5}\left(\frac{q_{2}}{\theta_{5}}\right)^{\prime}}}.
\end{align*}
Since we have the assumption that $a\geq\Big\{\frac{q_{1}}{\left(\theta_{4}\left({q_{1}}/{\theta_{4}}\right)^{\prime}\right)^{\prime}},\frac{q_{2}}{\left(\theta_{5}\left({q_{2}}/{\theta_{2}}\right)^{\prime}\right)^{\prime}}\Big\}>1$, it can be inferred from the H\"{o}lder's inequality that
$$
\left(\fint_{3Q}w_{1}(y)^{-\theta_{4}\left(\frac{q_{1}}{\theta_{4}}\right)^{\prime}}dy\right)^{\frac{1}{\theta_{4}\left(\frac{q_{1}}{\theta_{4}}\right)^{\prime}}}\leq\left(\fint_{3Q}w_{1}(y)^{-\left(\frac{q_{1}}{a}\right)^{\prime}}dy\right)^{\frac{1}{\left(\frac{q_{1}}{a}\right)^{\prime}}}
$$
and
$$
\left(\fint_{3Q}w_{2}(z)^{-\theta_{5}\left(\frac{q_{2}}{\theta_{5}}\right)^{\prime}}dz\right)^{\frac{2}{\theta_{5}\left(\frac{q_{2}}{\theta_{5}}\right)^{\prime}}}\leq\left(\fint_{3Q}w_{2}(z)^{-\left(\frac{q_{2}}{a}\right)^{\prime}}dz\right)^{\frac{1}{\left(\frac{q_{2}}{a}\right)^{\prime}}}.
$$
Then we have
\begin{align}\label{f,g controlled by BMO norm}
&\fint_{3Q}\prod_{i\in\bar{A}}|b_{i}(y)-\lambda_{i}|f(y)dy\fint_{3Q}\prod_{i\in\bar{B}}|b_{i}(z)-\lambda_{i}|g(z)dz\nonumber\\
&\leq C|3Q|^{-\frac{1}{p}}\prod_{i\in\bar{A}\cup\bar{B}}\|b_{i}\|_{\rm BMO}\left(\fint_{3Q}w_{1}(y)^{-\left(\frac{q_{1}}{a}\right)^{\prime}}dy\right)^{\frac{1}{\left(\frac{q_{1}}{a}\right)^{\prime}}}\left(\fint_{3Q}w_{2}(z)^{-\left(\frac{q_{2}}{a}\right)^{\prime}}dz\right)^{\frac{1}{\left(\frac{q_{2}}{a}\right)^{\prime}}}.
\end{align}
By H\"{o}lder's inequality for $a\geq\theta_{3}>1$ as in Lemma $\ref{lemma_3}$, we obtain
\begin{align*}
&\left(\int_{Q_{0}}\prod_{i\in A\cup B}|b_{i}(x)-\lambda_{i}|^{\frac{t}{1-t}}v(x)^{\frac{t}{1-t}}dx\right)^{1-t}\\
&\leq\left(\int_{Q_{0}}\prod_{i\in A\cup B}|b_{i}(x)-\lambda_{i}|^{\frac{\theta_{3}^{\prime}t}{1-t}}dx\right)^{\frac{1-t}{\theta_{3}^{\prime}}}\left(\int_{Q_{0}}v(x)^{\frac{t\theta_{3}}{1-t}}dx\right)^{\frac{1-t}{\theta_{3}}}\\
&\leq C\prod_{i\in A\cup B}\left(\fint_{3Q_{0}}|b_{i}(x)-\lambda_{i}|^{\frac{\theta_{3}^{\prime}t|A\cup B|}{1-t}}dx\right)^{\frac{1-t}{\theta_{3}^{\prime}|A\cup B|}}\left(\fint_{Q_{0}}v(x)^{\frac{t\theta_{3}}{1-t}}dx\right)^{\frac{1-t}{\theta_{3}}}|Q_{0}|^{1-t}\\
&\leq C\prod_{i\in A\cup B}\left(\fint_{3Q_{0}}|b_{i}(x)-\lambda_{i}|^{\frac{\theta_{3}^{\prime}t|A\cup B|}{1-t}}dx\right)^{\frac{1-t}{\theta_{3}^{\prime}|A\cup B|}}\left(\fint_{Q_{0}}v(x)^{\frac{at}{1-t}}dx\right)^{\frac{1-t}{a}}|Q_{0}|^{1-t}.
\end{align*}
We now evaluate $|b_{i}-\lambda_{i}|$. If $Q\supsetneq Q_{0}$ and $Q\in\mathcal{D}(\mathbb{R}^{n})$, then there exists $k=1,2,\cdots$ such that $Q_{k}:=Q$, $Q_{j}\in\mathcal{D}(\mathbb{R}^{n})$, $Q_{j}\supsetneq Q_{j-1}$ and $|Q_{j}|=2^{n}|Q_{j-1}|$ $(j=1,2,\cdots,k)$. By Lemma $\ref{lemma_5}$, we have
\begin{equation}\label{controlled by BMO norm}
|b_{i}-\lambda_{i}|\leq|b_{i}(x)-m_{3Q_{0}}(b_{i})|+k2^{n}\|b_{i}\|_{\rm BMO}.
\end{equation}
Using the inequality $\eqref{controlled by BMO norm}$, Minkowski's inequality and Lemma $\ref{lemma_1}$, we obtain
\begin{align*}
&\left(\fint_{3Q_{0}}|b_{i}(x)-\lambda_{i}|^{\frac{\theta_{3}^{\prime}t|A\cup B|}{1-t}}dx\right)^{\frac{1-t}{\theta_{3}^{\prime}|A\cup B|}}\\
&\leq
\left(\fint_{3Q_{0}}\left(|b_{i}(x)-m_{3Q_{0}}(b_{i})|+k2^{n}\|b_{i}\|_{\rm BMO}\right)^{\frac{\theta_{3}^{\prime}t|A\cup B|}{1-t}}dx\right)^{\frac{1-t}{\theta_{3}^{\prime}|A\cup B|}}\\
&\leq\left(\fint_{3Q_{0}}|b_{i}(x)-m_{3Q_{0}}(b_{i})|^{\frac{\theta_{3}^{\prime}t|A\cup B|}{1-t}}dx\right)^{\frac{1-t}{\theta_{3}^{\prime}|A\cup B|}}+k2^{n}\|b_{i}\|_{\rm BMO}\\
&\leq C(1+k2^{n})\|b_{i}\|_{\rm BMO}.
\end{align*}
Hence,
\begin{align}\label{v controlled by BMO norm}
\left(\int_{Q_{0}}\prod_{i\in A\cup B}|b_{i}(x)-\lambda_{i}|^{\frac{t}{1-t}}v(x)^{\frac{t}{1-t}}dx\right)^{1-t}\leq Ck^{t}\prod_{i\in A\cup B}\|b_{i}\|_{*}^{t}\left(\fint_{Q_{0}}v(x)^{\frac{at}{1-t}}dx\right)^{\frac{1-t}{a}}|Q_{0}|^{1-t}.
\end{align}
Next, we consider the Morrey norm of $\mathscr{T}_{2}(x)$ as follows,
$$|Q_{0}|^{\frac{1}{s}}\left(\fint_{Q_{0}}|\mathscr{T}_{2}(x)v(x)|^{t}dx\right)^{\frac{1}{t}}=\left(|Q_{0}|^{\frac{t}{s}}\fint_{Q_{0}}|\mathscr{T}_{2}(x)v(x)|^{t}dx\right)^{\frac{1}{t}}.$$
Combining the inequalities $\eqref{f,g controlled by BMO norm}$, $\eqref{v controlled by BMO norm}$ with the condition $\eqref{two weight condition_1}$, we have
\begin{align*}
&|Q_{0}|^{\frac{t}{s}}\fint_{Q_{0}}|\mathscr{T}_{2}(x)v(x)|^{t}dx\leq C\sum_{A\subseteq\{1,\cdots,m\}}\sum_{B\subseteq\{m+1,\cdots,N\}}\prod_{i\in\bar{A}\cup\bar{B}}\|b_{i}\|_{*}^{t}\mathop{\sum_{Q\supsetneq Q_{0}}}_{Q\in\mathcal{D}(\mathbb{R}^{n})}|Q|^{(\frac{\alpha}{n}+1)t}|Q_{0}|^{\frac{t}{s}-1}|3Q|^{-\frac{t}{p}}\\
&\times\left(\fint_{3Q}w_{1}(y)^{-\left(\frac{q_{1}}{a}\right)^{\prime}}dy\right)^{\frac{1}{\left(\frac{q_{1}}{a}\right)^{\prime}}}\left(\fint_{3Q}w_{2}(z)^{-\left(\frac{q_{2}}{a}\right)^{\prime}}dz\right)^{\frac{1}{\left(\frac{q_{2}}{a}\right)^{\prime}}}\left(\int_{Q_{0}}\prod_{i\in A\cup B}|b_{i}(x)-\lambda_{i}|^{\frac{t}{1-t}}v(x)^{\frac{t}{1-t}}dx\right)^{1-t}\\
&\leq C\sum_{A\subseteq\{1,\cdots,m\}}\sum_{B\subseteq\{m+1,\cdots,N\}}\prod_{i\in\bar{A}\cup\bar{B}}\|b_{i}\|_{*}^{t}\prod_{i\in{A}\cup{B}}\|b_{i}\|_{*}^{t}\sum_{k=1}^{\infty}k^{t}\mathop{\sum_{Q_{k}\in\mathcal{D}(\mathbb{R}^{n})}}_{Q_{k}\supseteq Q,\,|Q_{k}|=2^{kn}|Q_{0}|}|Q_{k}|^{(\frac{\alpha}{n}+1)t}|Q_{0}|^{\frac{t}{s}-t}\\
&\qquad\times|3Q_{k}|^{-\frac{t}{p}}\left(\fint_{3Q_{k}}w_{1}(y)^{-\left(\frac{q_{1}}{a}\right)^{\prime}}dy\right)^{\frac{1}{\left(\frac{q_{1}}{a}\right)^{\prime}}}\left(\fint_{3Q_{k}}w_{2}(z)^{-\left(\frac{q_{2}}{a}\right)^{\prime}}dz\right)^{\frac{1}{\left(\frac{q_{2}}{a}\right)^{\prime}}}\left(\fint_{Q_{0}}v(x)^{\frac{at}{1-t}}dx\right)^{\frac{1-t}{a}}\\
&\leq C\|\vec{b}\|_{{\rm BMO}^{N}}^{t}\sum_{k=1}^{\infty}k^{t}\mathop{\sum_{Q_{k}\in\mathcal{D}(\mathbb{R}^{n})}}_{Q_{k}\supseteq Q,\,|Q_{k}|=2^{kn}|Q_{0}|}\left[\left(\frac{|Q_{0}|}{|3Q_{k}|}\right)^{\frac{1-s}{as}}|3Q_{k}|^{\frac{1}{r}}\left(\fint_{3Q_{k}}w_{1}(y)^{-\left(\frac{q_{1}}{a}\right)^{\prime}}dy\right)^{\frac{1}{\left(\frac{q_{1}}{a}\right)^{\prime}}}\right.\\
&\qquad\qquad\qquad\times\left.\left(\fint_{3Q_{k}}w_{2}(z)^{-\left(\frac{q_{2}}{a}\right)^{\prime}}dz\right)^{\frac{1}{\left(\frac{q_{2}}{a}\right)^{\prime}}}\left(\fint_{Q_{0}}v(x)^{\frac{at}{1-t}}dx\right)^{\frac{1-t}{at}}\right]^{t}\left(\frac{|Q_{0}|}{|3Q_{k}|}\right)^{\frac{(1-s)t}{s}(1-\frac{1}{a})}\\
&\leq C\|\vec{b}\|_{{\rm BMO}^{N}}^{t}\left([v,\vec{w}]_{t,\vec{q}/a}^{r,as}\right)^{t}\sum_{k=1}^{\infty}k^{t}\mathop{\sum_{Q_{k}\in\mathscr{D}}}_{Q_{k}\supseteq Q,\,|Q_{k}|=2^{kn}|Q_{0}|}2^{\frac{(1-s)ntk}{s}(1-\frac{1}{a})}\\
&\leq C\|\vec{b}\|_{{\rm BMO}^{N}}^{t}\left([v,\vec{w}]_{t,\vec{q}/a}^{r,as}\right)^{t}\sum_{k=1}^{\infty}k^{t}2^{\frac{(1-s)ntk}{s}(1-\frac{1}{a})}\\
&\leq C\|\vec{b}\|_{{\rm BMO}^{N}}^{t}\left([v,\vec{w}]_{t,\vec{q}/a}^{r,as}\right)^{t}.
\end{align*}
Then for $0<s<1$, there holds
\begin{equation}\label{final estimate_2}
|Q_{0}|^{\frac{1}{s}}\left(\fint_{Q_{0}}|\mathscr{T}_{2}(x)v(x)|^{t}dx\right)^{\frac{1}{t}}\leq C\|\vec{b}\|_{{\rm BMO}^{N}}[v,\vec{w}]_{t,\vec{q}/a}^{r,as}.
\end{equation}
For $1\leq s<\infty$, applying condition $\eqref{two weight condition_1_1}$, the inequality $\eqref{final estimate_2}$ also holds.

Combining $\eqref{final estimate}$ with $\eqref{final estimate_2}$, it yields that
\begin{align*}
&\|B_{\alpha}(f,g)v\|_{\mathcal{M}_{t}^{s}}\\
&\leq C\|\vec{b}\|_{{\rm BMO}^{N}}[v,\vec{w}]_{t,\vec{q}/a}^{r,as}\sup_{Q\in\mathcal{D}(\mathbb{R}^{n})}|Q|^{\frac{1}{p}}\left(\fint_{Q}|f(y)w_{1}(y)|^{q_{1}}dy\right)^{\frac{1}{q_{1}}}\left(\fint_{Q}|g(z)w_{z}(z)|^{q_{2}}dz\right)^{\frac{1}{q_{2}}}.
\end{align*}
This completes the proof of Theorem $\ref{main_2}$. $\hfill$ $\Box$

\section{The proof Theorem $\ref{main_3}$}
Let us keep using the notations $\mathscr{T}_{1}$ and $\mathscr{T}_{2}$ as in Theorem $\ref{main_2}$.

{\noindent}{\bf Eestimate of $\mathscr{T}_{1}$.} By the duality argument, we have
$$\left(\int_{Q_{0}}|\mathscr{T}_{1}(x)v(x)|^{t}dx\right)^{\frac{1}{t}}=\sup_{\|h\|_{L^{t^{\prime}}(Q_{0})}=1}\left(\int_{Q_{0}}\mathscr{T}_{1}(x)v(x)h(x)dx\right).$$

Without loss of generality, we may still assume that $f$ and $g$ are nonnegative, bounded and compactly supported. From $\eqref{estimate for T_{1}}$, we have
\begin{align*}
&\int_{Q_{0}}\mathscr{T}_{1}(x)v(x)h(x)dx\leq \sum_{A\subseteq\{1,\cdots,m\}}\sum_{B\subseteq\{m+1,\cdots,N\}}\sum_{Q\in\mathscr{D}(Q_{0})}|Q|^{\frac{\alpha}{n}-1}\int_{Q}\mathop{\int}_{|y|_{\infty}\leq \ell(Q)}\prod_{i\in\bar{A}}|b_{i}(x-y)-\lambda_{i}|\\
&\qquad\qquad\qquad\qquad\qquad\quad\times\prod_{i\in\bar{B}}|b_{i}(x+y)-\lambda_{i}|f(x-y)g(x+y)dy\prod_{i\in A\cup B}|b_{i}(x)-\lambda_{i}|v(x)h(x)dx.
\end{align*}
If we use H\"{o}lder's inequality with the pair $(r_{1},r_{2})$ for the inner integral, and then perform a change in variables, we will arrive at
\begin{align*}
&\int_{Q_{0}}\mathscr{T}_{1}(x)v(x)h(x)dx\leq \sum_{A\subseteq\{1,\cdots,m\}}\sum_{B\subseteq\{m+1,\cdots,N\}}\sum_{Q\in\mathscr{D}(Q_{0})}|Q|^{\frac{\alpha}{n}-1}\left[\int_{3Q}\prod_{i\in\bar{A}}(|b_{i}(y)-\lambda_{i}|f(y))^{r_{1}}dy\right]^{\frac{1}{r_{1}}}\\
&\qquad\qquad\qquad\qquad\qquad\quad\times\left[\int_{3Q}\prod_{i\in\bar{B}}(|b_{i}(z)-\lambda_{i}|g(z))^{r_{2}}dz\right]^{\frac{1}{r_{2}}}\int_{Q}\prod_{i\in A\cup B}|b_{i}(x)-\lambda_{i}|v(x)h(x)dx.
\end{align*}
We now use H\"{o}lder's inequality for $\vartheta_{1},\vartheta_{2},\vartheta_{3}>1$ as in Lemma $\ref{lemma_4}$ and Lemma $\ref{lemma_2}$ to get the following estimates,
\begin{align}\label{estimate for f}
\fint_{3Q}\prod_{i\in\bar{A}}(|b_{i}(y)-\lambda_{i}|f(y))^{r_{1}}dy&\leq\prod_{i\in\bar{A}}\left(\fint_{3Q}|b_{i}(y)-\lambda_{i}|^{\vartheta_{1}^{\prime}r_{1}|\bar{A}|}dy\right)^{\frac{1}{\vartheta_{1}^{\prime}|\bar{A}|}}\left(\fint_{3Q}|f(y)|^{\vartheta_{1}r_{1}}dy\right)^{\frac{1}{\vartheta_{1}}}\nonumber\\
&\leq C\prod_{i\in\bar{A}}\|b_{i}\|_{\rm BMO}^{r_{1}}\left(\fint_{3Q}|f(y)|^{\vartheta_{1}r_{1}}dy\right)^{\frac{1}{\vartheta_{1}}},
\end{align}

\begin{align}\label{estimate for g}
\fint_{3Q}\prod_{i\in\bar{B}}(|b_{i}(z)-\lambda_{i}|g(z))^{r_{2}}dz&\leq\prod_{i\in\bar{B}}\left(\fint_{3Q}|b_{i}(z)-\lambda_{i}|^{\vartheta_{2}^{\prime}r_{2}|\bar{B}|}dz\right)^{\frac{1}{\vartheta_{2}^{\prime}|\bar{B}|}}\left(\fint_{3Q}|g(z)|^{\vartheta_{2}r_{2}}dz\right)^{\frac{1}{\vartheta_{2}}}\nonumber\\
&\leq C\prod_{i\in\bar{B}}\|b_{i}\|_{\rm BMO}^{r_{2}}\left(\fint_{3Q}|g(z)|^{\vartheta_{2}r_{2}}dz\right)^{\frac{1}{\vartheta_{2}}}
\end{align}
and
\begin{align}\label{estimate for h}
&\fint_{Q}\prod_{i\in A\cup B}|b_{i}(x)-\lambda_{i}|v(x)h(x)dx\nonumber\\
&\leq\prod_{i\in A\cup B}\left(\fint_{Q}|b_{i}(x)-\lambda_{i}|^{\vartheta_{3}^{\prime}|A\cup B|}dz\right)^{\frac{1}{\vartheta_{3}^{\prime}|A\cup B|}}\left(\fint_{Q}(v(x)h(x))^{\vartheta_{3}}dx\right)^{\frac{1}{\vartheta_{3}}}\nonumber\\
&\leq C\prod_{i\in A\cup B}\|b_{i}\|_{\rm BMO}\left(\fint_{Q}(v(x)h(x))^{\vartheta_{3}}dx\right)^{\frac{1}{\vartheta_{3}}}.
\end{align}
By the above three estimates $\eqref{estimate for f}$, $\eqref{estimate for g}$ and $\eqref{estimate for h}$, we come up with the following estimate:
\begin{align}\label{T_1 for t>1}
&\int_{Q_{0}}\mathscr{T}_{1}(x)v(x)h(x)\leq C\|\vec{b}\|_{BMO^{N}}\sum_{Q\in\mathcal{D}(Q_{0})}|Q|^{\frac{\alpha}{n}+1}\nonumber\\
&\qquad\times\left(\fint_{3Q}|f(y)|^{\vartheta_{1}r_{1}}dy\right)^{\frac{1}{\vartheta_{1}r_{1}}}\left(\fint_{3Q}|g(z)|^{\vartheta_{2}r_{2}}dz\right)^{\frac{1}{\vartheta_{2}r_{2}}}\left(\fint_{Q}(v(x)h(x))^{\vartheta_{3}}dx\right)^{\frac{1}{\vartheta_{3}}}.
\end{align}
To estimate $\eqref{T_1 for t>1}$, we will consider $\mathcal{D}(Q_{0})$ divided into two parts. Let
$$\mathcal{D}_{0}^{\prime}(Q_{0}):=\left\{Q\in\mathcal{D}(Q_{0}):\,\left(\fint_{3Q}|f(y)|^{\vartheta_{1}r_{1}}dy\right)^{\frac{1}{\vartheta_{1}r_{1}}}\left(\fint_{3Q}|g(z)|^{\vartheta_{2}r_{2}}dz\right)^{\frac{1}{\vartheta_{2}r_{2}}}\leq \gamma^{\prime} A^{\prime}\right\}$$
and ${\mathcal{D}_{j}^{k}}^{\prime}(Q_{0})$ denote the following family of cubes $Q$:
$$\left\{Q\in\mathcal{D}(Q_{0}):\,Q\subseteq {Q_{j}^{k}}^{\prime},\gamma^{\prime} {A^{\prime}}^{k}<\left(\fint_{3Q}|f(y)|^{\vartheta_{1}r_{1}}dy\right)^{\frac{1}{\vartheta_{1}r_{1}}}\left(\fint_{3Q}|g(z)|^{\vartheta_{2}r_{2}}dz\right)^{\frac{1}{\vartheta_{2}r_{2}}}\leq \gamma^{\prime} {A^{\prime}}^{k+1}\right\},$$
where ${Q_{j}^{k}}^{\prime}$ is as in Lemma $\ref{lemma_2}$. Then we have
\begin{equation}\label{divide into two parts}
\mathcal{D}(Q_{0})=\mathcal{D}_{0}^{\prime}(Q_{0})\bigcup\Big(\bigcup_{k,j}{\mathcal{D}_{j}^{k}}^{\prime}(Q_{0})\Big).
\end{equation}
By the inequality $\eqref{T_1 for t>1}$ and the equality $\eqref{divide into two parts}$, we obtain
\begin{align*}
\int_{Q_{0}}\mathscr{T}_{1}(x)v(x)h(x)&\leq C\|\vec{b}\|_{{\rm BMO}^{N}}\Big(\sum_{Q\in\mathcal{D}_{0}^{\prime}(Q_{0})}+\sum_{k,j}\sum_{Q\in{\mathcal{D}_{j}^{k}}^{\prime}(Q_{0})}\Big)|Q|^{\frac{\alpha}{n}+1}\nonumber\\
&\times\left(\fint_{3Q}|f(y)|^{\vartheta_{1}r_{1}}dy\right)^{\frac{1}{\vartheta_{1}r_{1}}}\left(\fint_{3Q}|g(z)|^{\vartheta_{2}r_{2}}dz\right)^{\frac{1}{\vartheta_{2}r_{2}}}\left(\fint_{Q}(v(x)h(x))^{\vartheta_{3}}dx\right)^{\frac{1}{\vartheta_{3}}}\\
&=:C\|\vec{b}\|_{{\rm BMO}^{N}}\Big(\mathcal{T}_{0}+\sum_{k,j}\mathcal{T}_{j}^{k}\Big).
\end{align*}
Next, we will estimate for $\mathcal{T}_{j}^{k}$. If $Q\in{\mathcal{D}_{j}^{k}}^{\prime}(Q_{0})$, one has
$$
\left(\fint_{3Q}|f(y)|^{\vartheta_{1}r_{1}}dy\right)^{\frac{1}{\vartheta_{1}r_{1}}}\left(\fint_{3Q}|g(z)|^{\vartheta_{2}r_{2}}dz\right)^{\frac{1}{\vartheta_{2}r_{2}}}\leq \gamma^{\prime} {A^{\prime}}^{k+1}.
$$
Hence, we obtain
\begin{align*}
\mathcal{T}_{j}^{k}&\leq\gamma^{\prime} {A^{\prime}}^{k+1}\sum_{Q\in{\mathcal{D}_{j}^{k}}^{\prime}(Q_{0})}|Q|^{\frac{\alpha}{n}+1}\left(\fint_{Q}(v(x)h(x))^{\vartheta_{3}}dx\right)^{\frac{1}{\vartheta_{3}}}\\
&\leq\gamma^{\prime} {A^{\prime}}^{k+1}|{Q_{j}^{k}}^{\prime}|^{\frac{\alpha}{n}}\sum_{Q\in{\mathcal{D}_{j}^{k}}^{\prime}(Q_{0})}|Q|\left(\fint_{Q}(v(x)h(x))^{\vartheta_{3}}dx\right)^{\frac{1}{\vartheta_{3}}}\\
&\leq\gamma^{\prime} {A^{\prime}}^{k+1}|{Q_{j}^{k}}^{\prime}|^{\frac{\alpha}{n}}\sum_{Q\in{\mathcal{D}_{j}^{k}}^{\prime}(Q_{0})}\int_{Q}\left(\fint_{Q}(v(y)h(y))^{\vartheta_{3}}dy\right)^{\frac{1}{\vartheta_{3}}}dx.
\end{align*}
Since
$$\gamma^{\prime} {A^{\prime}}^{k}<\left(\fint_{3{Q_{j}^{k}}^{\prime}}|f(y)|^{\vartheta_{1}r_{1}}dy\right)^{\frac{1}{\vartheta_{1}r_{1}}}\left(\fint_{3{Q_{j}^{k}}^{\prime}}|g(z)|^{\vartheta_{2}r_{2}}dz\right)^{\frac{1}{\vartheta_{2}r_{2}}},$$
we have
\begin{align*}
\mathcal{T}_{j}^{k}&\leq A^{\prime}|{Q_{j}^{k}}^{\prime}|^{\frac{\alpha}{n}}\left(\fint_{3{Q_{j}^{k}}^{\prime}}|f(y)|^{\vartheta_{1}r_{1}}dy\right)^{\frac{1}{\vartheta_{1}r_{1}}}\left(\fint_{3{Q_{j}^{k}}^{\prime}}|g(z)|^{\vartheta_{2}r_{2}}dz\right)^{\frac{1}{\vartheta_{2}r_{2}}}\\
&\quad\times\sum_{Q\in{\mathcal{D}_{j}^{k}}^{\prime}(Q_{0})}\int_{Q}M[(v_{j}^{k}h)^{\vartheta_{3}}](x)^{\frac{1}{\vartheta_{3}}}dx,
\end{align*}
where $v_{j}^{k}=v\chi_{{Q_{j}^{k}}^{\prime}}$ and the symbol $M$ is the ordinary Hardy-Littlewood maximal operator.
By Lemma $\ref{lemma_2}$, we obtain the following estimate for $\mathcal{T}_{j}^{k}$
\begin{align*}
&\mathcal{T}_{j}^{k}\leq A^{\prime}|{Q_{j}^{k}}^{\prime}|^{\frac{\alpha}{n}+1}\left(\fint_{3{Q_{j}^{k}}^{\prime}}|f(y)|^{\vartheta_{1}r_{1}}dy\right)^{\frac{1}{\vartheta_{1}r_{1}}}\left(\fint_{3{Q_{j}^{k}}^{\prime}}|g(z)|^{\vartheta_{2}r_{2}}dz\right)^{\frac{1}{\vartheta_{2}r_{2}}}\fint_{{Q_{j}^{k}}^{\prime}}M[(v_{j}^{k}h)^{\vartheta_{3}}](x)^{\frac{1}{\vartheta_{3}}}dx\\
&\leq2A^{\prime}|{E_{j}^{k}}^{\prime}||{Q_{j}^{k}}^{\prime}|^{\frac{\alpha}{n}}\left(\fint_{3{Q_{j}^{k}}^{\prime}}|f(y)|^{\vartheta_{1}r_{1}}dy\right)^{\frac{1}{\vartheta_{1}r_{1}}}\left(\fint_{3{Q_{j}^{k}}^{\prime}}|g(z)|^{\vartheta_{2}r_{2}}dz\right)^{\frac{1}{\vartheta_{2}r_{2}}}\fint_{{Q_{j}^{k}}^{\prime}}M[(v_{j}^{k}h)^{\vartheta_{3}}](x)^{\frac{1}{\vartheta_{3}}}dx\\
&=2A^{\prime}\int_{{E_{j}^{k}}^{\prime}}|{Q_{j}^{k}}^{\prime}|^{\frac{\alpha}{n}}\left(\fint_{3{Q_{j}^{k}}^{\prime}}|f(y_{1})|^{\vartheta_{1}r_{1}}dy_{1}\right)^{\frac{1}{\vartheta_{1}r_{1}}}\left(\fint_{3{Q_{j}^{k}}^{\prime}}|g(z)|^{\vartheta_{2}r_{2}}dz\right)^{\frac{1}{\vartheta_{2}r_{2}}}\fint_{{Q_{j}^{k}}^{\prime}}M[(v_{j}^{k}h)^{\vartheta_{3}}](x)^{\frac{1}{\vartheta_{3}}}dxdy.
\end{align*}
We may take $c\in(t,r)$ and $L>1$ as in Lemma $\ref{lemma_4}$. Using H\"{o}lder's inequality for $c>1$, we have
$$
M[(v_{j}^{k}h)^{\vartheta_{3}}](x)^{\frac{1}{\vartheta_{3}}}\leq M[(v_{j}^{k})^{c\vartheta_{3}}](x)^{\frac{1}{c\vartheta_{3}}}M[h^{c^{\prime}\vartheta_{3}}](x)^{\frac{1}{c^{\prime}\vartheta_{3}}}.
$$
By H\"{o}lder's inequality for $Lt>1$, we obtain the following inequality:
\begin{align*}
\fint_{{Q_{j}^{k}}^{\prime}}M[(v_{j}^{k}h)^{\vartheta_{3}}](x)^{\frac{1}{\vartheta_{3}}}dx\leq\left(\fint_{{E_{j}^{k}}^{\prime}}M[(v_{j}^{k})^{c\vartheta_{3}}](x)^{\frac{Lt}{c\vartheta_{3}}}dx\right)^{\frac{1}{Lt}}\left(\fint_{{E_{j}^{k}}^{\prime}}M[h^{c^{\prime}\vartheta_{3}}](x)^{\frac{(Lt)^{\prime}}{c^{\prime}\vartheta_{3}}}dx\right)^{\frac{1}{(Lt)^{\prime}}}.
\end{align*}
Since $c\vartheta_{3}<Lt$, the boundedness of $M:\,L^{\frac{Lt}{c\vartheta_{3}}}\rightarrow L^{\frac{Lt}{c\vartheta_{3}}}$ gives us the following inequality
\begin{align*}
\fint_{{Q_{j}^{k}}^{\prime}}M[(v_{j}^{k}h)^{\vartheta_{3}}](x)^{\frac{1}{\vartheta_{3}}}dx\leq\left(\frac{1}{|{Q_{j}^{k}}^{\prime}|}\int_{\mathbb{R}^{n}}(v_{j}^{k}(x))^{Lt}dx\right)^{\frac{1}{Lt}}\left(\fint_{{E_{j}^{k}}^{\prime}}M[h^{c^{\prime}\vartheta_{3}}](x)^{\frac{(Lt)^{\prime}}{c^{\prime}\vartheta_{3}}}dx\right)^{\frac{1}{(Lt)^{\prime}}}.
\end{align*}
Since $a\geq L>1$, it follows from H\"{o}lder's inequality for $\frac{a}{L}\geq1$ that
\begin{align*}
\fint_{{Q_{j}^{k}}^{\prime}}M[(v_{j}^{k}h)^{\vartheta_{3}}](x)^{\frac{1}{\vartheta_{3}}}dx\leq C\left(\fint_{{Q_{j}^{k}}^{\prime}}v(x)^{at}dx\right)^{\frac{1}{at}}\left(\fint_{{E_{j}^{k}}^{\prime}}M[h^{c^{\prime}\vartheta_{3}}](x)^{\frac{(Lt)^{\prime}}{c^{\prime}\vartheta_{3}}}dx\right)^{\frac{1}{(Lt)^{\prime}}}.
\end{align*}
By Lemma $\ref{lemma_4}$, it implies that
\begin{align*}
&\mathcal{T}_{j}^{k}\leq2A^{\prime}\int_{{E_{j}^{k}}^{\prime}}M_{\alpha,at}^{r_{\vartheta_{1}},r_{\vartheta_{2}}}(f,g,v)(x)\cdot M[M[h^{c^{\prime}\vartheta_{3}}]^{\frac{(Lt)^{\prime}}{c^{\prime}\vartheta_{3}}}](x)^{\frac{1}{(Lt)^{\prime}}}dx,
\end{align*}
where $r_{\vartheta_{1}}=r_{1}\vartheta_{1}$, $r_{\vartheta_{2}}=r_{2}\vartheta_{2}$ and
\begin{align*}
&M_{\alpha,at}^{r_{\vartheta_{1}},r_{\vartheta_{2}}}(f,g,v)(x)\\
&=\sup_{Q\ni x}|Q|^{\frac{\alpha}{n}}\left(\fint_{3Q}|f(y)|^{\vartheta_{1}r_{1}}dy\right)^{\frac{1}{\vartheta_{1}r_{1}}}\left(\fint_{3Q}|g(z)|^{\vartheta_{2}r_{2}}dz\right)^{\frac{1}{\vartheta_{2}r_{2}}}\left(\fint_{Q}v(x)^{at}dx\right)^{\frac{1}{at}}.
\end{align*}
A similar argument gives us the following estimate
$$\mathcal{T}_{0}\leq2A^{\prime}\int_{{E_{0}}^{\prime}}M_{\alpha,at}^{r_{\vartheta_{1}},r_{\vartheta_{2}}}(f,g,v)(x)\cdot M[M[h^{c^{\prime}\vartheta_{3}}]^{\frac{(Lt)^{\prime}}{c^{\prime}\vartheta_{3}}}](x)^{\frac{1}{(Lt)^{\prime}}}dx.$$
Summing up $\mathcal{T}_{0}$ and $\mathcal{T}_{j}^{k}$, one can obtain
$$\mathcal{T}_{0}+\sum_{k,j}\mathcal{T}_{j}^{k}\leq 2A^{\prime}\int_{Q_{0}}M_{\alpha,at}^{r_{\vartheta_{1}},r_{\vartheta_{2}}}(f,g,v)(x)\cdot M[M[h^{c^{\prime}\vartheta_{3}}]^{\frac{(Lt)^{\prime}}{c^{\prime}\vartheta_{3}}}](x)^{\frac{1}{(Lt)^{\prime}}}dx.$$
By H\"{o}lder's inequality for $t>1$, we have
\begin{align*}
&\int_{Q_{0}}M_{\alpha,at}^{r_{\vartheta_{1}},r_{\vartheta_{2}}}(f,g,v)(x)\cdot M[M[h^{c^{\prime}\vartheta_{3}}]^{\frac{(Lt)^{\prime}}{c^{\prime}\vartheta_{3}}}](x)^{\frac{1}{(Lt)^{\prime}}}dx\\
&\leq\left(\int_{Q_{0}}M_{\alpha,at}^{r_{\vartheta_{1}},r_{\vartheta_{2}}}(f,g,v)(x)^{t}dx\right)^{\frac{1}{t}}\left(\int_{Q_{0}}M[M[h^{c^{\prime}\vartheta_{3}}]^{\frac{(Lt)^{\prime}}{c^{\prime}\vartheta_{3}}}](x)^{\frac{t^{\prime}}{(Lt)^{\prime}}}dx\right)^{\frac{1}{t^{\prime}}}.
\end{align*}
Since $(Lt)^{\prime}<t^{\prime}$, the boundedness of $M:\,L^{\frac{t^{\prime}}{(Lt)^{\prime}}}\rightarrow L^{\frac{t^{\prime}}{(Lt)^{\prime}}}$ tells us that
\begin{align*}
\left(\int_{Q_{0}}M[M[h^{c^{\prime}\vartheta_{3}}]^{\frac{(Lt)^{\prime}}{c^{\prime}\vartheta_{3}}}](x)^{\frac{t^{\prime}}{(Lt)^{\prime}}}dx\right)^{\frac{1}{t^{\prime}}}&\leq C\left(\int_{\mathbb{R}^{n}}M[h^{c^{\prime}\vartheta_{3}}](x)^{\frac{(Lt)^{\prime}}{c^{\prime}\vartheta_{3}}\cdot\frac{t^{\prime}}{(Lt)^{\prime}}}dx\right)^{\frac{1}{t^{\prime}}}\\
&=C\left(\int_{\mathbb{R}^{n}}M[h^{c^{\prime}\vartheta_{3}}](x)^{\frac{t^{\prime}}{c^{\prime}\vartheta_{3}}}dx\right)^{\frac{1}{t^{\prime}}}.
\end{align*}
Noting $c^{\prime}\vartheta_{3}<t^{\prime}$ and the boundedness of $M:\,L^{\frac{t^{\prime}}{c^{\prime}\vartheta_{3}}}\rightarrow L^{\frac{t^{\prime}}{c^{\prime}\vartheta_{3}}}$, we get
\begin{align*}
\left(\int_{\mathbb{R}^{n}}M[h^{c^{\prime}\vartheta_{3}}](x)^{\frac{t^{\prime}}{c^{\prime}\vartheta_{3}}}dx\right)^{\frac{1}{t^{\prime}}}\leq C\left(\int_{Q_{0}}|h(x)|^{c^{\prime}\vartheta_{3}\cdot\frac{t^{\prime}}{c^{\prime}\vartheta_{3}}}dx\right)^{\frac{1}{t^{\prime}}}=C\left(\int_{Q_{0}}|h(x)|^{t^{\prime}}dx\right)^{\frac{1}{t^{\prime}}}\leq C.
\end{align*}
Therefore, we have
\begin{align*}
\int_{Q_{0}}M_{\alpha,at}^{r_{\vartheta_{1}},r_{\vartheta_{2}}}(f,g,v)(x)\cdot M[M[h^{c^{\prime}\vartheta_{3}}]^{\frac{(Lt)^{\prime}}{c^{\prime}\vartheta_{3}}}](x)^{\frac{1}{(Lt)^{\prime}}}dx\leq C\left(\int_{Q_{0}}M_{\alpha,at}^{r_{\vartheta_{1}},r_{\vartheta_{2}}}(f,g,v)(x)^{t}dx\right)^{\frac{1}{t}}.
\end{align*}
Applying H\"{o}lder's inequality for $\frac{q_{1}}{r_{\vartheta_{1}}a_{*}}>1$ and $\frac{q_{2}}{r_{\vartheta_{2}}a_{*}}>1$ as in Lemma $\ref{lemma_4}$, we conclude that
\begin{align*}
&M_{\alpha,at}^{r_{\vartheta_{1}},r_{\vartheta_{2}}}(f,g,v)(x)\leq C\sup_{Q\ni x}|Q|^{\frac{\alpha}{n}}\left(\fint_{3Q}|f(y)w_{1}(y)|^{\frac{q_{1}}{a_{*}}}dy\right)^{\frac{a_{*}}{q_{1}}}\left(\fint_{3Q}|g(z)w_{2}(z)|^{\frac{q_{2}}{a_{*}}}dz\right)^{\frac{a_{*}}{q_{2}}}\times\\
&\left(\fint_{3Q}w_{1}(y)^{-r_{\vartheta_{1}}\big(\frac{q_{1}}{r_{\vartheta_{1}}a_{*}}\big)^{\prime}}dy\right)^{\frac{1}{r_{\vartheta_{1}}\big(\frac{q_{1}}{r_{\vartheta_{1}}a_{*}}\big)^{\prime}}}\left(\fint_{3Q}w_{2}(z)^{-r_{\vartheta_{2}}\big(\frac{q_{2}}{r_{\vartheta_{2}}a_{*}}\big)^{\prime}}dz\right)^{\frac{1}{r_{\vartheta_{2}}\big(\frac{q_{2}}{r_{\vartheta_{2}}a_{*}}\big)^{\prime}}}\left(\fint_{Q}v(x)^{at}dx\right)^{\frac{1}{at}}.
\end{align*}
By Lemma $\ref{lemma_4}$, we have $r_{\vartheta_{1}}\big(\frac{q_{1}}{r_{\vartheta_{1}}a_{*}}\big)^{\prime}\leq\left(\frac{q_{1}}{a}\right)^{\prime}$ and $r_{\vartheta_{2}}\big(\frac{q_{2}}{r_{\vartheta_{2}}a_{*}}\big)^{\prime}\leq\left(\frac{q_{2}}{a}\right)^{\prime}$. If we use H\"{o}lder's inequality and we have
\begin{align*}
&M_{\alpha,at}^{r_{\vartheta_{1}},r_{\vartheta_{2}}}(f,g,v)(x)\leq C\sup_{Q\ni x}|Q|^{\frac{\alpha}{n}}\left(\fint_{3Q}|f(y)w_{1}(y)|^{\frac{q_{1}}{a_{*}}}dy\right)^{\frac{a_{*}}{q_{1}}}\left(\fint_{3Q}|g(z)w_{2}(z)|^{\frac{q_{2}}{a_{*}}}dz\right)^{\frac{a_{*}}{q_{2}}}\left(\frac{|3Q|}{|Q|}\right)^{\frac{1}{as}}\\
&\times|3Q|^{-\frac{1}{r}}\left(\frac{|Q|}{|3Q|}\right)^{\frac{1}{as}}|3Q|^{\frac{1}{r}}\left(\fint_{3Q}w_{1}(y)^{-\left(\frac{q_{1}}{a}\right)^{\prime}}dy\right)^{\frac{1}{\left(\frac{q_{1}}{a}\right)^{\prime}}}\left(\fint_{3Q}w_{2}(z)^{-\left(\frac{q_{2}}{a}\right)^{\prime}}dz\right)^{\frac{1}{\left(\frac{q_{2}}{a}\right)^{\prime}}}\left(\fint_{Q}v(x)^{at}dx\right)^{\frac{1}{at}}.
\end{align*}
By condition $\eqref{two weight condition_2}$, we obtain
$$M_{\alpha,at}^{r_{\vartheta_{1}},r_{\vartheta_{2}}}(f,g,v)(x)\leq C[v,\vec{w}]_{t,\vec{q}/a}^{r,as}M_{\alpha-n/r,\vec{q}/a_{*}}(fw_{1},gw_{2})(x),$$
which implies that
$$|Q_{0}|^{\frac{1}{s}}\left(\fint_{Q_{0}}M_{\alpha,at}^{r_{\vartheta_{1}},r_{\vartheta_{2}}}(f,g,v)(x)^{t}dx\right)^{\frac{1}{t}}\leq C[v,\vec{w}]_{t,\vec{q}/a}^{r,as}\|M_{\alpha-n/r,\vec{q}/a_{*}}(fw_{1},gw_{2})\|_{\mathcal{M}_{t}^{s}}.$$
Since
$$\frac{1}{s}=\frac{1}{p}-\frac{\alpha-n/r}{n}\quad\text{and}\quad\frac{t}{s}=\frac{q}{p},$$
by Lemma $\ref{lemma_5}$, we get
$$\|M_{\alpha-n/r,\vec{q}/a_{*}}(fw_{1},gw_{2})\|_{\mathcal{M}_{t}^{s}}\leq C\sup_{Q\in\mathcal{D}(\mathbb{R}^{n})}|Q|^{1/p}\left(\fint_{Q}(|f|w_{1})^{q_{1}}\right)^{1/q_{1}}\left(\fint_{Q}(|g|w_{2})^{q_{2}}\right)^{1/q_{2}}.$$
Therefore we have
$$\|\mathscr{T}_{1}v\|_{\mathcal{M}_{t}^{s}}\leq C[v,\vec{w}]_{t,\vec{q}/a}^{r,as}\|\vec{b}\|_{{\rm BMO}^{N}}\sup_{Q\in\mathcal{D}(\mathbb{R}^{n})}|Q|^{1/p}\left(\fint_{Q}(|f|w_{1})^{q_{1}}\right)^{1/q_{1}}\left(\fint_{Q}(|g|w_{2})^{q_{2}}\right)^{1/q_{2}}.$$
\quad

{\noindent}{\bf Estimate of $\mathscr{T}_{2}$.} Keeping on assuming that $\eqref{normalization}$ holds.  By the defintion of $\mathscr{T}_{2}$, trangle inequality, the equalities $\eqref{define as average of BMO function}$, $\eqref{part one}$ and $\eqref{part two}$, we obtain
\begin{align*}
\mathscr{T}_{2}(x)&\leq \sum_{A\subseteq\{1,\cdots,m\}}\sum_{B\subseteq\{m+1,\cdots,N\}}\mathop{\sum_{Q\supsetneq Q_{0}}}_{Q\in\mathcal{D}(\mathbb{R}^{n})}|Q|^{\frac{\alpha}{n}-1}\mathop{\int}_{|y|_{\infty}\leq \ell(Q)}\prod_{i\in A\cup B}|b_{i}(x)-\lambda_{i}|\prod_{i\in\bar{A}}|b_{i}(x-y)-\lambda_{i}|\nonumber\\
&\quad\times\prod_{i\in\bar{B}}|b_{i}(x+y)-\lambda_{i}|f(x-y)g(x+y)dy\chi_{Q}(x).
\end{align*}
Using H\"{o}lder's inequality with the pair $(r_{1},r_{2})$ and making a change viariables, we will obtain the following inequality
\begin{align*}
\mathscr{T}_{2}(x)&\leq \sum_{A\subseteq\{1,\cdots,m\}}\sum_{B\subseteq\{m+1,\cdots,N\}}\mathop{\sum_{Q\supsetneq Q_{0}}}_{Q\in\mathcal{D}(\mathbb{R}^{n})}|Q|^{\frac{\alpha}{n}-1}\left(\prod_{i\in A\cup B}|b_{i}(x)-\lambda_{i}|\right)\left[\int_{3Q}\prod_{i\in\bar{A}}(|b_{i}(y)-\lambda_{i}|f(y))^{r_{1}}dy\right]^{\frac{1}{r_{1}}}\nonumber\\
&\quad\times\left[\int_{3Q}\prod_{i\in\bar{B}}(|b_{i}(z)-\lambda_{i}|g(z))^{r_{2}}dz\right]^{\frac{1}{r_{2}}}\chi_{Q}(x).
\end{align*}
For $\vartheta_{1}>1$ and $\vartheta_{2}>1$, using H\"{o}lder's inequality for $\frac{q_{1}}{r_{\vartheta_{1}}}>1$ and $\frac{q_{2}}{r_{\vartheta_{2}}}>1$ as in Lemma $\eqref{lemma_4}$, the inequalities $\eqref{estimate for f}$ and $\eqref{estimate for g}$, and the assumption $\eqref{normalization}$, we have
\begin{align*}
&\mathscr{T}_{2}(x)\leq \sum_{A\subseteq\{1,\cdots,m\}}\sum_{B\subseteq\{m+1,\cdots,N\}}\prod_{i\bar{A}\cup\bar{B}}\|b_{i}\|_{\rm BMO}\mathop{\sum_{Q\supsetneq Q_{0}}}_{Q\in\mathcal{D}(\mathbb{R}^{n})}|Q|^{\frac{\alpha}{n}}\left[\fint_{3Q}|f(y)|^{r_{\vartheta_{1}}}dy\right]^{\frac{1}{r_{\vartheta_{1}}}}\left[\fint_{3Q}|g(z)|^{r_{\vartheta_{2}}}dz\right]^{\frac{1}{r_{\vartheta_{2}}}}\nonumber\\
&\quad\times\left(\prod_{i\in A\cup B}|b_{i}(x)-\lambda_{i}|\right)\chi_{Q}(x)\\
&\leq\sum_{A\subseteq\{1,\cdots,m\}}\sum_{B\subseteq\{m+1,\cdots,N\}}\prod_{i\bar{A}\cup\bar{B}}\|b_{i}\|_{\rm BMO}\mathop{\sum_{Q\supsetneq Q_{0}}}_{Q\in\mathcal{D}(\mathbb{R}^{n})}\chi_{Q}(x)|Q|^{\frac{\alpha}{n}}\left[\fint_{3Q}|f(y)w_{1}(y)|^{q_{1}}dy\right]^{\frac{1}{q_{1}}}\left[\fint_{3Q}|g(z)w_{2}(z)|^{q_{2}}dz\right]^{\frac{1}{q_{2}}}\nonumber\\
&\quad\times\left(\prod_{i\in A\cup B}|b_{i}(x)-\lambda_{i}|\right)\left(\fint_{3Q}w_{1}(y)^{-r_{\vartheta_{1}}\big(\frac{q_{1}}{r_{\vartheta_{1}}}\big)^{\prime}}dy\right)^{\frac{1}{r_{\vartheta_{1}}\big(\frac{q_{1}}{r_{\vartheta_{1}}}\big)^{\prime}}}\left(\fint_{3Q}w_{2}(z)^{-r_{\vartheta_{2}}\big(\frac{q_{2}}{r_{\vartheta_{2}}}\big)^{\prime}}dz\right)^{\frac{1}{r_{\vartheta_{1}}\big(\frac{q_{2}}{r_{\vartheta_{1}}}\big)^{\prime}}}\\
&\leq C\sum_{A\subseteq\{1,\cdots,m\}}\sum_{B\subseteq\{m+1,\cdots,N\}}\prod_{i\bar{A}\cup\bar{B}}\|b_{i}\|_{\rm BMO}\mathop{\sum_{Q\supsetneq Q_{0}}}_{Q\in\mathcal{D}(\mathbb{R}^{n})}|Q|^{\frac{\alpha}{n}}|3Q|^{-\frac{1}{p}}\left(\prod_{i\in A\cup B}|b_{i}(x)-\lambda_{i}|\right)\chi_{Q}(x)\\
&\qquad\qquad\times\left(\fint_{3Q}w_{1}(y)^{-r_{\vartheta_{1}}\big(\frac{q_{1}}{r_{\vartheta_{1}}}\big)^{\prime}}dy\right)^{\frac{1}{r_{\vartheta_{1}}\big(\frac{q_{1}}{r_{\vartheta_{1}}}\big)^{\prime}}}\left(\fint_{3Q}w_{2}(z)^{-r_{\vartheta_{2}}\big(\frac{q_{2}}{r_{\vartheta_{2}}}\big)^{\prime}}dz\right)^{\frac{1}{r_{\vartheta_{2}}\big(\frac{q_{2}}{r_{\vartheta_{2}}}\big)^{\prime}}}.
\end{align*}
By Minkowski's inequality for $t>1$, we obain
\begin{align*}
&\left(\fint_{Q_{0}}(\mathscr{T}_{2}(x)v(x))^{t}dx\right)^{\frac{1}{t}}\\
&\leq C\sum_{A\subseteq\{1,\cdots,m\}}\sum_{B\subseteq\{m+1,\cdots,N\}}\prod_{i\bar{A}\cup\bar{B}}\|b_{i}\|_{\rm BMO}\mathop{\sum_{Q\supsetneq Q_{0}}}_{Q\in\mathcal{D}(\mathbb{R}^{n})}|3Q|^{\frac{\alpha}{n}-\frac{1}{p}}\left(\fint_{Q_{0}}\prod_{i\in A\cup B}|b_{i}(x)-\lambda_{i}|^{t}v(x)^{t}dx\right)^{\frac{1}{t}}\\
&\qquad\qquad\qquad\qquad\quad\times\left(\fint_{3Q}w_{1}(y)^{-r_{\vartheta_{1}}\big(\frac{q_{1}}{r_{\vartheta_{1}}}\big)^{\prime}}dy\right)^{\frac{1}{r_{\vartheta_{1}}\big(\frac{q_{1}}{r_{\vartheta_{1}}}\big)^{\prime}}}\left(\fint_{3Q}w_{2}(z)^{-r_{\vartheta_{2}}\big(\frac{q_{2}}{r_{\vartheta_{2}}}\big)^{\prime}}dz\right)^{\frac{1}{r_{\vartheta_{1}}\big(\frac{q_{2}}{r_{\vartheta_{1}}}\big)^{\prime}}}\\
&\leq C\sum_{A\subseteq\{1,\cdots,m\}}\sum_{B\subseteq\{m+1,\cdots,N\}}\prod_{i\bar{A}\cup\bar{B}}\|b_{i}\|_{\rm BMO}\sum_{k=1}^{\infty}\mathop{\sum_{Q_{k}\in\mathcal{D}(\mathbb{R}^{n})}}_{Q_{k}\supseteq Q_{0},\,|Q_{k}|=2^{kn}|Q_{0}|}\left(\fint_{Q_{0}}\prod_{i\in A\cup B}|b_{i}(x)-\lambda_{i}|^{t}v(x)^{t}dx\right)^{\frac{1}{t}}\\
&\qquad\qquad\times|3Q_{k}|^{\frac{\alpha}{n}-\frac{1}{p}}\left(\fint_{3Q_{k}}w_{1}(y)^{-r_{\vartheta_{1}}\big(\frac{q_{1}}{r_{\vartheta_{1}}}\big)^{\prime}}dy\right)^{\frac{1}{r_{\vartheta_{1}}\big(\frac{q_{1}}{r_{\vartheta_{1}}}\big)^{\prime}}}\left(\fint_{3Q_{k}}w_{2}(z)^{-r_{\vartheta_{2}}\big(\frac{q_{2}}{r_{\vartheta_{2}}}\big)^{\prime}}dz\right)^{\frac{1}{r_{\vartheta_{2}}\big(\frac{q_{2}}{r_{\vartheta_{2}}}\big)^{\prime}}}.
\end{align*}
For $a\geq\vartheta_{3}>1$, similar to the estimate $\eqref{v controlled by BMO norm}$, we have
\begin{equation}\label{v controlled by BMO norm for t>1}
\left(\fint_{Q_{0}}\prod_{i\in A\cup B}|b_{i}(x)-\lambda_{i}|^{t}v(x)^{t}dx\right)^{\frac{1}{t}}\leq Ck\prod_{i\in A\cup B}\|b_{i}\|_{\rm BMO}\left(\fint_{Q_{0}}v(x)^{at}dx\right)^{\frac{1}{at}}.
\end{equation}
Again using H\"{o}lder's inequality for $r_{\vartheta_{1}}\big(\frac{q_{1}}{r_{\vartheta_{1}}}\big)^{\prime}\leq\left(\frac{q_{1}}{a}\right)^{\prime}$ and $r_{\vartheta_{2}}\big(\frac{q_{2}}{r_{\vartheta_{2}}}\big)^{\prime}\leq\left(\frac{q_{2}}{a}\right)^{\prime}$ as in Lemma $\ref{lemma_4}$, one has
\begin{equation}\label{estimate for w_1}
\left(\fint_{3Q_{k}}w_{1}(y)^{-r_{\vartheta_{1}}\big(\frac{q_{1}}{r_{\vartheta_{1}}}\big)^{\prime}}dy\right)^{\frac{1}{r_{\vartheta_{1}}\big(\frac{q_{1}}{r_{\vartheta_{1}}}\big)^{\prime}}}\leq\left(\int_{3Q_{k}}w_{1}(y)^{-\left(\frac{q_{1}}{a}\right)^{\prime}}dy\right)^{\frac{1}{\left(\frac{q_{1}}{a}\right)^{\prime}}}
\end{equation}
and
\begin{equation}\label{estimate for w_2}
\left(\fint_{3Q_{k}}w_{2}(z)^{-r_{\vartheta_{2}}\big(\frac{q_{2}}{r_{\vartheta_{2}}}\big)^{\prime}}dz\right)^{\frac{1}{r_{\vartheta_{2}}\big(\frac{q_{2}}{r_{\vartheta_{2}}}\big)^{\prime}}}\leq\left(\int_{3Q_{k}}w_{2}(z)^{-\left(\frac{q_{2}}{a}\right)^{\prime}}dz\right)^{\frac{1}{\left(\frac{q_{2}}{a}\right)^{\prime}}}.
\end{equation}
Combining $\eqref{v controlled by BMO norm for t>1}$ with $\eqref{estimate for w_1}$ and $\eqref{estimate for w_2}$, it yields that
\begin{align*}
&|Q_{0}|^{\frac{1}{s}}\left(\fint_{Q_{0}}(\mathscr{T}_{2}(x)v(x))^{t}dx\right)^{\frac{1}{t}}\\
&\leq C\|\vec{b}\|_{{\rm BMO}^{N}}\sum_{k=1}^{\infty}k\mathop{\sum_{Q_{k}\in\mathcal{D}(\mathbb{R}^{n})}}_{Q_{k}\supseteq Q_{0},\,|Q_{k}|=2^{kn}|Q_{0}|}|Q_{0}|^{\frac{1}{s}}|3Q_{k}|^{\frac{\alpha}{n}-\frac{1}{p}}\\
&\qquad\times\left(\fint_{Q_{0}}v(x)^{at}dx\right)^{\frac{1}{at}}\left(\int_{3Q_{k}}w_{1}(y)^{-\left(\frac{q_{1}}{a}\right)^{\prime}}dy\right)^{\frac{1}{\left(\frac{q_{1}}{a}\right)^{\prime}}}\left(\int_{3Q_{k}}w_{2}(z)^{-\left(\frac{q_{2}}{a}\right)^{\prime}}dz\right)^{\frac{1}{\left(\frac{q_{2}}{a}\right)^{\prime}}}\\
&= C\|\vec{b}\|_{{\rm BMO}^{N}}\sum_{k=1}^{\infty}k\mathop{\sum_{Q_{k}\in\mathcal{D}(\mathbb{R}^{n})}}_{Q_{k}\supseteq Q_{0},\,|Q_{k}|=2^{kn}|Q_{0}|}|Q_{0}|^{\frac{1}{s}}|3Q_{k}|^{\frac{1}{r}-\frac{1}{s}}\left(\frac{|3Q_{k}|}{|Q_{0}|}\right)^{\frac{1}{as}}\left(\frac{|Q_{0}|}{|3Q_{k}|}\right)^{\frac{1}{as}}\\
&\qquad\times\left(\fint_{Q_{0}}v(x)^{at}dx\right)^{\frac{1}{at}}\left(\int_{3Q_{k}}w_{1}(y)^{-\left(\frac{q_{1}}{a}\right)^{\prime}}dy\right)^{\frac{1}{\left(\frac{q_{1}}{a}\right)^{\prime}}}\left(\int_{3Q_{k}}w_{2}(z)^{-\left(\frac{q_{2}}{a}\right)^{\prime}}dz\right)^{\frac{1}{\left(\frac{q_{2}}{a}\right)^{\prime}}}.
\end{align*}
By condition $\eqref{two weight condition_2}$, we get
\begin{align*}
&|Q_{0}|^{\frac{1}{s}}\left(\fint_{Q_{0}}(\mathscr{T}_{2}(x)v(x))^{t}dx\right)^{\frac{1}{t}}\\
&\leq C[v,\vec{w}]_{t,\vec{q}/a}^{r,as}\|\vec{b}\|_{{\rm BMO}^{N}}\sum_{k=1}^{\infty}k\mathop{\sum_{Q_{k}\in\mathcal{D}(\mathbb{R}^{n})}}_{Q_{k}\supseteq Q_{0},\,|Q_{k}|=2^{kn}|Q_{0}|}\left(\frac{|Q_{0}|}{|3Q_{k}|}\right)^{(1-\frac{1}{a})\frac{1}{s}}\\
&\leq C[v,\vec{w}]_{t,\vec{q}/a}^{r,as}\|\vec{b}\|_{{\rm BMO}^{N}}\sum_{k=1}^{\infty}k\mathop{\sum_{Q_{k}\in\mathcal{D}(\mathbb{R}^{n})}}_{Q_{k}\supseteq Q_{0},\,|Q_{k}|=2^{kn}|Q_{0}|}2^{-\frac{kn}{s}(1-\frac{1}{a})}\\
&\leq C[v,\vec{w}]_{t,\vec{q}/a}^{r,as}\|\vec{b}\|_{{\rm BMO}^{N}}.
\end{align*}
Therefore, we have
$$\|\mathscr{T}_{2}v\|_{\mathcal{M}_{t}^{s}}\leq C[v,\vec{w}]_{t,\vec{q}/a}^{r,as}\|\vec{b}\|_{{\rm BMO}^{N}}\sup_{Q\in\mathcal{D}(\mathbb{R}^{n})}|Q|^{1/p}\left(\fint_{Q}(|f|w_{1})^{q_{1}}\right)^{1/q_{1}}\left(\fint_{Q}(|g|w_{2})^{q_{2}}\right)^{1/q_{2}}.$$
We obtain the desired result. $\hfill$ $\Box$

\section{The proof of Theorem $\ref{main_5}$}

Again, without loss of generality, we may restrict ourselves onto working with $f$ and $g$ that are non-negative, bounded and compactly supported. Thanks to Theorem $\ref{lemma_8}$, we only need to verify the inequality for a certain $q_{0}\in(0,\infty)$ and an arbitrary weight
$w\in A_{\infty}$. We will work with $q_{0}=1$. By mimicking what we did in the proof of Theorem $\ref{main_3}$, we have
\begin{align}\label{controlled by maximal operator with the weight w}
&\int_{\mathbb{R}^{n}}|[\vec{b},B_{\alpha}]_{\vec{\beta}}(f,g)(x)|w(x)dx\\
&\lesssim\|\vec{b}\|_{{\rm BMO}^{N}}\sum_{k,j}|Q_{j}^{k}|^{\frac{\alpha}{n}+1}\left(\fint_{3Q_{j}^{k}}|f(y)|^{r_{\vartheta_{1}}}dy\right)^{\frac{1}{r_{\vartheta_{1}}}}\left(\fint_{3Q_{j}^{k}}|g(z)|^{r_{\vartheta_{2}}}dz\right)^{\frac{1}{r_{\vartheta_{2}}}}\left(\fint_{Q_{j}^{k}}w(x)^{\vartheta_{3}}dx\right)^{\frac{1}{\vartheta_{3}}}\nonumber.
\end{align}
Since $w\in A_{\infty}$, there exist, by Lemma $\ref{lemma_7}$, a index $\vartheta_{3}>1$ such that
$$\left(\fint_{Q_{j}^{k}}w(x)^{\vartheta_{3}}dx\right)^{\frac{1}{\vartheta_{3}}}\leq C\fint_{Q_{j}^{k}}w(x)dx.$$
Substituting this result into $\eqref{controlled by maximal operator with the weight w}$, we have
\begin{align*}
&\int_{\mathbb{R}^{n}}|[\vec{b},B_{\alpha}]_{\vec{\beta}}(f,g)(x)|w(x)dx\\
&\lesssim \|\vec{b}\|_{{\rm BMO}^{N}}\sum_{k,j}|Q_{j}^{k}|^{\frac{\alpha}{n}}\left(\fint_{3Q_{j}^{k}}|f(y)|^{r_{\vartheta_{1}}}dy\right)^{\frac{1}{r_{\vartheta_{1}}}}\left(\fint_{3Q_{j}^{k}}|g(z)|^{r_{\vartheta_{2}}}dz\right)^{\frac{1}{r_{\vartheta_{2}}}}w(Q_{j}^{k})\\
&\lesssim \|\vec{b}\|_{{\rm BMO}^{N}}\sum_{k,j}|Q_{j}^{k}|^{\frac{\alpha}{n}}\left(\fint_{3Q_{j}^{k}}|f(y)|^{r_{\vartheta_{1}}}dy\right)^{\frac{1}{r_{\vartheta_{1}}}}\left(\fint_{3Q_{j}^{k}}|g(z)|^{r_{\vartheta_{2}}}dz\right)^{\frac{1}{r_{\vartheta_{2}}}}w(E_{j}^{k})\\
&\lesssim\|\vec{b}\|_{{\rm BMO}^{N}}\sum_{k,j}\int_{E_{j}^{k}}M_{\alpha,\vec{R_{\vartheta}}}(f,g)(x)w(x)dx\\
&\leq\|\vec{b}\|_{{\rm BMO}^{N}}\int_{\mathbb{R}^{n}}M_{\alpha,\vec{R_{\vartheta}}}(f,g)(x)w(x)dx,
\end{align*}
where $\vec{R_{\vartheta}}=(r_{\vartheta_{1}},r_{\vartheta_{2}})$, and the second inequality is due to Lemma $\ref{lemma_7}$ and the fact that $w\in A_{\infty}$. $\hfill$ $\Box$

\section{The proof of Theorem $\ref{main_6}$}

$[Condition~\eqref{two weight condition_3}~\Rightarrow~the~weak ~type ~boundedness]$

For a fixed cube $Q_{0}\in\mathcal{D}(\mathbb{R}^{n})$, let $(f^{0},g^{0})=(\chi_{3Q_{0}}f,\chi_{3Q_{0}}g)$. A normalization allows us to assume that
\begin{equation}\label{normalization_1}
\mathop{\sup_{Q\in\mathcal{D}(\mathbb{R}^{n})}}_{Q\supseteq Q_{0}}|Q|^{1/p}\left(\fint_{Q}(|f(y)|w_{1}(y))^{q_{1}}dy\right)^{1/q_{1}}\left(\fint_{Q}(|g(z)|w_{2}(z))^{q_{2}}dz\right)^{1/q_{2}}=1.
\end{equation}
Then by a standard argument we have, for $x\in Q_{0}$,
\begin{equation}\label{maxinal operator controlled by two parts}
M_{\alpha,\vec{R}}(f,g)(x)\leq M_{\alpha,\vec{R}}(f^{0},g^{0})(x)+Cc_{\infty},
\end{equation}
where
\begin{equation}\label{define c}
c_{\infty}=\mathop{\sup_{Q\in\mathcal{D}(\mathbb{R}^{n})}}_{Q\supset Q_{0}}|Q|^{\frac{\alpha}{n}}\left(\fint_{Q}|f(y)|^{r_{1}}dy\right)^{\frac{1}{r_{1}}}\left(\fint_{Q}|g(z)|^{r_{2}}dz\right)^{\frac{1}{r_{2}}}.
\end{equation}
For $x\in Q_{0}$, we notice that the following pointwise equivalence holds:
\begin{equation}\label{equivalence formula}
M_{\alpha,\vec{R}}(f^{0},g^{0})(x)\approx\hat{M}_{\alpha,\vec{R}}(f,g)(x),
\end{equation}
where
$$\hat{M}_{\alpha,\vec{R}}(f,g)(x)=\sup_{\mathcal{D}(Q_{0})\ni {Q}\ni x}|{Q}|^{\frac{\alpha}{n}}\left(\fint_{{Q}}|f(y)|^{r_{1}}dy\right)^{\frac{1}{r_{1}}}\left(\fint_{{Q}}|g(z)|^{r_{2}}dz\right)^{\frac{1}{r_{2}}}.$$
Combining $\eqref{maxinal operator controlled by two parts}$ with $\eqref{equivalence formula}$, it yields that
\begin{equation*}
\{x\in Q_{0}:\,M_{\alpha,\vec{R}}(f,g)(x)>\lambda\}\subseteq\{x\in Q_{0}:\,\hat{M}_{\alpha,\vec{R}}(f,g)(x)>\lambda/2\}\bigcup\{x\in Q_{0}:\,Cc_{\infty}>\lambda/2\}.
\end{equation*}
Define
\begin{eqnarray}\label{define E_1 and E_2}
E_{\lambda}^{1}:=\{x\in Q_{0}:\,\hat{M}_{\alpha,\vec{R}}(f,g)(x)>\lambda/2\}\quad\text{and}\quad E_{\lambda}^{2}:=\{x\in Q_{0}:\,Cc_{\infty}>\lambda/2\}.
\end{eqnarray}
By performing the Calder\'{o}n-Zygmund decomposition algorithm, we have
\begin{equation}\label{decomposition E_{1}}
E_{\lambda}^{1}=\bigcup_{j}Q_{j},
\end{equation}
where $Q_{j}$'s are pairwise disjoint maximal dyadic cubes that satisfy
\begin{equation}\label{maximal dyadic cubes}
|Q_{j}|^{\frac{\alpha}{n}}\left(\fint_{Q_{j}}|f(y)|^{r_{1}}dy\right)^{\frac{1}{r_{1}}}\left(\fint_{Q_{j}}|g(z)|^{r_{2}}dz\right)^{\frac{1}{r_{2}}}>\lambda/2.
\end{equation}
Using the usual notation for weight, we write
\begin{equation}\label{notation for weight}
v^{t}(E)=\int_{E}v(x)^{t}dx.
\end{equation}
From $\eqref{decomposition E_{1}}$ and $\eqref{maximal dyadic cubes}$ we have
\begin{align*}
v^{t}(E_{\lambda}^{1})&=\sum_{j}\int_{Q_{j}}v(x)^{t}dx=\sum_{j}|Q_{j}|\fint_{Q_{j}}v(x)^{t}dx\\
&\leq\frac{2^{t}}{\lambda^{t}}\sum_{j}|Q_{j}|^{\frac{\alpha t}{n}+1}\left(\fint_{Q_{j}}v(x)^{t}dx\right)\left(\fint_{Q_{j}}|f(y)|^{r_{1}}dy\right)^{\frac{t}{r_{1}}}\left(\fint_{Q_{j}}|g(z)|^{r_{2}}dz\right)^{\frac{t}{r_{2}}}.
\end{align*}
Using H\"{o}lder's inequality for  the second and the third dashed integrals, we obtain
\begin{align*}
v^{t}(E_{\lambda}^{1})&\leq\frac{2^{t}}{\lambda^{t}}\sum_{j}|Q_{j}|^{\frac{\alpha t}{n}+1}\left(\fint_{Q_{j}}v(x)^{t}dx\right)\prod_{i=1}^{2}\left(\fint_{Q_{j}}w_{i}(y)^{-r_{i}\big(\frac{q_{i}}{r_{i}}\big)^{\prime}}dy\right)^{\frac{t}{r_{i}\big(\frac{q_{i}}{r_{i}}\big)^{\prime}}}\\
&\qquad\qquad\qquad\qquad\times\left(\fint_{Q_{j}}(|f(y)|w_{1}(y))^{q_{1}}dy\right)^{\frac{t}{q_{1}}}\left(\fint_{Q_{j}}(|g(z)|w_{2}(z))^{q_{2}}dz\right)^{\frac{t}{q_{2}}}.
\end{align*}
By condition $\eqref{two weight condition_3}$, we have
$$
c_{1}=\sup_{Q\in\mathcal{D}(\mathbb{R}^{n})}|Q|^{\frac{1}{r}}\left(\fint_{Q}v(x)^{t}dx\right)^{\frac{1}{t}}\prod_{i=1}^{2}\left(\fint_{Q}w_{i}(y)^{-r_{i}\big(\frac{q_{i}}{r_{i}}\big)^{\prime}}dy\right)^{\frac{1}{r_{i}\big(\frac{q_{i}}{r_{i}}\big)^{\prime}}}
\leq[v,\vec{w}]_{t,\vec{q}}^{r,s}.$$
Therefore,
\begin{align*}
v^{t}(E_{\lambda}^{1})&\leq C \frac{c_{1}}{\lambda^{t}}\sum_{j}|Q_{j}|^{\frac{t}{p}-\frac{t}{s}+1-\frac{t}{q}}\left(\int_{Q_{j}}(|f(y)|w_{1}(y))^{q_{1}}dy\right)^{\frac{t}{q_{1}}}\left(\int_{Q_{j}}(|g(z)|w_{2}(z))^{q_{2}}dz\right)^{\frac{t}{q_{2}}}dx\\
&\leq C \frac{c_{1}}{\lambda^{t}}\sum_{j}\left(\int_{Q_{j}}(|f(y)|w_{1}(y))^{q_{1}}dy\right)^{\frac{t}{q_{1}}}\left(\int_{Q_{j}}(|g(z)|w_{2}(z))^{q_{2}}dz\right)^{\frac{t}{q_{2}}}\\
&\leq C\frac{c_{1}}{\lambda^{t}}\left[\sum_{j}\left(\int_{Q_{j}}(|f(y)|w_{1}(y))^{q_{1}}dy\right)^{\frac{q}{q_{1}}}\left(\int_{Q_{j}}(|g(z)|w_{2}(z))^{q_{2}}dz\right)^{\frac{q}{q_{2}}}dx\right]^{\frac{t}{q}}\\
&\leq C\frac{c_{1}}{\lambda^{t}}\left(\sum_{j}\int_{Q_{j}}(|f(y)|w_{1}(y))^{q_{1}}dy\right)^{\frac{t}{q_{1}}}\left(\sum_{j}\int_{Q_{j}}(|g(z)|w_{2}(z))^{q_{2}}dz\right)^{\frac{t}{q_{2}}}\\
&\leq C\frac{c_{1}}{\lambda^{t}}\left(\int_{Q_{0}}(|f(y)|w_{1}(y))^{q_{1}}dy\right)^{\frac{t}{q_{1}}}\left(\int_{Q_{0}}(|g(z)|w_{2}(z))^{q_{2}}dz\right)^{\frac{t}{q_{2}}},
\end{align*}
where the condition $\frac{1}{p}-\frac{1}{s}=\frac{1}{q}-\frac{1}{t}$ is used in the second step. Noting the assumption $\eqref{normalization_1}$, we have
\begin{align}\label{estimate for E_1}
&|Q_{0}|^{\frac{1}{s}-\frac{1}{t}}{\lambda}\left(v^{t}(E_{\lambda}^{1})\right)^{\frac{1}{t}}\nonumber\\
&\leq C c_{1}|Q_{0}|^{\frac{1}{p}-\frac{1}{q}}\left(\int_{Q_{0}}(|f(y)|w_{1}(y))^{q_{1}}dy\right)^{\frac{1}{q_{1}}}\left(\int_{Q_{0}}(|g(z)|w_{2}(z))^{q_{2}}dz\right)^{\frac{1}{q_{2}}}\leq C [v,\vec{w}]_{t,\vec{q}}^{r,s}.
\end{align}
Next we estimate for $ E_{\lambda}^{2}$. In the case $\lambda>0$, by Chebyshev's inequality, we have
\begin{align*}
\sup_{\lambda>0}|Q_{0}|^{\frac{1}{s}-\frac{1}{t}}\frac{1}{\lambda}(v^{t}( E_{\lambda}^{2}))^{\frac{1}{t}}&=\sup_{\lambda>0}|Q_{0}|^{\frac{1}{s}-\frac{1}{t}}\frac{1}{\lambda}\left(\int_{ E_{\lambda}^{2}}v(x)^{t}dx\right)^{\frac{1}{t}}\\
&\leq C|Q_{0}|^{\frac{1}{s}}\left(\fint_{Q_{0}}v(x)^{t}dx\right)^{\frac{1}{t}}c_{\infty}=:CK.
\end{align*}
By the definition of $c_{\infty}$ and the assumption $\eqref{normalization_1}$ , we can estimate
\begin{align}\label{estimate for E_2}
K&\leq |Q_{0}|^{\frac{1}{s}}\left(\fint_{Q_{0}}v(x)^{t}dx\right)^{\frac{1}{t}}\mathop{\sup_{Q\in\mathcal{D}(\mathbb{R}^{n})}}_{Q\supset Q_{0}}|Q|^{\frac{\alpha}{n}}\left(\fint_{Q}(|f(y)|w_{1}(y)^{q_{1}}dy\right)^{\frac{1}{q_{1}}}\left(\fint_{Q}(|g(z)|w_{2}(z)^{q_{2}}dz\right)^{\frac{1}{q_{2}}}\nonumber\\
&\qquad\qquad\qquad\qquad\qquad\times\left(\fint_{Q}w_{1}(y)^{-r_{1}\big(\frac{q_{1}}{r_{1}}\big)^{\prime}}dy\right)^{\frac{1}{r_{1}\big(\frac{q_{1}}{r_{1}}\big)^{\prime}}}\left(\fint_{Q}w_{2}(z)^{-r_{2}\big(\frac{q_{2}}{r_{2}}\big)^{\prime}}dz\right)^{\frac{1}{r_{2}\big(\frac{q_{2}}{r_{2}}\big)^{\prime}}}\nonumber\\
&\leq\mathop{\sup_{Q\in\mathcal{D}(\mathbb{R}^{n})}}_{Q\supset Q_{0}}\left(\fint_{Q}v(x)^{t}dx\right)^{\frac{1}{t}}\left(\fint_{Q}w_{1}(y)^{-r_{1}\left(\frac{q_{1}}{r_{1}}\right)^{\prime}}dy\right)^{\frac{1}{r_{1}\big(\frac{q_{1}}{r_{1}}\big)^{\prime}}}\left(\fint_{Q}w_{2}(z)^{-r_{2}\big(\frac{q_{2}}{r_{2}}\big)^{\prime}}dz\right)^{\frac{1}{r_{2}\big(\frac{q_{2}}{r_{2}}\big)^{\prime}}}\nonumber\\
&\quad\quad\times\left(\frac{|Q_{0}|}{|Q|}\right)^{\frac{1}{s}}|Q|^{\frac{1}{r}}\leq C[v,\vec{w}]_{t,\vec{q}}^{r,s}.
\end{align}
Combining $\eqref{maxinal operator controlled by two parts}$, $\eqref{define c}$, $\eqref{equivalence formula}$, $\eqref{define E_1 and E_2}$, $\eqref{estimate for E_1}$ with $\eqref{estimate for E_2}$ and Minkowski's inequality, we conclude that for every $Q_{0}\in\mathcal{D}(\mathbb{R}^{n})$,
\begin{align*}
&\sup_{\lambda>0}|Q_{0}|^{\frac{1}{s}-\frac{1}{t}}\lambda (v^{t}(\{x\in Q_{0}:\,|M_{\alpha,\vec{R}}(f,g)(x)|>\lambda\}))^{\frac{1}{t}}\\
&\qquad\qquad\leq C[v,\vec{w}]_{t,\vec{q}/a}^{r,s} \mathop{\sup_{Q\in\mathcal{D}(\mathbb{R}^{n})}}_{Q\supseteq Q_{0}}|Q^{1/p}\left(\fint_{Q}(|f|w_{1})^{q_{1}}\right)^{1/q_{1}}\left(\fint_{Q}(|g|w_{2})^{q_{2}}\right)^{1/q_{2}}.
\end{align*}

\quad

{\noindent}[$The ~weak ~boundedness~\Rightarrow ~condition~ \eqref{two weight condition_3}$]

We assume to the contrary that
\begin{equation}\label{the reverse assumption_1}
[v,\vec{w}]_{t,\vec{q}/a}^{r,s}=\mathop{\sup_{Q,Q^{\prime}\in\mathcal{D}(\mathbb{R}^{n})}}_{Q\subseteq Q^{\prime}}\left(\frac{|Q|}{|Q^{\prime}|}\right)^{\frac{1}{s}}|Q^{\prime}|^{\frac{1}{r}}\left(\fint_{Q}v(x)^{t}dx\right)^{\frac{1}{t}}\prod_{i=1}^{2}\left(\fint_{Q^{\prime}}w_{i}(y)^{-r_{i}\left(\frac{q_{i}}{r_{i}}\right)^{\prime}}dy\right)^{\frac{1}{r_{i}\left(\frac{q_{i}}{r_{i}}\right)^{\prime}}}=\infty.
\end{equation}
By $\eqref{the reverse assumption_1}$ we can select two cubes $Q\supseteq Q^{\prime}$ such that for a large $N$,
\begin{equation}\label{the reverse assumption_2}
\left(\frac{|Q|}{|Q^{\prime}|}\right)^{\frac{1}{s}}|Q^{\prime}|^{\frac{1}{r}}\left(\fint_{Q}v(x)^{t}dx\right)^{\frac{1}{t}}\prod_{i=1}^{2}\left(\fint_{Q^{\prime}}w_{i}(y)^{-r_{i}\left(\frac{q_{i}}{r_{i}}\right)^{\prime}}dy\right)^{\frac{1}{r_{i}\left(\frac{q_{i}}{r_{i}}\right)^{\prime}}}>N.
\end{equation}
Select a smaller cube $Q^{\prime}$ in $\eqref{the reverse assumption_2}$ (if necessary). Without loss of generality, we may assume that $Q^{\prime}$ is minimal in the sense that
\begin{equation}\label{the condition assumption}
\mathop{\sup_{R\in\mathcal{D}(\mathbb{R}^{n})}}_{Q\subseteq R\subseteq Q^{\prime}}\fint_{R}w_{i}(y)^{-l_{i}\big(\frac{q_{i}}{r_{i}}\big)^{\prime}}dy=\fint_{Q^{\prime}}w_{i}(y)^{-r_{i}\big(\frac{q_{i}}{r_{i}}\big)^{\prime}}dy\quad \text{for}\quad i=1,2.
\end{equation}
By noticing $\frac{1}{p}-\frac{1}{q}=\frac{1}{p}-\frac{1}{q_{1}}-\frac{1}{q_{2}}\leq 0$, the equality $\eqref{the condition assumption}$ yields
\begin{equation}\label{the condition conclude}
\mathop{\sup_{R\in\mathcal{D}(\mathbb{R}^{n})}}_{R\supseteq Q}|R|^{\frac{1}{p}}\prod_{i=1}^{2}\left(\fint_{R}\chi_{Q^{\prime}}w_{i}(y)^{-r_{i}\big(\frac{q_{i}}{r_{i}}\big)^{\prime}}dy\right)^{\frac{1}{q_{i}}}=|Q^{\prime}|^{\frac{1}{p}}\prod_{i=1}^{2}\left(\fint_{Q^{\prime}}w_{i}(y)^{-r_{i}\big(\frac{q_{i}}{r_{i}}\big)^{\prime}}dy\right)^{\frac{1}{q_{i}}}.
\end{equation}
For the cube $Q^{\prime}$, let
$$
f=\chi_{Q^{\prime}}w_{1}^{-\frac{q_{1}}{q_{1}-r_{1}}}\quad\text{and}\quad g=\chi_{Q^{\prime}}w_{2}^{-\frac{q_{2}}{q_{2}-r_{2}}}.
$$
Then by choosing
$$\lambda=\frac{1}{2}|Q^{\prime}|^{\frac{\alpha}{n}}\left(\fint_{Q^{\prime}}f(y)^{r_{1}}dy\right)^{\frac{1}{r_{1}}}\left(\fint_{Q^{\prime}}g(z)^{r_{2}}dz\right)^{\frac{1}{r_{2}}},$$
from the weak type boundedness of $\eqref{the weak type boundedness}$ we have
\begin{align*}
&|Q|^{\frac{1}{s}-\frac{1}{t}}\left(\int_{Q}v(x)^{t}dx\right)^{\frac{1}{t}}|Q^{\prime}|^{\frac{\alpha}{n}}\left(\fint_{Q^{\prime}}f(y)^{r_{1}}dy\right)^{\frac{1}{r_{1}}}\left(\fint_{Q^{\prime}}g(z)^{r_{2}}dz\right)^{\frac{1}{r_{2}}}\\
&\qquad\qquad\leq 2C\mathop{\sup_{R\in\mathcal{D}(\mathbb{R}^{n})}}_{R\supseteq Q}|R|^{\frac{1}{p}}\left(\fint_{R}(|f(y)|w_{1}(y))^{q_{1}}dy\right)^{\frac{1}{q_{1}}}\left(\fint_{R}(|g(z)|w_{2}(z))^{q_{2}}dz\right)^{\frac{1}{q_{2}}}.
\end{align*}
Now, if we substitute our specific choices of $f$ and $g$ into the expression, we have
\begin{align*}
&|Q|^{\frac{1}{s}}\left(\fint_{Q}v(x)^{t}dx\right)^{\frac{1}{t}}|Q^{\prime}|^{\frac{\alpha}{n}}\left(\fint_{Q^{\prime}}w_{1}(y)^{-\frac{q_{1}r_{1}}{q_{1}-r_{1}}}dy\right)^{\frac{1}{r_{1}}}\left(\fint_{Q^{\prime}}w_{2}(z)^{-\frac{q_{2}r_{2}}{q_{2}-r_{2}}}dz\right)^{\frac{1}{r_{2}}}\\
&\leq 2C\mathop{\sup_{R\in\mathscr{D}}}_{R\supseteq Q}|R|^{\frac{1}{p}}\left(\fint_{R}\chi_{Q^{\prime}}w_{1}(y)^{-\frac{q_{1}r_{1}}{q_{1}-r_{1}}}dy\right)^{\frac{1}{q_{1}}}\left(\fint_{R}\chi_{Q^{\prime}}w_{2}(z)^{-\frac{q_{2}r_{2}}{q_{2}-r_{2}}}dz\right)^{\frac{1}{q_{2}}}\\
&=2C|Q^{\prime}|^{\frac{1}{p}}\left(\fint_{Q^{\prime}}w_{1}(y)^{-\frac{q_{1}r_{1}}{q_{1}-r_{1}}}dy\right)^{\frac{1}{q_{1}}}\left(\fint_{Q^{\prime}}w_{2}(z)^{-\frac{q_{2}r_{2}}{q_{2}-r_{2}}}dz\right)^{\frac{1}{q_{2}}},
\end{align*}
which yields a contradiction
\begin{align*}
N&<\left(\frac{|Q|}{|Q^{\prime}|}\right)^{\frac{1}{s}}|Q^{\prime}|^{\frac{1}{r}}\left(\fint_{Q}v(x)^{t}dx\right)^{\frac{1}{t}}\\
&\quad\times\left(\fint_{Q^{\prime}}w_{1}(y)^{-\frac{q_{1}r_{1}}{q_{1}-r_{1}}}dy\right)^{\frac{q_{1}-r_{1}}{q_{1}r_{1}}}\left(\fint_{Q^{\prime}}w_{2}(z)^{-\frac{q_{2}r_{2}}{q_{2}-r_{2}}}dz\right)^{\frac{q_{2}-r_{2}}{q_{2}r_{2}}}\leq 2C.
\end{align*}
This produces the desired results. $\hfill$ $\Box$

\section{The proof of Theorem $\ref{main_7}$}

In what follows we always assume that $\eqref{normalization}$ holds. Let $B=(4\cdot 18^{n})^{\frac{1}{r_{1}}+\frac{1}{r_{2}}}$ and
$$\gamma_{0}:=|Q_{0}|^{\frac{\alpha}{n}}\left(\fint_{3Q_{0}}|f(y)|^{r_{1}}dy\right)^{\frac{1}{r_{1}}}\left(\fint_{3Q_{0}}|g(z)|^{r_{2}}dz\right)^{\frac{1}{r_{2}}}.$$
For $k=1,2,\cdots$, we set
$$D_{k}=\bigcup\{Q:\,Q\in\mathcal{D}(Q_{0}), |Q|^{\frac{\alpha}{n}}\left(\fint_{3Q}|f(y)|^{r_{1}}dy\right)^{\frac{1}{r_{1}}}\left(\fint_{3Q}|g(z)|^{r_{2}}dz\right)^{\frac{1}{r_{2}}}>\gamma_{0}B^{k}\}.$$
Considering the maximal cubes with respect to inclusion, we can write, with nonoverlapping cubes
$$D_{k}=\bigcup_{j}Q_{j}^{k}.$$
By the maximality of $Q_{j}^{k}$, we have
\begin{equation}\label{eq:8.1}
\gamma_{0}B^{k}\leq|Q_{j}^{k}|^{\frac{\alpha}{n}}\left(\fint_{3Q_{j}^{k}}|f(y)|^{r_{1}}dy\right)^{\frac{1}{r_{1}}}\left(\fint_{3Q_{j}^{k}}|g(z)|^{r_{2}}dz\right)^{\frac{1}{r_{2}}}\leq 2^{n(\frac{1}{r_{1}}+\frac{1}{r_{2}})}\gamma_{0}B^{k}.
\end{equation}\label{eq:8.2}
Let $E_{0}=Q_{0}\backslash D_{1}$ and $E_{j}^{k}=Q_{j}^{k}\backslash D_{k+1}$ as before, and then we have to check
\begin{equation}
|Q_{0}|\leq 2|E_{0}|\quad\text{and}\quad|Q_{j}^{k}|\leq2|E_{j}^{k}|.
\end{equation}
For fixed $Q_{j}^{k}$, we take
\begin{equation*}
B_{1}:=\left[\left(\int_{3Q_{j}^{k}}|f(y)|^{r_{1}}dy\right)^{\frac{1}{r_{1}}}\left(\int_{3Q_{j}^{k}}|g(z)|^{r_{2}}dz\right)^{\frac{1}{r_{2}}}\right]^{-\frac{r_{2}}{r_{1}+r_{2}}}\left(\frac{\gamma_{0}A^{k+1}}{|Q_{j}^{k}|^{\frac{\alpha}{n}}} \right)^{\frac{r_{2}}{r_{1}+r_{2}}}\left(\int_{3Q_{j}^{k}}|f(y)|^{r_{1}}dy\right)^{\frac{1}{r_{1}}}
\end{equation*}	
and
\begin{equation*}
B_{2}:=\left[\left(\int_{3Q_{j}^{k}}|f(y)|^{r_{1}}dy\right)^{\frac{1}{r_{1}}}\left(\int_{3Q_{j}^{k}}|g(z)|^{r_{2}}dz\right)^{\frac{1}{r_{2}}}\right]^{-\frac{r_{1}}{r_{1}+r_{2}}}\left(\frac{\gamma_{0}A^{k+1}}{|Q_{j}^{k}|^{\frac{\alpha}{n}}} \right)^{\frac{r_{1}}{r_{1}+r_{2}}}\left(\int_{3Q_{j}^{k}}|g(z)|^{r_{2}}dz\right)^{\frac{1}{r_{2}}}.
\end{equation*}	
Observe that $B_{1}B_{2}=\gamma_{0} B^{k+1}$. If $Q_{j^{\prime}}^{k}\subset Q_{j}^{k}$, then
\begin{align}\label{eq:8.3}
\gamma_{0}B^{k+1}&<|Q_{j^{\prime}}^{k+1}|^{\frac{\alpha}{n}}\left(\fint_{3Q_{j^{\prime}}^{k+1}}|f(y)|^{r_{1}}dy\right)^{\frac{1}{r_{1}}}\left(\fint_{3Q_{j^{\prime}}^{k+1}}|g(z)|^{r_{2}}dz\right)^{\frac{1}{r_{2}}}\nonumber\\
&<|Q_{j}^{k}|^{\frac{\alpha}{n}}\left(\fint_{3Q_{j^{\prime}}^{k+1}}|f(y)|^{r_{1}}dy\right)^{\frac{1}{r_{1}}}\left(\fint_{3Q_{j^{\prime}}^{k+1}}|g(z)|^{r_{2}}dz\right)^{\frac{1}{r_{2}}}.
\end{align}
The inequalities $\eqref{eq:8.1}$ and $\eqref{eq:8.3}$ tell us
\begin{align*}
Q_{j}^{k}\cap D_{k+1}&\subset\left\{x\in Q_{j}^{k}:\, M_{\vec{R}}(\chi_{Q_{j}^{k}}f,\chi_{Q_{j}^{k}}g)(x)>\frac{\gamma_{0}B^{k+1}}{|Q_{j}^{k}|^{\frac{\alpha}{n}}}\right\}\\
&\subset\left\{x\in Q_{j}^{k}:\,M_{r_{1}}(\chi_{Q_{j}^{k}}f)(x)M_{r_{2}}(\chi_{Q_{j}^{k}}g)(x)>\frac{\gamma_{0}B^{k+1}}{|Q_{j}^{k}|^{\frac{\alpha}{n}}}\right\}\\
&\subset\left\{x\in\mathbb{R}^{n}:\, M(\chi_{Q_{j}^{k}}f^{r_{1}})(x)>B_{1}^{r_{1}}\right\}\bigcup\left\{x\in\mathbb{R}^{n}:\, M(\chi_{Q_{j}^{k}}g^{r_{2}})(x)>B_{2}^{r_{2}}\right\}
\end{align*}
Using the weak-$(1,1)$ boundedness of $M$, we obtain
\begin{align*}
|Q_{j}^{k}\cap D_{k+1}|&\leq\left|\left\{x\in\mathbb{R}^{n}:\, M(\chi_{Q_{j}^{k}}f^{r_{1}})(x)>B_{1}^{r_{1}}\right\}\right|+\left|\left\{x\in\mathbb{R}^{n}:\, M(\chi_{Q_{j}^{k}}g^{r_{2}})(x)>B_{2}^{r_{2}}\right\}\right|\\
&\leq\frac{3^{n}}{B_{1}^{r_{1}}}\int_{3Q_{j}^{k}}|f(y)|^{r_{1}}dy+\frac{3^{n}}{B_{2}^{r_{2}}}\int_{3Q_{j}^{k}}|g(z)|^{r_{2}}dz\\
&=2\cdot 3^{n}\left[\frac{|Q_{j}^{k}|^{\frac{\alpha}{n}}}{\gamma_{0}B^{k+1}}\left(\int_{3Q_{j}^{k}}|f(y)|^{r_{1}}dy\right)^{\frac{1}{r_{1}}}\left(\int_{3Q_{j}^{k}}|g(z)|^{r_{2}}dz\right)^{\frac{1}{r_{2}}}\right]^{\frac{r_{1}r_{2}}{r_{1}+r_{2}}},
\end{align*}
where we have used the definitions of $B_{1}$ and $B_{2}$. It follows from $\eqref{eq:8.1}$ that
\begin{equation}\label{eq:8.4}
|Q_{j}^{k}\cap D_{k+1}|\leq 2\cdot 3^{n}\left[\frac{|Q_{j}^{k}|^{\frac{\alpha}{n}}}{\gamma_{0}B^{k+1}}\left(\fint_{3Q_{j}^{k}}|f(y)|^{r_{1}}dy\right)^{\frac{1}{r_{1}}}\left(\fint_{3Q_{j}^{k}}|g(z)|^{r_{2}}dz\right)^{\frac{1}{r_{2}}}\right]^{\frac{r_{1}r_{2}}{r_{1}+r_{2}}}|3Q_{j}^{k}|\leq\frac{1}{2}|Q_{j}^{k}|.
\end{equation}
Similary, we can show that
\begin{equation}\label{eq:8.5}
|D_{1}|\leq\frac{1}{2}|Q_{0}|.
\end{equation}
Consequently, we deduce $\eqref{eq:8.2}$ from $\eqref{eq:8.4}$ and $\eqref{eq:8.5}$.

We now return to the proof of Theorem $\ref{main_7}$. Keeping in mind $\eqref{maxinal operator controlled by two parts}$-$\eqref{equivalence formula}$ and $\eqref{notation for weight}$, it follows that
\begin{align*}
&\int_{Q_{0}}\hat{M}_{\alpha,\vec{R}}(f,g)(x)^{t}v(x)^{t}dx\\
&=\int_{E_{0}}\hat{M}_{\alpha,\vec{R}}(f,g)(x)^{t}v(x)^{t}dx+\sum_{k,j}\int_{E_{j}^{k}}\hat{M}_{\alpha,\vec{R}}(f,g)(x)^{t}v(x)^{t}dx\\
&\leq v^{t}(E_{0})(\gamma_{0}B)^{t}+\sum_{k,j}v^{t}(E_{j}^{k})(\gamma_{0}B^{k+1})^{t}\\
&\leq C\left\{|Q_{0}|^{\frac{\alpha}{n}}\left(\fint_{3Q_{0}}|f(y)|^{r_{1}}dy\right)^{\frac{1}{r_{1}}}\left(\fint_{3Q_{0}}|g(z)|^{r_{2}}dz\right)^{\frac{1}{r_{2}}}\left(\fint_{Q_{0}}v(x)^{t}dx\right)^{\frac{1}{t}}\right\}^{t}|Q_{0}|\\
&\quad+C\left\{|Q_{j}^{k}|^{\frac{\alpha}{n}}\left(\fint_{3Q_{j}^{k}}|f(y)|^{r_{1}}dy\right)^{\frac{1}{r_{1}}}\left(\fint_{3Q_{j}^{k}}|g(z)|^{r_{2}}dz\right)^{\frac{1}{r_{2}}}\left(\fint_{Q_{j}^{k}}v(x)^{t}dx\right)^{\frac{1}{t}}\right\}^{t}|Q_{j}^{k}|.
\end{align*}
Using the properties $\eqref{eq:8.2}$, we see further that
\begin{align}\label{eq:8.6}
&\int_{Q_{0}}\hat{M}_{\alpha,\vec{R}}(f,g)(x)^{t}v(x)^{t}dx\nonumber\\
&\leq C\left(\int_{E_{0}}\widetilde{M}_{\alpha,\vec{R}}^{t}(f^{0},g^{0})(x)^{t}dx+\sum_{k,j}\int_{E_{j}^{k}}\widetilde{M}_{\alpha,\vec{R}}^{t}(f^{0},g^{0})(x)^{t}dx\right)\nonumber\\
&\leq C\int_{Q_{0}}\widetilde{M}_{\alpha,\vec{R}}^{t}(f^{0},g^{0})(x)^{t}dx,
\end{align}
where
$$\widetilde{M}_{\alpha,\vec{R}}^{t}(f^{0},g^{0})(x)=\sup_{\mathcal{D}(\mathbb{R}^{n})\ni {Q}\ni x}|{Q}|^{\frac{\alpha}{n}}\left(\fint_{{Q}}|f^{0}(y)|^{r_{1}}dy\right)^{\frac{1}{r_{1}}}\left(\fint_{{Q}}|g^{0}(z)|^{r_{2}}dz\right)^{\frac{1}{r_{2}}}\left(\fint_{Q}v(x)^{t}dx\right)^{\frac{1}{t}}.$$
By condition $\eqref{two weight condition_4}$, we show that
$$c_{0}=\sup_{Q\in\mathcal{D}(\mathbb{R}^{n})}|Q|^{\frac{1}{r}}\left(\fint_{Q}v(x)^{t}dx\right)^{\frac{1}{t}}\prod_{i=1}^{2}\left(\fint_{Q}w_{i}(y)^{-r_{i}\left(\frac{q_{i}}{ar_{i}}\right)^{\prime}}dy\right)^{\frac{1}{r_{i}\left(\frac{q_{i}}{ar_{i}}\right)^{\prime}}}\leq [v,\vec{w}]_{t,\vec{q}/a}^{r,s}.$$
For any $Q\in\mathcal{D}(\mathbb{R}^{n})$, we can calculate that
\begin{align*}
&\left(\fint_{{Q}}|f(y)|^{r_{1}}dy\right)^{\frac{1}{r_{1}}}\left(\fint_{{Q}}|g(z)|^{r_{2}}dz\right)^{\frac{1}{r_{2}}}\left(\fint_{Q}v(x)^{t}dx\right)^{\frac{1}{t}}\\
&\leq\left(\fint_{{Q}}(|f(y)|w_{1}(y))^{\frac{q_{1}}{a}}dy\right)^{\frac{a}{q_{1}}}\left(\fint_{{Q}}(|g(z)|w_{2}(z))^{\frac{q_{2}}{a}}dz\right)^{\frac{a}{q_{2}}}\left(\fint_{Q}v(x)^{t}dx\right)^{\frac{1}{t}}\\
&\quad\times\prod_{i=1}^{2}\left(\fint_{Q}w_{i}(y)^{-r_{i}\left(\frac{q_{i}}{ar_{i}}\right)^{\prime}}dy\right)^{\frac{1}{r_{i}\left(\frac{q_{i}}{ar_{i}}\right)^{\prime}}}\\
&\leq c_{0}|Q|^{-\frac{1}{r}}\left(\fint_{{Q}}(|f(y)|w_{1}(y))^{\frac{q_{1}}{a}}dy\right)^{\frac{a}{q_{1}}}\left(\fint_{{Q}}(|g(z)|w_{2}(z))^{\frac{q_{2}}{a}}dz\right)^{\frac{a}{q_{2}}},
\end{align*}
which implies, for $x\in Q_{0}$,
\begin{equation}\label{eq:8.7}
\widetilde{M}_{\alpha,\vec{R}}^{t}(f^{0},g^{0})(x)\leq C_{0}M_{\alpha-n/r,\vec{q}/a}(f,g)(x).
\end{equation}
The inequalities $\eqref{eq:8.6}$, $\eqref{eq:8.7}$ and Lemma $\ref{lemma_6}$ yield
\begin{equation}\label{eq:8.8}
|Q_{0}|^{\frac{1}{s}}\left(\fint_{Q_{0}}M_{\alpha,\vec{R}}(f^{0},g^{0})(x)^{t}v(x)^{t}dx\right)^{\frac{1}{t}}\leq Cc_{0}\|M_{\alpha-n/r,\vec{q}/a}(f,g)\|_{\mathcal{M}_{t}^{s}}\leq Cc_{0},
\end{equation}
where we have used the assumption
$$\frac{1}{s}=\frac{1}{p}-\frac{\alpha-n/r}{n},\quad\frac{t}{s}=\frac{q}{p}$$
and $\eqref{normalization}$. From H\"{o}lder's inequality, $\eqref{normalization}$ and the fact that
$$\frac{1}{s}=\frac{1}{p}+\frac{1}{r}-\frac{\alpha}{n},$$
it follows that
\begin{align}\label{eq:8.9}
&|Q_{0}|^{\frac{1}{s}}\left(\fint_{Q_{0}}v(x)^{t}dx\right)^{\frac{1}{t}}c_{\infty}\nonumber\\
&\leq|Q_{0}|^{\frac{1}{s}}\left(\fint_{Q_{0}}v(x)^{t}dx\right)^{\frac{1}{t}}\mathop{\sup_{Q\in\mathcal{D}(\mathbb{R}^{n})}}_{Q\supset Q_{0}}|Q|^{\frac{\alpha}{n}}\left(\fint_{Q}(|f(y)|w_{1}(y))^{\frac{q_{1}}{a}}dy\right)^{\frac{a}{q_{1}}}\nonumber\\
&\qquad\times\left(\fint_{Q}(|g(z)|w_{2}(z))^{\frac{q_{2}}{a}}dz\right)^{\frac{a}{q_{2}}}\prod_{i=1}^{2}\left(\fint_{Q}w_{i}(y)^{-r_{i}\left(\frac{q_{i}}{ar_{i}}\right)^{\prime}}dy\right)^{\frac{1}{r_{i}\left(\frac{q_{i}}{ar_{i}}\right)^{\prime}}}\nonumber\\
&\leq\mathop{\sup_{Q\in\mathcal{D}(\mathbb{R}^{n})}}_{Q\supset Q_{0}}\left(\frac{|Q_{0}|}{|Q|}\right)^{\frac{1}{s}}|Q|^{\frac{1}{r}+\frac{1}{p}}\left(\fint_{Q_{0}}v(x)^{t}dx\right)^{\frac{1}{t}}\prod_{i=1}^{2}\left(\fint_{Q}w_{i}(y)^{-r_{i}\left(\frac{q_{i}}{ar_{i}}\right)^{\prime}}dy\right)^{\frac{1}{r_{i}\left(\frac{q_{i}}{ar_{i}}\right)^{\prime}}}\nonumber\\
&\qquad\times\left(\fint_{Q}(|f(y)|w_{1}(y))^{{q_{1}}}dy\right)^{\frac{1}{q_{1}}}\left(\fint_{Q}(|g(z)|w_{2}(z))^{{q_{2}}}dz\right)^{\frac{1}{q_{2}}}\nonumber\\
&\leq\mathop{\sup_{Q\in\mathcal{D}(\mathbb{R}^{n})}}_{Q\supset Q_{0}}\left(\frac{|Q_{0}|}{|Q|}\right)^{\frac{1}{s}}|Q|^{\frac{1}{r}}\left(\fint_{Q_{0}}v(x)^{t}dx\right)^{\frac{1}{t}}\prod_{i=1}^{2}\left(\fint_{Q}w_{i}(y)^{-r_{i}\left(\frac{q_{i}}{ar_{i}}\right)^{\prime}}dy\right)^{\frac{1}{r_{i}\left(\frac{q_{i}}{ar_{i}}\right)^{\prime}}}\leq[v,\vec{w}]_{t,\vec{q}/a}^{r,s}.
\end{align}
Combining $\eqref{eq:8.8}$ together with $\eqref{eq:8.9}$, we obtain the desired result. $\hfill$ $\Box$

\section{The proof of Theorem $\ref{main_8}$}

Using the condition $t<s$ and H\"{o}lder's inequality, the weight condition $\eqref{two weight condition_4}$ holds only if
\begin{equation}\label{eq:9.1}
\sup_{Q\in\mathcal{D}(\mathbb{R}^{n})}|Q|^{\frac{1}{r}}\left(\fint_{Q}v^{s}\right)^{\frac{1}{s}}\prod_{i=1}^{2}\left(\fint_{Q}w_{i}^{-r_{i}\left(\frac{q_{i}}{ar_{i}}\right)^{\prime}}\right)^{\frac{1}{r_{i}\left(\frac{q_{i}}{ar_{i}}\right)^{\prime}}}<\infty
\end{equation}
holds. Taking $v=u_{1}^{\frac{1}{q_{1}}}u_{2}^{\frac{1}{q_{2}}}$ and  $r=\infty$,  the condition $\eqref{eq:9.1}$ becomes as
\begin{equation}\label{eq:3.118}
\sup_{Q\in\mathcal{D}(\mathbb{R}^{n})}\left(\fint_{Q}u_{1}^{\frac{s}{q_{1}}}u_{2}^{\frac{s}{q_{2}}}\right)^{\frac{1}{s}}\prod_{i=1}^{2}\left(\fint_{Q}u_{i}^{-\frac{r_{i}}{q_{i}}\left(\frac{q_{i}}{ar_{i}}\right)^{\prime}}\right)^{\frac{1}{r_{i}\left(\frac{q_{i}}{ar_{i}}\right)^{\prime}}}<\infty.
\end{equation}
Now, we assume that the weight condition $\eqref{two weight condition_5}$ holds. By Lemma $\ref{Characterization of a multiple weights}$, we give the following relationship
$$u_{1}^{\frac{s}{q_{1}}}u_{2}^{\frac{s}{q_{2}}}\in A_{2s},\quad u_{1}^{-\frac{r_{1}}{q_{1}-r_{1}}}\in A_{\frac{2q_{1}r_{1}}{q_{1}-r_{1}}}\quad\text{and}\quad u_{2}^{-\frac{r_{2}}{q_{2}-r_{2}}}\in A_{\frac{2q_{2}r_{2}}{q_{2}-r_{2}}}.$$
By means of Lemma $\ref{lemma_7}$, there exists $a>1$ such that
$$
\left(\fint_{Q}u_{1}^{-\frac{r_{1}}{q_{1}}\left(\frac{q_{1}}{ar_{1}}\right)^{\prime}}\right)^{\frac{1}{r_{1}\left(\frac{q_{1}}{ar_{1}}\right)^{\prime}}}\leq \left(\fint_{Q}u_{1}^{-\frac{r_{1}}{q_{1}-r_{1}}}\right)^{\frac{q_{1}-r_{1}}{r_{1}q_{1}}}
$$
and
$$
\left(\fint_{Q}u_{2}^{-\frac{r_{2}}{q_{2}}\left(\frac{q_{2}}{ar_{2}}\right)^{\prime}}\right)^{\frac{1}{r_{2}\left(\frac{q_{2}}{ar_{2}}\right)^{\prime}}}\leq\left(\fint_{Q}u_{2}^{-\frac{r_{2}}{q_{2}-r_{2}}}\right)^{\frac{q_{2}-r_{1}}{r_{2}q_{2}}}.
$$
This implies that $\eqref{eq:3.118}$ holds. In other words, the weight condition $\eqref{two weight condition_4}$ holds. Hence, it is direct to imply Theorem $\ref{main_8}$ from Theorem $\ref{main_7}$.
\section{Applications and examples}
\noindent 7.1. {\bf{A bilinear Stein-Weiss inequality}}. Given $0<\alpha<n$, let $T_{\alpha}$ be define by
$$T_{\alpha}f(x)=I_{n-\alpha}f(x)=\int_{\mathbb{R}^{n}}\frac{f(y)}{|x-y|^{n-\alpha}}dy.$$
Stein and Weiss\cite{SW1958} proved the following weighted inequality for $T_{\alpha}$:
$$\left(\int_{\mathbb{R}^{n}}\left(\frac{T_{\alpha}f(x)}{|x|^{\beta}}\right)^{t}dx\right)^{1/t}\leq\left(\int_{\mathbb{R}^{n}}(f(x)|x|^{\gamma})^{q}\right)^{1/q},$$
where $\alpha$, $\beta$, $\gamma$ are positive numbers that depend on $p$ and $q$.
Below, we consider a bilinear Stein-Weiss inequality on Morrey spaces for case $1<t\leq s<\infty$. For $0<\alpha<n$, let $BT_{\alpha}$ be the bilinear operator defined by
$$BT_{\alpha}(f,g)(x)=B_{n-\alpha}(f,g)(x)=\int_{\mathbb{R}^{n}}\frac{f(x-y)g(x+y)}{|y|^{\alpha}}dy.$$
The case for $0<t\leq s\leq1$ was obtained by the authors \cite{HY2018}.

\begin{theorem}\label{bilinear Stein-Weiss inequality on Morrey spaces}
	Assume that $1<q_{1}\leq p_{1}<\infty$, $1<q_{2}\leq p_{2}<\infty$ and $1<t\leq s<\infty$, $\frac{n}{n-\alpha}<r\leq\infty$. Here, $p$ and $q$ are given by
	$$\frac{1}{p}=\frac{1}{p_{1}}+\frac{1}{p_{2}}\quad\text{and}\quad\frac{1}{q}=\frac{1}{q_{1}}+\frac{1}{q_{2}}.$$
	Suppose that
	$$\frac{1}{s}=\frac{1}{p}+\frac{1}{r}-\frac{n-\alpha}{n},\quad\frac{1}{t}=\frac{1}{q}+\frac{1}{r}-\frac{n-\alpha}{n}$$
	and $\alpha$, $\beta$, $\gamma_{1}$, $\gamma_{2}$ satisfy the conditons
	\begin{equation}\label{bilinear Stein-Weiss inequality condition}
	\left\{
	\begin{aligned}
	&\beta<\frac{n}{s},\qquad \gamma_{1}<\frac{n}{q_{1}^{\prime}},\qquad \gamma_{2}<\frac{n}{q_{2}^{\prime}}, \\
	&\quad\alpha+\beta+\gamma_{1}+\gamma_{2}=n+\frac{n}{t}-\frac{n}{q},\\
	&\qquad\qquad\quad\beta+\gamma_{1}+\gamma_{2}\geq0.
	\end{aligned}
	\right.
	\end{equation}
	Then the following inequality holds for all $f,g\geq 0$,
	\begin{align*}
	&\sup_{Q\in\mathcal{D}(\mathbb{R}^{n})}|Q|^{\frac{1}{s}}\left(\fint_{Q}\left(\frac{BT_{\alpha}(f,g)(x)}{|x|^{\beta}}\right)^{t}\right)^{\frac{1}{t}}\\
	&\qquad\qquad\leq C\sup_{Q\in\mathcal{D}(\mathbb{R}^{n})}|Q|^{\frac{1}{p_{1}}}\left(\fint_{Q}\left(f(x)|x|^{\gamma_{1}}\right)^{q_{1}}\right)^{\frac{1}{q_{1}}}\sup_{Q\in\mathcal{D}(\mathbb{R}^{n})}|Q|^{\frac{1}{p_{2}}}\left(\fint_{Q}\left(g(x)|x|^{\gamma_{2}}\right)^{q_{2}}\right)^{\frac{1}{q_{2}}}.
	\end{align*}
\end{theorem}
\begin{remark}
	The first condition of $\eqref{bilinear Stein-Weiss inequality condition}$ corresponds to the condition of Theorem 6.1
	$$\beta<(1-s)\frac{n}{s},\quad \gamma_{1}<\frac{n}{q_{1}^{\prime}},\quad \gamma_{2}<\frac{n}{q_{2}^{\prime}}$$
	as in \cite{HY2018}. The interesting phenomenon here is that the factor $1-s$ (for the case $t\leq s\leq1$) has become $1$ (for the case $1<t\leq s$). This may reveal some clues about the case $t\leq 1<s$.
\end{remark}

\begin{remark}
	If we consider the bilinear case on Lebesgue space, Hoang and Moen \cite{HM2016} obtain that the conditions $\eqref{bilinear Stein-Weiss inequality condition}$ can be reduced to
	\begin{equation*}
	\left\{
	\begin{aligned}
	&\beta<\frac{n}{s},\qquad \gamma_{1}<(q-1)\frac{n}{q_{1}},\qquad \gamma_{2}<(q-1)\frac{n}{q_{2}}, \\
	&\quad\alpha+\beta+\gamma_{1}+\gamma_{2}=n+\frac{n}{t}-\frac{n}{q},\\
	&\qquad\qquad\quad\beta+\gamma_{1}+\gamma_{2}\geq0.
	\end{aligned}
	\right.
	\end{equation*}
	For the linear case, Stein and Weiss \cite{SW1958} obtained the conditions $\eqref{bilinear Stein-Weiss inequality condition}$ which can just drop $q_{2}$ and $\gamma_{2}$.
\end{remark}

{\noindent}$Proof$ $of$ $Theorem$ $\ref{bilinear Stein-Weiss inequality on Morrey spaces}$. Taking
$$1<t_{0}\leq s<\infty\quad\text{and}\quad \frac{t_{0}}{s}=\frac{q_{1}}{p_{1}}=\frac{q_{2}}{p_{2}}.$$
Similar to the estimate for Theorem 1.4 in \cite{HY2018}, we have
\begin{equation}\label{eq:9.2}
\sup_{Q\in\mathcal{D}(\mathbb{R}^{n})}|Q|^{\frac{1}{s}}\left(\fint_{Q}\left(\frac{BT_{\alpha}(f,g)(x)}{|x|^{\beta}}\right)^{t}dx\right)^{\frac{1}{t}}\leq\sup_{Q\in\mathcal{D}(\mathbb{R}^{n})}|Q|^{\frac{1}{s}}\left(\fint_{Q}\left(\frac{BT_{\alpha}(f,g)(x)}{|x|^{\beta}}\right)^{t_{0}}dx\right)^{\frac{1}{t_{0}}}.
\end{equation}
Using Theorem $\ref{main_1}$, H\"{o}lder's inequality, $\eqref{two weight condition_2}$ and $\eqref{eq:9.2}$, we only need to prove that there exists $1<a<\min(q_{1},q_{2})$ such that
\begin{equation}\label{eq:9.3}
\sup_{Q\in\mathcal{D}(\mathbb{R}^{n})}|Q|^{\frac{1}{r}}\left(\fint_{Q}|x|^{-a\beta s}dx\right)^{\frac{1}{as}}\prod_{i=1}^{2}\left(\fint_{Q}|x|^{-\gamma_{i}(q_{i}/a)^{\prime}}dx\right)^{\frac{1}{(q_{i}/a)^{\prime}}}<\infty.
\end{equation}
For any cube $Q$, we set $Q_{0}=Q(0,\ell(Q))$, and then there holds either of the two cases: $Q\cap Q_{0}=\emptyset$ and $Q\cap Q_{0}\neq\emptyset$.

First we consider the case for $Q\cap Q_{0}=\emptyset$. By geometry property, we have $|x|\sim |x|_{\infty}\geq \ell(Q)$ for all $x\in Q$, which implies that the left hand side of $\eqref{eq:9.3}$ can be estimated as
$$\sup_{Q\in\mathcal{D}(\mathbb{R}^{n})}|\ell(Q)|^{n-\alpha+\frac{n}{t}-\frac{n}{q}}|\ell(Q)|^{-\beta-\gamma_{1}-\gamma_{2}}=1.$$

Next, we consider the case for $Q\cap Q_{0}\neq\emptyset$. Let $1<a\min(q_{1},q_{2})$ be such that $\beta<1/as$, $\gamma_{1}<1/(q_{1}/a)^{\prime}$ and $\gamma_{2}<1/(q_{2}/a)^{\prime}$.  If $Q\cap Q_{0}\neq\emptyset$, by geometry property, we have $|x|\leq\sqrt{n}|x|_{\infty}\leq 2\sqrt{n}\ell(Q)$ for all $x\in Q$. This implies $Q\subset B=B(0,2\sqrt{n}\ell(Q))$. Hence, the left hand side of $\eqref{eq:9.3}$ can be estimated as
\begin{align*}
&\sup_{Q\in\mathcal{D}(\mathbb{R}^{n})}|Q|^{\frac{1}{r}}\left(\fint_{Q}|x|^{-a\beta s}dx\right)^{\frac{1}{as}}\prod_{i=1}^{2}\left(\fint_{Q}|x|^{-\gamma_{i}(q_{i}/a)^{\prime}}dx\right)^{\frac{1}{(q_{i}/a)^{\prime}}}\\
&\leq C\sup_{Q\in\mathcal{D}(\mathbb{R}^{n})} |\ell(Q)|^{\beta+\gamma_{1}+\gamma_{2}}\left(\fint_{B}|x|^{-a\beta s}dx\right)^{\frac{1}{as}}\prod_{i=1}^{2}\left(\fint_{B}|x|^{-\gamma_{i}(q_{i}/a)^{\prime}}dx\right)^{\frac{1}{(q_{i}/a)^{\prime}}}\\
&\leq C.
\end{align*}
This completes the proof of Theorem $\ref{bilinear Stein-Weiss inequality on Morrey spaces}$. $\hfill$ $\Box$

\bigskip\bigskip

{\noindent}{\bf Acknowledgement}. This work was supported by the National Natural Science Foundation of China (11561062 and 11871452) and Project of Henan Provincial Department of Education(18A110028).

 \end{document}